%% file: 03112026.tex
\documentclass{article}

\usepackage{graphicx} 

\usepackage[a4paper, total={6in, 8in}]{geometry}
\usepackage{fancyhdr}
\usepackage{pdfpages}
\usepackage{hyperref}
\usepackage{caption}
\hypersetup{colorlinks = true, linkcolor = blue!70!black, citecolor= blue!70!black}
\usepackage[toc,page]{appendix}
\usepackage{dsfont}
\usepackage[utf8]{inputenc}
\usepackage{bbm}
\usepackage{amsmath}
\usepackage[normalem]{ulem}
\usepackage{amsfonts}
\usepackage{amssymb}
\usepackage{amsthm}
\usepackage{graphicx} 
\usepackage{subcaption}
\usepackage[toc,page]{appendix}
\newtheorem*{theorem*}{Theorem}
\newtheorem{theorem}{Theorem}[section]
\newtheorem{corollary}[theorem]{Corollary}
\newtheorem{lemma}[theorem]{Lemma}
\newtheorem{definition}[theorem]{Definition}
\newtheorem{remark}[theorem]{Remark}
\newtheorem*{remark*}{Remark}
\newtheorem{proposition}[theorem]{Proposition}
\usepackage{enumitem}
\usepackage[section]{placeins}
\usepackage{float}
\usepackage[toc,page]{appendix}
\numberwithin{equation}{section}

\usepackage{standalone}
\usepackage{tikzit}
\usepackage[utf8]{inputenc}
\usepackage{pgfplots}
\pgfplotsset{width=7cm, compat=1.10}
\usepgfplotslibrary{fillbetween}
\pgfplotsset{compat=1.18}
\usetikzlibrary{patterns}

\usepackage[utf8]{inputenc}
\usepackage{newunicodechar}
\newunicodechar{̀}{\`{}}

\DeclareUnicodeCharacter{2217}{*}

\usepackage[
  backend=biber,
  style=numeric,
  sortcites=true,
  defernumbers=true,
  giveninits=true,   
  eprint=false,
  url=false,         
  doi=false,         
  isbn=false         
]{biblatex}
\renewcommand{\epsilon}{\varepsilon}

\usepackage{todonotes}

\title{Minimal time of null-controllability of the multi-dimensional singular Grushin equation}
\author{Roman Vanlaere\footnote{Dep. EDAN and IMUS, Universidad de Sevilla, 41080 Sevilla, Spain. Email: rvanlaere@us.es}}

\date{October 2026}

\begin{document}

\maketitle

\begin{abstract}
    In this paper, we study the null-controllability of the multi-dimensional Grushin equation with an inverse-square potential on cylinders. In the critically degenerate case, we compute the minimal time of null-controllability from vertical strips, obtaining an explicit formula depending on the distance of the strip to the singular set and on the strength of the singular potential. To do so, after a Fourier and spherical harmonic decomposition of our control system, we use the moments method for radial problems to construct controls on annuli. For this purpose, we prove a uniform spectral gap and uniform localization estimates for the eigenfunctions of the infinite family of radial operators that arise. We then reduce the control zone using spectral inequalities and the Lebeau-Robbiano method. In the weakly degenerate case, we prove null-controllability in arbitrarily small time from any open set. This is done using the transmutation method and quantitative observability estimates for radial waves, which yield a cost of null-controllability from annuli for the Fourier components of our system that is uniform with respect to the frequency. We conclude again with spectral inequalities and the Lebeau-Robbiano method.
\end{abstract}

{\small

\medskip

\par}

\tableofcontents

\section{Introduction}

The goal of this paper is to compute the minimal time of null-controllability of 
\begin{align}\label{eqn: control equation intro}
    \partial_t f - \Delta_xf - |x|^{2\gamma} \Delta_y f + \displaystyle\frac{\nu^2-\mathsf{H}}{|x|^2}f  =  \mathbbm{1}_\omega u(t,x,y), \quad (t,x,y) \in (0,T) \times \Omega_x \times \Omega_y.
\end{align}

It was first proved in \cite{cannarsa2014null,anh2016null}, for $\gamma = 1$, $\nu > 0$, and $\omega = \omega_x \times \Omega_y$, that a positive minimal time of null-controllability exists. We went further in \cite{vanlaere2026observability} by giving an explicit characterization of this minimal time under some geometric restriction on $\omega_x$, which is automatically satisfied when $\Omega_x = (0,L)$. This characterization highlights the dependence of the minimal time on both the observation set and the strength of the singularity, and reveals that the minimal time for the Laplace-Beltrami operator arising in almost-Riemannian geometry depends on the dimension of the singular set, an observation that we strengthen in the present paper. We also showed in \cite{vanlaere2026observability} that if $\gamma > 1$ and $\nu > 0$, then \eqref{eqn: control equation intro} is never null-controllable, and we left the case $\gamma \in (0,1)$ open.\\

We go even further in the present paper, in the case $\Omega_x = B(0,1)$, by 
\begin{enumerate}[label= (\alph*)]
    \item obtaining an explicit characterization of the minimal time of null-controllability when $\gamma = 1$ without the geometric restriction imposed on $\omega_x$ in \cite{vanlaere2026observability},
    \item showing that null-controllability holds in arbitrarily small times from any open set $\omega \subset \Omega$ when $\gamma \in (0,1)$, 
    \item treating the critical Hardy case $\nu = 0$ that was left open in \cite{anh2016null,vanlaere2026observability, cannarsa2014null},
    \item obtaining uniform spectral estimates for the Fourier and spherical decompositions of the spatial operator appearing in \eqref{eqn: control equation intro} when $\gamma = 1$.
\end{enumerate}
Our approach when $\gamma = 1$ relies on the moments method, thus on a very careful spectral analysis, and a Lebeau-Robbiano strategy as used in \cite{DardeTrabut2026}. When $\gamma \in (0,1)$, our approach is based on the transmutation method, similarly to \cite{beauchard20152d}, along with a concluding application of the Lebeau-Robbiano method similarly to \cite[Theorem 1 (1)]{beauchard2014null}. \\

Our analysis also provides an alternative approach to the computation of the minimal time of null-controllability obtained in \cite{DardeTrabut2026} in the non-singular case (\textit{i.e.} $\nu^2 = \mathsf{H}$) for $\gamma = 1$, and extends their results to $\gamma \in (0,1)$. When $\gamma = 1$, we employ the moments method because we were not able to extend the Carleman machinery of \cite{DardeTrabut2026} to our setting. The Carleman route to obtain the results of this paper nevertheless remains of interest and open.

\subsection{Setting and main results}

Let $\Omega = B(0,1) \times \Omega_y$, with $B(0,1) := \{x \in \mathbb{R}^{d_x}, \ |x| < 1\}$, $d_x \geq 3$\footnote{See \cite[Remark 1.1]{vanlaere2026observability} for the case $d_x = 2$.}, and $\Omega_y \cup \partial \Omega_y$ a connected compact Riemannian manifold of dimension $d_y \geq 1$ possibly with boundary. We consider the following control system. Let $T > 0$, 
\begin{align} \label{stm : control system}\tag{$\mathrm{H}$} 
\left\{ \begin{array}{lcll}
     \partial_t f - \Delta_xf - |x|^{2\gamma} \Delta_y f + \displaystyle\frac{\nu^2-\mathsf{H}}{|x|^2}f & = & \mathbbm{1}_\omega u(t,x,y), & (t,x,y) \in (0,T) \times \Omega,\\[6pt]
     f(t,x,y) & = & 0, & (t,x,y) \in (0,T) \times \partial \Omega, \\[6pt]
     f(0,x,y) & = & f_0(x,y), & (x,y) \in \Omega,
    \end{array} \right.
\end{align}
where $f_0 \in L^2(\Omega)$, $f$ is the state, $u \in L^2((0,T) \times \Omega)$ is the control supported in a nonempty open set $\omega \subset \Omega$. Moreover, $\gamma > 0$ and $\nu \geq 0$ are fixed parameters (the case $\gamma = 0$ being now understood, see \cite{ervedoza2008control}), and 
\begin{align}\label{eqn: definition hardy constant}
     \mathsf{H}  = (d_x - 2)^2/4
\end{align}
is the largest constant such that (see \textit{e.g.} \cite{hardy1952inequalities,cazacu2014controllability}),
\begin{align}\label{eqn: hardy inequality}
    \mathsf{H}\int_{B(0,1)} \frac{|f|^2}{|x|^2} \ dx \leq \int_{B(0,1)}  |\nabla_x f|^2 \ dx, \quad \text{for every } f \in H_0^1(B(0,1)).
\end{align}
This Hardy inequality also holds for $d_x = 1$ on $H_0^1(0,1)$ with $\mathsf{H} = 1/4$. \\

The computation of the minimal time of null-controllability (see Definition \ref{def: null-controllability} below) for system \eqref{stm : control system} when $\gamma = 1$ is motivated by the following. Let $\sigma$ be a Riemannian metric on $\Omega_y$, and let $(E_j)_{1 \leq j \leq d_y}$ be a local orthonormal frame for $\sigma$. A model sub-Riemannian structure on $\Omega$ is generated by the vector fields $X_i = \partial_{x_i}$, and $Y_{ij} = x_i E_j$, $1 \leq i \leq d_x$, $1 \leq j \leq d_y$. Outside $\{0\} \times \Omega_y$, the corresponding Riemannian metric is given by
\begin{align*}
    \mathsf{g} = dx^2 + |x|^{-2}\sigma,
\end{align*}
which is singular at $\{x = 0\}$. The corresponding Laplace-Beltrami operator is given by
\begin{align}\label{eqn: def laplace beltrami}
    -\Delta_{\mathsf{g}} f = -\Delta_x f - |x|^2\Delta_y f
    + \frac{d_y}{|x|^2}x \cdot \nabla_x f,
\end{align}
on $L^2(\Omega, |x|^{-d_y}dx\,d\mathrm{vol}_\sigma)$, with $\Delta_y$ the Laplace-Beltrami operator associated with $\sigma$ on $\Omega_y$, and $d\mathrm{vol}_\sigma$ the corresponding Riemannian volume form. The change of variable $f = |x|^{d_y/2}g$ leads to
\begin{align*}
    L g = -\Delta_x g - |x|^2\Delta_y g
    + \frac{d_y}{2}\left(\frac{d_y}{2} - d_x + 2\right)
    \frac{g}{|x|^2},
\end{align*}
on $L^2(\Omega, dxd\mathrm{vol}_\sigma)$. Since $\mathsf{H} = (d_x-2)^2/4$, we have
\begin{align*}
    \frac{d_y}{2}\left(\frac{d_y}{2} - d_x + 2\right)
    = \frac{(d_y-d_x+2)^2}{4} - \mathsf{H} \geq - \mathsf{H}.
\end{align*}
In system \eqref{stm : control system}, the operator $L$ therefore corresponds to $\nu = |d_y-d_x+2|/2$. For the corresponding Friedrichs realizations, the operator $-\Delta_x - |x|^2\Delta_y + (\nu^2-\mathsf{H})/|x|^2$ is thus a generalization of $L$ in which the degenerate term and the singular term are decoupled. This allows us to treat a broader class of operators of the same type as $L$.\\

The study of the null-controllability properties of subelliptic heat equations (and of the corresponding wave and Schrödinger equations) is a long-standing line of research, in which this work is situated. We review the state of the art in Section \ref{section: state of the art}. 

\begin{definition}[Null-controllability]\label{def: null-controllability}
Let $T>0$ and $\omega \subset \Omega$ be a nonempty open set. We say that system \eqref{stm : control system} is null-controllable from $\omega$ in time $T > 0$ if, for every $f_0 \in L^2(\Omega)$, there exists $u \in L^2((0,T) \times \Omega)$ such that the associated solution satisfies $f(T) = 0$.
\end{definition}

For any $\omega \subset \Omega$, we define
\begin{align}
    T(\omega) := \inf \{T > 0, \text{ such that the system is null-controllable in time $T$ from $\omega$}\}.
\end{align}

Our first main result is the following.

\begin{theorem}\label{thm: main thm minimal time gamma = 1}
    Let $\gamma = 1$. Let $\nu \geq 0$, and $\omega_x \subset B(0,1)$ be a nonempty open set. The minimal time of null-controllability of system \eqref{stm : control system} from $\omega_x \times \Omega_y$ satisfies
    \begin{align}\label{eqn: def minimal time null-control gamma = 1}
        T(\omega_x \times \Omega_y) = \frac{\operatorname{dist} \left(\omega_x, 0_{\mathbb{R}^{d_x}} \right)^2}{4(1+\nu)}.
    \end{align}
\end{theorem}

Theorem \ref{thm: main thm minimal time gamma = 1} is proved in Sections \ref{section: Proof of Theorem 1 positive} and \ref{section: Proof of Theorem 1 negative}, and its strategy of proof is given in Section \ref{section: Strategy}. Proving \eqref{eqn: def minimal time null-control gamma = 1} requires two main steps. First, we prove that $T(\omega_x \times \Omega_y) \leq \operatorname{dist} \left(\omega_x, 0_{\mathbb{R}^{d_x}} \right)^2/(4(1+\nu)) $, \textit{i.e.}, that for any $T > \operatorname{dist} \left(\omega_x, 0_{\mathbb{R}^{d_x}} \right)^2/(4(1+\nu))$, system \eqref{stm : control system} is null-controllable (Section \ref{section: Proof of Theorem 1 positive}). Second, we prove that $T(\omega_x \times \Omega_y) \geq \operatorname{dist} \left(\omega_x, 0_{\mathbb{R}^{d_x}} \right)^2/(4(1+\nu))$, \textit{i.e.}, that for any $T < \operatorname{dist} \left(\omega_x, 0_{\mathbb{R}^{d_x}} \right)^2/(4(1+\nu))$, system \eqref{stm : control system} is not null-controllable (Section \ref{section: Proof of Theorem 1 negative}). \\

As a corollary of Theorem \ref{thm: main thm minimal time gamma = 1}, for the Friedrichs realization of the Laplace-Beltrami operator $\Delta_{\mathsf{g}}$ given in \eqref{eqn: def laplace beltrami}, the minimal time of null-controllability is given by
\begin{align}
    T(\omega_x \times \Omega_y) = \frac{\operatorname{dist} \left(\omega_x, 0_{\mathbb{R}^{d_x}} \right)^2}{4 + 2|d_y-d_x+2|}.
\end{align}

When $\nu^2 = \mathsf{H}$, we recover the sub-Laplacian operator for the proposed sub-Riemannian structure associated with the product measure $dx d\mathrm{vol}_\sigma$, induced by the metric $\tilde{\mathsf{g}} = dx^2 + \sigma$. The corresponding minimal time is
\begin{align}
    T(\omega_x \times \Omega_y) = \frac{\operatorname{dist} \left(\omega_x, 0_{\mathbb{R}^{d_x}} \right)^2}{2d_x},
\end{align}
which was obtained in \cite{DardeTrabut2026}.\\

Our second main result is the following. 

\begin{theorem}\label{thm: main thm minimal time gamma < 1}
    Let $\gamma \in (0,1)$. Let $\nu \geq 0$, and $\omega \subset B(0,1) \times \Omega_y$ a nonempty open set. The minimal time of null-controllability of system \eqref{stm : control system} from $\omega$ satisfies
    \begin{align}
        T(\omega) = 0.
    \end{align}
\end{theorem}

We obtained in \cite[Theorem 1.5]{vanlaere2026observability} that if $\gamma > 1$ and $\nu > 0$, then, for any open set $\omega \subset \Omega$ such that $\overline{\omega} \cap \{0\} \times \Omega_y = \emptyset$, system \eqref{stm : control system} is never null-controllable. This result can be extended to $\nu = 0$ by adapting \cite[Proposition 5.1]{vanlaere2026observability} to this case. Thus, the present paper completes the picture, recovering the three thresholds of null-controllability with respect to $\gamma > 0$ first obtained in the non-singular case on $(-1,1) \times (0,1)$ in the pioneering work \cite{beauchard2014null}. 

\begin{remark}
    The methods employed in the present work extend naturally, and in a much simpler manner, to the setting of $d_x = 1$ and $\Omega = (0,1) \times \Omega_y$. In particular, Theorems \ref{thm: main thm minimal time gamma = 1} and \ref{thm: main thm minimal time gamma < 1} are true in this one-dimensional setting.  
\end{remark}

\subsection{State of the art}\label{section: state of the art}

The present work, like those reviewed below, belongs to a line of research that aims at understanding the null-controllability properties of subelliptic heat equations.\\

The study of the null-controllability properties of the singular Grushin equation has not attracted as much attention in the literature as the non-singular case. Approximate controllability has been obtained in \cite{morancey2015approximate} when $d_x = 1$, with an interior singularity. When $\gamma = 1$, existence of a minimal time of null-controllability has been proved when $d_x = 1$ with a boundary singularity in \cite{cannarsa2014null}, and the case $d_x\geq 3$ is treated in \cite{anh2016null}. We complemented these two works in \cite{vanlaere2026observability}, by providing a characterization of this minimal time under some geometric restriction on $\omega_x$ when $d_x \geq 3$, a restriction that is automatically satisfied when $d_x = 1$, and showing that null-controllability fails in any time when $\gamma > 1$. In \cite{morancey2015approximate,cannarsa2014null,anh2016null, vanlaere2026observability}, it is always assumed that $\nu > 0$. Theorem \ref{thm: main thm minimal time gamma = 1} thus completes these results, and as emphasized, also extends \cite{DardeTrabut2026}. \\

When the singular term vanishes, we recover the classical non-singular Grushin equation which now enjoys a vast literature. The operator at hand below is $-\Delta_x - |x|^{2\gamma}\Delta_y$. The pioneering work on the subject is \cite{beauchard2014null}, in the two-dimensional setting $\Omega = (-1,1) \times (0,1)$, where the control zone is a vertical strip (\textit{i.e.} $\omega  =\omega_x \times \Omega_y$). It is shown that the minimal time depends on the strength of the degeneracy. If $\gamma < 1$, null-controllability occurs in any time, if $\gamma > 1$ it is never null-controllable, while in the critical case $\gamma = 1$, a minimal amount of time is required. The minimal time in the two-dimensional setting is obtained, when the domain and the control zone are symmetric with respect to $x=0$, in \cite{beauchard20152d} using the transmutation method. This result is extended in \cite{allonsius2021analysis} using the moments method for control zones consisting of a single vertical strip. The equation was then generalized, by replacing $x$ with $q(x)$, and when the control is any vertical strip at the boundary of the domain, the minimal time is obtained in \cite{beauchard2020minimal}. In addition to these results, it has been shown in \cite{LISSYprolate2025, duprez2020control,darde2023null,CRMATH_2017__355_12_1215_0} that a geometric condition on $\omega$ is necessary when $\gamma = 1$. If $\omega$ is the complement of a horizontal strip, then null-controllability always fails. The analysis in higher dimensions was first provided in \cite{beauchard2014inverse}, where the authors only proved existence of a minimal time, which has been computed very recently in \cite{DardeTrabut2026}. We mention that the work \cite{beauchard2020minimal} also treats the multi-dimensional case with a control supported on the whole vertical boundary $\partial\Omega_x \times \Omega_y$. We also refer the reader to \cite{huang2026observability,alphonse2024quantitative,alphonse2025unique,jaming2025null} in non-compact settings. \\

We also have results in more general settings. In \cite{tamekue2022null}, the heat equation is studied on the Grushin sphere, and in this setting the minimal time is obtained in \cite{tamekue2026minimal}. On almost-Riemannian surfaces, obstruction to null-controllability in small times under a local assumption of separation of variables was obtained in \cite{vanlaere2025grushinlike}. Removing this assumption, the minimal time of null-controllability on almost-Riemannian surfaces has been obtained recently in \cite{boulard2026minimal}.  \\

Concerning other null-controllability results for a class of degenerate parabolic equations, we can also cite \cite{beauchard2017heat}, which studies the heat equation on the Heisenberg group, or results concerning the Kolmogorov equation \cite{LEROUSSEAU20163193, beauchardKolmog2015, koenigfractionalheat, beauchardkolmogorov, DardeRoyerKolmogorov}. We also refer to \cite{laurent2022tunneling} for approximate controllability of hypoelliptic
equations on sub-Riemannian manifolds. \\

For subelliptic wave equations, we refer to \cite{letrouit2023subelliptic}. For the minimal time of null-controllability of the Grushin Schrödinger equation, we refer to \cite{burq2022time}, and on almost-Riemannian surfaces to \cite{boulard2026minimal}. \\

Concerning the heat equation in the presence of a singular potential in the interior without degenerate term ($\gamma = 0$ in \eqref{stm : control system}), we refer to \cite{vancostenoble2008null,ervedoza2008control}, and to \cite{cazacu2014controllability} when the singular potential is localized at the boundary.

\subsection{Structure of the paper}

This paper is organized as follows: 

\begin{itemize}
     \item In Section \ref{section: Strategy}, we present proof strategies. In particular, we state in Section \ref{section: Fourier decomposition and the control cost from annuli} an intermediate result concerning the cost of control from annuli for the Fourier components of system \eqref{stm : control system} when $\gamma = 1$ (Theorem \ref{thm: uniform control annulus gamma = 1}), and another one when $\gamma \in (0,1)$ (Theorem \ref{thm: uniform control annulus gamma < 1}). Section \ref{section: Spherical harmonic decomposition} is devoted to the spherical harmonic decomposition of the Fourier components of system \eqref{stm : control system}. Section \ref{section: Moments method for radial problems} presents the strategy of proof of Theorem \ref{thm: uniform control annulus gamma = 1}, and Section \ref{section: strategy proof transmuation method} the one of Theorem \ref{thm: uniform control annulus gamma < 1}. 
    \item Section \ref{section: Uniform spectral estimates} is devoted to the spectral analysis of the Fourier and spherical components of our singular Grushin operators. We first recall some eigenvalue estimates in Section \ref{section: some estimates on the eigenvalues}. Then, for $\gamma = 1$, we obtain uniform spectral gaps in Section \ref{section: A uniform spectral gap via the factorization method}, and uniform localization estimates for the eigenfunctions in Section \ref{section: Uniform localization of eigenfunctions}. 
    \item Section \ref{section: radial observability of waves} is devoted to proving observability inequalities for radial waves, in the presence of both a degenerate and singular potential (Theorem \ref{thm: radial wave observability}).
    \item Section \ref{section: proof of the controllability results in the critical degenerate case} is devoted to the proofs of first Theorem \ref{thm: uniform control annulus gamma = 1}, and then Theorem \ref{thm: main thm minimal time gamma = 1}.
    \item Section \ref{section: Proof of null-controllability results in the weakly degenerate case} is devoted to the proofs of first Theorem  \ref{thm: uniform control annulus gamma < 1}, and then Theorem \ref{thm: main thm minimal time gamma < 1}.
\end{itemize}

The proofs of our main theorems rely on the analysis of several auxiliary systems, introduced at the beginning of Section \ref{section: Strategy}, which arise from the Fourier decomposition of system \eqref{stm : control system} and from the spherical harmonic decomposition of the resulting Fourier components. The proofs of Theorems \ref{thm: main thm minimal time gamma = 1} and \ref{thm: main thm minimal time gamma < 1} thus rely on numerous intermediate results. To help the reader navigate the present work, we propose the following diagram. 

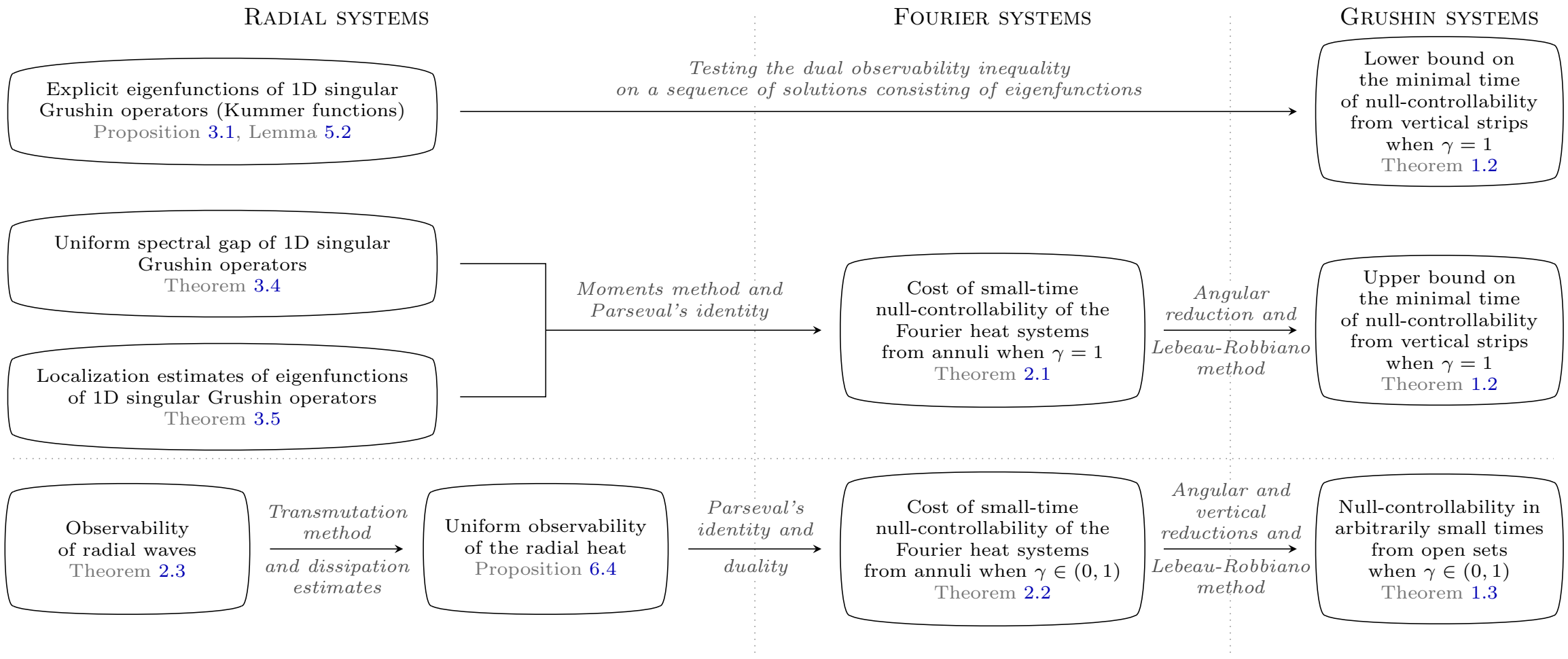
\begin{figure}[H]
    \centering
    \resizebox{\textwidth}{!}{\input{structurepaper.tex}}
    \caption{Diagram of the structure of the paper.}
    \label{fig:placeholder}
\end{figure}

\section{Proof strategies}\label{section: Strategy}

\subsection{Systems and operators of interest}

For the sake of clarity, we start by introducing the five systems of interest in the present work, as well as the Grushin operators under study. In what follows, we let $S,T > 0$, $\xi > 0$, and $\gamma, \nu,\mu \geq 0$.

\subsubsection{Systems}
First, given a solution of system \eqref{stm : control system}, from the decomposition of its solutions in a basis of eigenfunctions of $\Delta_y$, we are led to the following systems, understood as the Fourier components of system \eqref{stm : control system},
\begin{align}\label{stm : Fourier heat control system} \tag{$\mathrm{H}_{\mathrm{Fou}}^{\mathrm{ctrl}}$} 
    \left\{ \begin{array}{lcll}
     \partial_t f_\xi - \Delta_x f_\xi + \left( \xi^2|x|^{2\gamma} + \displaystyle\frac{\nu^2 - \mathsf{H}}{|x|^2} \right) f_\xi & = & \mathbbm{1}_\mathcal{U} u_\xi(t,x), & (t,x) \in (0,T) \times B(0,1),\\[6pt]
     f_\xi(t,x) & = & 0, & t \in (0,T), \ x \in \partial B(0,1), \\[6pt]
     f_\xi(0,x) & = & f_\xi^0(x), & x \in B(0,1),
    \end{array} \right.
\end{align}
where $\mathcal{U} \subset B(0,1)$. Second, we will need the forward adjoint systems to systems \eqref{stm : Fourier heat control system}, given by
\begin{align}\label{stm: Fourier heat adjoint system}\tag{$\mathrm{H}_{\mathrm{Fou}}^{\mathrm{adj}}$} 
    \left\{ \begin{array}{lcll}
     \partial_t f_\xi -\Delta_x f_\xi + \left( \xi^2|x|^{2\gamma} + \displaystyle\frac{\nu^2 - \mathsf{H}}{|x|^2} \right) f_\xi & = & 0, & (t,x) \in (0,T) \times B(0,1),\\[6pt]
     f_\xi(t,x) & = & 0, & t \in (0,T), \ x \in \partial B(0,1), \\[6pt]
     f_\xi(0,x) & = & f_\xi^0(x), & x \in B(0,1).
    \end{array} \right.
\end{align}

Third, setting $x = r\theta$ in systems \eqref{stm : Fourier heat control system} and \eqref{stm: Fourier heat adjoint system}, the decomposition of their solutions in a basis of spherical harmonics gives the following radial systems,  
\begin{align}\label{stm: radial heat control system}\tag{$\mathrm{H}_{\mathrm{Rad}}^{\mathrm{ctrl}}$} 
    \left\{ \begin{array}{lcll}
        \partial_t z - \partial_r^2 z + \left( \xi^2r^{2\gamma} + \displaystyle \frac{\mu^2 - \frac{1}{4}}{r^2} \right) z &=& \mathbbm{1}_{(a,b)}\mathsf{u}, & (t,r) \in (0,T) \times (0,1), \\
        z(t,0) = z(t,1) &=& 0, & t \in (0,T), \\
        z(0,r) &=& z_0(r), & r \in (0,1),
    \end{array} \right.
\end{align}
where $0 \leq a < b \leq 1$, and 
\begin{align}\label{stm: radial heat adjoint system}\tag{$\mathrm{H}_{\mathrm{Rad}}^{\mathrm{adj}}$} 
    \left\{ \begin{array}{lcll}
        \partial_t z - \partial_r^2 z + \left( \xi^2r^{2\gamma} + \displaystyle \frac{\mu^2 - \frac{1}{4}}{r^2} \right) z &=& 0, & (t,r) \in (0,T) \times (0,1), \\
        z(t,0) = z(t,1) &=& 0, & t \in (0,T), \\
        z(0,r) &=& z_0(r), & r \in (0,1).
    \end{array} \right.
\end{align}
Notice that system \eqref{stm: radial heat adjoint system} is the forward adjoint system to system \eqref{stm: radial heat control system}. And finally, the homogeneous waves counterparts of systems \eqref{stm: radial heat adjoint system}, 
\begin{align}\label{stm: radial wave adjoint system}\tag{$\mathrm{W}_{\mathrm{Rad}}^{\mathrm{adj}}$} 
    \left\{ \begin{array}{lcll}
        \partial_s^2 w - \partial_r^2 w + \left( \xi^2r^{2\gamma} + \displaystyle \frac{\mu^2 - 1/4}{r^2} \right) w &=& 0, & (s,r) \in (-S,S) \times (0,1), \\
        w(s,0) = w(s,1) &=& 0, & s \in (-S,S), \\
        (w, \partial_s w)(0,r) &=& (w_0, w_1)(r), & r \in (0,1).
    \end{array} \right.
\end{align}

\subsubsection{Operators}

The spatial operators in the systems above are understood as their Friedrichs realizations in $L^2$ from $C_c^\infty$. Namely, in systems \eqref{stm : Fourier heat control system} and \eqref{stm: Fourier heat adjoint system}, the operators are 
\begin{align}\label{eqn: definition G_xi}
    \begin{array}{cll}
        D(\mathcal{G}_\xi) &=&  \left\{ f \in H_{0,\nu}^1(B(0,1)), \ \mathcal{G}_\xi f \in L^2(B(0,1)) \right\}, \\[6pt]
        \mathcal{G}_\xi &=& -\Delta_x + \xi^2|x|^{2\gamma} + \displaystyle\frac{\nu^2 - \mathsf{H}}{|x|^2},
    \end{array}
\end{align}
where $H_{0,\nu}^1(B(0,1))$ is the closure of $C_c^\infty(B(0,1))$ with respect to the norm (thanks to Hardy inequality \eqref{eqn: hardy inequality}) 
\begin{align}\label{eqn: def norm nu}
    \|f\|_\nu^2 := \int_{B(0,1)} |\nabla_x f|^2 + \left( \xi^2 |x|^{2\gamma} + \frac{\nu^2 - \mathsf{H}}{|x|^2}\right) |f|^2 \ dx.
\end{align}
We do not index $H_{0,\nu}^1(B(0,1))$ with respect to $\xi$ as this space does not depend on it, and we do not index with respect to $\gamma$ as no confusion shall arise. \\

Similarly, in systems \eqref{stm: radial heat control system}, \eqref{stm: radial heat adjoint system}, and \eqref{stm: radial wave adjoint system}, the operators are
\begin{align}\label{eqn: def G radial}
    \begin{array}{cll}
        D(G_{\xi,\mu}) &=&  \left\{ f \in H_{0,\mu}^1(0,1), \ G_{\xi,\mu} f \in L^2(0,1) \right\}, \\[6pt]
        G_{\xi,\mu} &=& - \partial_r^2 + \xi^2r^{2\gamma} + \displaystyle\frac{\mu^2 - 1/4}{r^2}.
    \end{array}
\end{align}
where $H_{0,\mu}^1(0,1)$ is the closure of $C_c^\infty(0,1)$ with respect to the norm (thanks to Hardy inequality \eqref{eqn: hardy inequality}) 
\begin{align}\label{eqn: def norm mu}
    \|f\|_\mu^2 := \int_0^1 |\partial_r f|^2 + \left( \xi^2 r^{2\gamma} + \frac{\mu^2 - 1/4}{r^2}\right) |f|^2 \ dr.
\end{align}
The operator $G_{\xi,\mu}$ is self-adjoint, nonnegative by the Hardy inequality \eqref{eqn: hardy inequality}, with compact resolvent. Its spectrum consists of an increasing sequence of simple eigenvalues, which we denote by 
\begin{align*}
    0 < \lambda_{\xi,\mu,0} < \lambda_{\xi,\mu,1} < \dots < \lambda_{\xi,\mu,k} < \lambda_{\xi,\mu,k+1} < \dots,
\end{align*}
and the associated eigenfunctions, normalized in $L^2(0,1)$, are denoted by $\varphi_{\xi,\mu,k}$. \\

Let us also introduce the energy form associated with $G_{\xi,\mu}$. For every $f \in C_c^\infty(0,1)$, we define
\begin{align}\label{eqn: def energy form Q smooth functions i}
    Q_{\xi,\mu}(f) := \int_0^1 (\partial_r f(r))^2 + \left( \xi^2 r^{2\gamma} + \frac{\mu^2 - 1/4}{r^2} \right)f(r)^2 \ dr,
\end{align}
and by $\mathcal{Q}_{\xi,\mu}$ its closure for the norm $\| \cdot \|_\mu$. For $\mu = 0$, the integral form \eqref{eqn: def energy form Q smooth functions i} does not extend to every function of $H_{0,\mu}^1(0,1)$. We thus introduce the change of variable $f = \sqrt{r} g$, and we have
\begin{align}\label{eqn: def energy form Q smooth functions ii}
    Q_{\xi,\mu}(f) = \int_0^1 rg'(r)^2 + \left(\xi^2 r^{2\gamma + 1} + \frac{\mu^2}{r} \right)g(r)^2 \ dr.
\end{align}
This form extends to every element of $H_{0,\mu}^1(0,1)$. Hence, 
\begin{align}\label{eqn: def energy form Q}
    \mathcal{Q}_{\xi,\mu}(f) = \int_0^1 rg'(r)^2 + \left(\xi^2 r^{2\gamma + 1} + \frac{\mu^2}{r} \right)g(r)^2 \ dr, \quad f \in H_{0,\mu}^1(0,1), \quad g = fr^{-1/2}.
\end{align}

We thus have the following characterization of the elements of $H_{0,\mu}^1(0,1)$, 
\begin{align}\label{eqn: characterization form domain}
    H_{0,\mu}^1(0,1) = \left\{ f \in H^1_{\mathrm{loc}}((0,1]), \ \mathcal{Q}_{\xi,\mu}(f) < \infty, \ f(1) = 0 \right\}.
\end{align}
Observe that it imposes the trace condition $f(0) = 0$.

\subsection{Fourier decomposition and the control cost from annuli}\label{section: Fourier decomposition and the control cost from annuli}

We first follow the usual strategy (introduced in \cite{beauchard2014null}) of Fourier decomposition of system \eqref{stm : control system}. Assume that $\omega = \omega_x \times \Omega_y$. Denoting by $\phi_n$ the eigenfunctions of $-\Delta_y$, with associated eigenvalues $\xi_n^2$, $n \in \mathbb{N} \setminus \{0\}$, we can write any solution of system \eqref{stm : control system} and the control as 
\begin{align}\label{eqn: fourier decomposition}
    \begin{aligned}
        f(t,x,y) &= \displaystyle \sum_n f_n(t,x) \phi_n(y), \quad \text{with} \quad f_n(t,x) = \displaystyle \int_{\Omega_y} f(t,x,y) \phi_n(y) \ dy, \\
        \mathbbm{1}_{\omega}u(t,x,y) &= \displaystyle \sum_n \mathbbm{1}_{\omega_x} u_n(t,x) \phi_n(y), \quad \text{with} \quad u_n(t,x) = \displaystyle \int_{\Omega_y} u(t,x,y) \phi_n(y) \ dy.
    \end{aligned}
\end{align}
Denoting $f_n = f_{\xi_n}$ and $u_n = u_{\xi_n}$, for every $n \geq 1$, the Fourier coefficients $f_n \in C([0,T], L^2(B(0,1)))\cap L^2((0,T), H_{0,\nu}^1(B(0,1)))$ are the unique weak solutions of (see \cite[Proposition 2]{cannarsa2014null} or \cite[Proposition 2.1]{anh2016null}) systems \eqref{stm : Fourier heat control system}, with $\xi = \xi_n$ and $\mathcal{U} = \omega_x$.\\

By Parseval's identity and existence of a minimal norm control, system \eqref{stm : control system} is null-controllable in time $T > 0$ from $\omega = \omega_x \times \Omega_y$ if and only if the Fourier systems \eqref{stm : Fourier heat control system} are uniformly null-controllable with respect to $n \geq 1$ in time $T > 0$ from $\omega_x$. That is, if and only if there exists $C > 0$ such that, for every $n \geq 1$, for every $f_n^0 \in L^2(B(0,1))$, there exists a control $u_n \in L^2((0,T) \times B(0,1))$ satisfying 
    \begin{align}\label{eqn: uniform cost fourier}
        \begin{aligned}
            \|u_n\|_{L^2((0,T) \times \omega_x)} &\leq C \|f_n^0\|_{L^2(B(0,1))}, \\[6pt]
            f_n(T) &= 0.
        \end{aligned}
    \end{align}
    
To prove the non null-controllability part of Theorem \ref{thm: main thm minimal time gamma = 1}, that is a lower bound on $T(\omega_x \times \Omega_y)$, we use the dual point of view of observability (see \textit{e.g.} \cite[Theorem 2.44]{coron2007control}) to disprove \eqref{eqn: uniform cost fourier}. Namely, we equivalently prove that if $T > 0$ is not large enough, there exists no constant $C > 0$ independent of $n$ such that, for any $f_n^0 \in L^2(B(0,1))$, the solutions $f_n$ of systems \eqref{stm: Fourier heat adjoint system} with $\gamma = 1$ and $\xi = \xi_n$ satisfy 
\begin{align}\label{eqn: obs ineq strategy of proof}
    \int_{B(0,1)} |f_n(T,x)|^2 \ dx \leq C \int_0^T \int_{\omega_x} |f_n(t,x)|^2 \ dx \ dt.
\end{align}
This is done in Section \ref{section: Proof of Theorem 1 negative} by testing \eqref{eqn: obs ineq strategy of proof} against a well-chosen sequence of solutions $f_n$ consisting of eigenfunctions of $\mathcal{G}_{\xi_n}$ that we carefully estimate. \\

We now focus on how to prove that \eqref{eqn: uniform cost fourier} is true when $T$ is sufficiently large. If $\partial \Omega_y = \emptyset$, then $\xi_1^2 = 0$ is an eigenvalue of $\Delta_y$, in which case we recover in system \eqref{stm : Fourier heat control system} the heat equation in the presence of a singular potential. For every $\nu \geq 0$, it is null-controllable thanks to \cite{ervedoza2008control}. Hence, it suffices for us to prove the uniform null-controllability \eqref{eqn: uniform cost fourier} only for $\xi_n > 0$. \\

To prove \eqref{eqn: uniform cost fourier} when $\gamma = 1$, for $\xi_n > 0$, we first estimate the cost of null-controllability in any time $T > 0$ from annuli. 

\begin{theorem}\label{thm: uniform control annulus gamma = 1}
     Let $\gamma = 1$. Let $0 \leq a < b \leq 1$ and $\mathcal{A} = \{x \in B(0,1) : a<|x|<b\}$. For every $\epsilon \in (0,(b-a)/2)$, $\delta \in (0,1/2)$, and $\beta > 1$, there exists a constant $C > 0$ such that, for every $T > 0$, for every $n \geq 1$ such that $\xi_n > 0$, and every $f_{n}^0\in L^2(B(0,1))$, there exists a control $u_n$ satisfying
    \begin{align}\label{eqn: uniform cost fourier annulus gamma = 1}
        \begin{aligned}
        &\|u_n\|_{L^2((0,T)\times\mathcal{A})} \leq C \exp\left(\frac{C}{T^\beta}\right) \exp\left( \frac{\xi_n}{2}(a+2\epsilon)^2 -2(1-2\delta)(1+\nu)T\xi_n \right) \|f_n^0\|_{L^2(B(0,1))},\\
        &f_n(T) = 0,
        \end{aligned}
    \end{align}
    where $f_n$ is the solution of system \eqref{stm : Fourier heat control system} with $\xi = \xi_n$ and $\mathcal{U} = \mathcal{A}$.
\end{theorem}

Theorem \ref{thm: uniform control annulus gamma = 1} is proved in Section \ref{section: Proof of Theorem 2}, and its strategy of proof is given in Section \ref{section: Moments method for radial problems}. Observe that Theorem \ref{thm: uniform control annulus gamma = 1} directly provides the sharp upper bound on $T(\mathcal{A} \times \Omega_y)$ of Theorem \ref{thm: main thm minimal time gamma = 1}, as the cost of control appearing in \eqref{eqn: uniform cost fourier annulus gamma = 1} is uniform for sufficiently large $T$. \\ 

We will then combine \eqref{eqn: uniform cost fourier annulus gamma = 1} with spectral estimates and the Lebeau-Robbiano method (see \textit{e.g.} \cite{lebeau1995controle, le2012carleman}) to restrict the control zone to
\begin{align*}
    \mathcal{A}_\Gamma = \{x = r \theta \in B(0,1), \ r \in [0,1), \ \theta \in \mathbb{S}^{d_x - 1}: \ 0 < a < |x| < b, \ \theta \in \Gamma \subset \mathbb{S}^{d_x - 1}\} \subset \omega_x,
\end{align*}
where $\Gamma$ is a nonempty open subset of $\mathbb{S}^{d_x-1}$, $a,b,\Gamma$ are chosen so that this inclusion holds, and on which a similar cost as \eqref{eqn: uniform cost fourier annulus gamma = 1} holds. We then obtain that, for sufficiently large times, 
\begin{align*}
    \|u_n\|_{L^2((0,T) \times \mathcal{A}_\Gamma)} \leq C \|f_n^0\|_{L^2(B(0,1))}, \quad \mathcal{A}_\Gamma \subset \omega_x. 
\end{align*}

This strategy was employed from the observability point of view by means of Carleman estimates in \cite{DardeTrabut2026} in the non-singular case (we also refer to \cite{miller2010direct} in an abstract setting). \\

To prove Theorem \ref{thm: main thm minimal time gamma < 1}, we also start from a uniform cost of null-controllability from annuli for systems \eqref{stm : Fourier heat control system} with $\gamma \in (0,1)$, which,  unlike Theorem \ref{thm: uniform control annulus gamma = 1}, is uniform in $n$ for any $T > 0$. 

\begin{theorem}\label{thm: uniform control annulus gamma < 1}
    Let $\gamma \in (0,1)$. Let $0 \leq a < b \leq 1$, and $\mathcal{A} = \{x \in B(0,1) : a<|x|<b\}$. There exists $C > 0$ such that, for every $T > 0$, for every $n \geq 1$ such that $\xi_n > 0$, for every $f_n^0 \in L^2(B(0,1))$, there exists a control $u_n$ satisfying  
    \begin{align}\label{eqn: uniform cost fourier annulus gamma < 1}
        \begin{aligned}
        &\|u_n\|_{L^2((0,T)\times\mathcal{A})} \leq C \exp \left( C T^{\frac{\gamma + 1}{ \gamma - 1}} \right) \|f_n^0\|_{L^2( B(0,1))},\\
        &f_n(T) = 0,
        \end{aligned}
    \end{align}
    where $f_n$ is the solution of system \eqref{stm : Fourier heat control system} with $\xi = \xi_n$ and $\mathcal{U} = \mathcal{A}$.
\end{theorem}

Theorem \ref{thm: uniform control annulus gamma < 1} is proved in Section \ref{section: From radial uniform observability of solutions of the heat equation to uniform null-controllability from annuli}, and its strategy of proof is presented in Section \ref{section: strategy proof transmuation method}. \\

Since $\omega$ is open, it contains a set of the form $\omega_x \times \omega_y$. We then reduce the control zone of Theorem \ref{thm: uniform control annulus gamma < 1} to 
\begin{align*}
    \mathcal{A}_\Gamma = \{x = r \theta \in B(0,1), \ r \in [0,1), \ \theta \in \mathbb{S}^{d_x - 1}: \ 0 < a < |x| < b, \ \theta \in \Gamma \subset \mathbb{S}^{d_x - 1}\} \subset \omega_x,
\end{align*}
where $\Gamma$ is a nonempty open subset of $\mathbb{S}^{d_x-1}$, and $a,b,\Gamma$ are chosen so that this inclusion holds, in the exact same manner as for Theorem \ref{thm: uniform control annulus gamma = 1}. We obtain a similar cost as \eqref{eqn: uniform cost fourier annulus gamma < 1} on $\mathcal{A}_\Gamma$, obtaining null-controllability in arbitrarily small times from $\mathcal{A}_\Gamma \times \Omega_y$. Since $T > 0$ can be taken arbitrarily small, we can then reapply this strategy to reduce the control zone in the $y$-direction, as done for \cite[Theorem 1]{beauchard2014null}, recovering the null-controllability in arbitrarily small times from $\omega_x \times \omega_y \subset \omega$, which proves Theorem \ref{thm: main thm minimal time gamma < 1}.

\subsection{Spherical harmonic decomposition}\label{section: Spherical harmonic decomposition}

Consider system \eqref{stm : Fourier heat control system} with $\mathcal{U} = \mathcal{A} := \{x \in B(0,1): \ 0 \leq a < |x| < b \leq 1\}$ an annulus. We write $x = r\theta$ with $r \in (0,1)$ and $\theta \in \mathbb{S}^{d_x - 1}$. In these coordinates,
\begin{align}
    -\Delta_x = -\partial_r^2 - \frac{d_x - 1}{r}\partial_r - \frac{1}{r^2}\Delta_{\mathbb{S}^{d_x - 1}}.
\end{align}
It acts on $L^2((0,1) \times \mathbb{S}^{d_x - 1}, r^{d_x - 1}drd\theta)$, where $d\theta$ is the surface measure on the sphere $\mathbb{S}^{d_x - 1}$. Our operator thus becomes
\begin{align}
    \mathcal{G}_\xi = -\partial_r^2 - \frac{d_x - 1}{r}\partial_r - \frac{1}{r^2}\Delta_{\mathbb{S}^{d_x - 1}} + \xi^2r^{2\gamma} + \frac{\nu^2 - \mathsf{H}}{r^2}.
\end{align}
We keep the notation $\mathcal{G}_\xi$ here for simplicity as no confusion shall arise. We decompose the state and the control in an orthonormal basis of eigenfunctions of $-\Delta_{\mathbb{S}^{d_x - 1}}$. Recall (see \textit{e.g.} \cite[Chapter 3.C]{berger2006spectre}) that its eigenvalues are given by $\ell(d_x + \ell - 2)$, $\ell \geq 0$, and the corresponding eigenspaces $V_\ell$ are of dimension $m_\ell$, spanned by the eigenfunctions $Y_{\ell,m}$, $1 \leq m \leq m_\ell$,
\begin{align*}
    -\Delta_{\mathbb{S}^{d_x - 1}}Y_{\ell,m} = \ell(d_x + \ell - 2)Y_{\ell,m}, \quad \text{for every } \ell \geq 0 \quad \text{and} \quad 1 \leq m \leq m_\ell.
\end{align*}
We write
\begin{align*}
    f_\xi(t,r\theta) = \sum_{\ell \geq 0}\sum_{m=1}^{m_\ell} F_{\xi,\ell,m}(t,r)Y_{\ell,m}(\theta), \quad \text{and} \quad u_\xi(t,r\theta) = \sum_{\ell \geq 0}\sum_{m=1}^{m_\ell} U_{\xi,\ell,m}(t,r)Y_{\ell,m}(\theta),
\end{align*}
where
\begin{align*}
    F_{\xi, \ell,m}(t,r) = \int_{\mathbb{S}^{d_x - 1}} f_\xi(t,r\theta)Y_{\ell,m}(\theta) \ d\theta, \quad U_{\xi,\ell,m}(t,r) = \int_{\mathbb{S}^{d_x - 1}} u_\xi(t,r\theta)Y_{\ell,m}(\theta) \ d\theta.
\end{align*}
Since $\mathbbm{1}_{\mathcal{A}}(r\theta) = \mathbbm{1}_{(a,b)}(r)$, the coefficients $F_{\xi, \ell,m}$ solve
\begin{align*}
    \partial_t F_{\xi,\ell,m} - \partial_r^2 F_{\xi,\ell,m} - \frac{d_x - 1}{r}\partial_r F_{\xi,\ell,m} + \left(\xi^2r^{2\gamma} + \frac{\nu^2 - \mathsf{H} + \ell(d_x + \ell - 2)}{r^2}\right)F_{\xi,\ell,m} = \mathbbm{1}_{(a,b)}U_{\xi,\ell,m},
\end{align*}
in $L^2((0,1), r^{d_x - 1}dr)$. Making the change of variables 
\begin{align*}
    F_{\xi,\ell,m} = r^{-(d_x - 1)/2}\mathsf{F}_{\xi,\ell,m}, \quad  U_{\xi,\ell,m} = r^{-(d_x - 1)/2}\mathsf{u}_{\xi,\ell,m},
\end{align*}
and developing the Hardy constant $\mathsf{H} = (d_x - 2)^2/4$, we obtain, for every $(\xi,\ell,m) \in (0,+\infty) \times \mathbb{N} \times \{1,\dots, m_\ell\}$, that $\mathsf{F}_{\xi,\ell,m}$ solves system \eqref{stm: radial heat control system} with $\xi = \xi_n$, $\mu^2 = \nu_\ell^2 = \nu^2 + \ell(d_x + \ell - 2)$, and initial state $\mathsf{F}_{\xi,\ell,m}(0) = \mathsf{F}_{\xi,\ell,m}^0$. \\

We briefly justify the reduction to system \eqref{stm: radial heat control system}, following \cite{beauchard2014null}. Denote by $\mathcal{Q}_{\xi,\mu}(\cdot,\cdot)$ the symmetric bilinear form associated with \eqref{eqn: def energy form Q}, that is,
\begin{align*}
    \mathcal{Q}_{\xi,\mu}(f,h) &:= \int_0^1 r g'(r) k'(r) + \left( \xi^2 r^{2\gamma+1} + \frac{\mu^2}{r} \right) g(r) k(r) \ dr, \\
    &\quad f,h \in H_{0,\mu}^1(0,1), \ g = r^{-1/2}f, \ k = r^{-1/2}h.
\end{align*}
For any $(\ell,m)$, system \eqref{stm: radial heat control system} with $\mu^2 = \nu_\ell^2$ admits a unique weak solution $\tilde{\mathsf{F}}_{\xi,\ell,m} \in C^0([0,T]; L^2(0,1)) \cap L^2(0,T;H_{0,\nu_\ell}^1(0,1))$, that satisfies $\tilde{\mathsf{F}}_{\xi,\ell,m}(0) = \mathsf{F}_{\xi,\ell,m}^0$ and, for every $\psi \in H_{0,\nu_\ell}^1(0,1)$,
\begin{align*}
    \frac{d}{dt}\int_0^1 \tilde{\mathsf{F}}_{\xi,\ell,m}(t,r) \psi(r) \ dr + \mathcal{Q}_{\xi,\nu_\ell}\left( \tilde{\mathsf{F}}_{\xi,\ell,m}(t), \psi \right) = \int_a^b \mathsf{u}_{\xi,\ell,m}(t,r) \psi(r) \ dr \quad \text{in } \mathcal{D}'(0,T).
\end{align*}

Given $f \in C^0([0,T], L^2(B(0,1))) \cap L^2((0,T),H_{0,\nu}^1(B(0,1)))$ a solution of system \eqref{stm : Fourier heat control system}, then one verifies (following for instance \cite[Proposition 2]{beauchard2014null}) that $\mathsf{F}_{\xi,\ell,m}$ is the unique weak solution of system \eqref{stm: radial heat control system} with $\mu^2 = \nu_\ell^2$, control $\mathsf{u}_{\xi,\ell,m}$, and initial state $\mathsf{F}_{\xi,\ell,m}^0$.

\subsection{Moments method for radial problems}\label{section: Moments method for radial problems}

In this section, we explain how to prove Theorem \ref{thm: uniform control annulus gamma = 1}, and we keep the notations of the radial decomposition given in the previous subsection. We apply the moments method at the level of system \eqref{stm: radial heat control system} with $\mu^2 = \nu_\ell^2 = \nu^2 + \ell(d_x + \ell - 2)$. For $\mu = \nu_\ell$, denote the eigenvalues of $G_{\xi,\nu_\ell}$ given in \eqref{eqn: def G radial} by $\lambda_{\xi,\ell,k} := \lambda_{\xi, \nu_\ell, k}$, $k\geq 0$.\\

Let $\mathsf{u}_{\xi,\ell,m} \in L^2((0,T) \times (a,b))$, and let $\mathsf{F}_{\xi,\ell,m}$ be a solution of system \eqref{stm: radial heat control system}. We test system \eqref{stm: radial heat control system} against the solutions $t \mapsto e^{-(T-t)\lambda_{\xi,\ell,k}}\varphi_{\xi, \ell, k}$ of the backward problem, and we obtain, for every $k \in \mathbb{N}$,
\begin{align*}
    \begin{aligned}
        \int_0^T \int_a^b \mathsf{u}_{\xi,\ell,m} (t,r) e^{-(T-t)\lambda_{\xi,\ell,k}}\varphi_{\xi, \ell, k}(r) \ dr \ dt  &= \int_0^1 \mathsf{F}_{\xi,\ell,m}(T,r) \varphi_{\xi,\ell,k}(r) \ dr \\
        &- \int_0^1 e^{-T \lambda_{\xi,\ell,k}} \mathsf{F}_{\xi,\ell,m}^0(r) \varphi_{\xi,\ell,k}(r) \ dr.
    \end{aligned}
\end{align*}

Since $(\varphi_{\xi,\ell,k})_{k \geq 0}$ is a Hilbert basis of $L^2(0,1)$, then $\mathsf{F}_{\xi,\ell,m}(T) = 0$ if and only if the control $\mathsf{u}_{\xi,\ell,m}$ satisfies, for every $k \in \mathbb{N}$,
\begin{align}\label{eqn: moment equation}
    \int_0^T \int_a^b \mathsf{u}_{\xi,\ell,m} (t,r) e^{-(T-t)\lambda_{\xi,\ell,k}}\varphi_{\xi, \ell, k}(r) \ dr \ dt  = - \int_0^1 e^{-T \lambda_{\xi,\ell,k}} \mathsf{F}_{\xi,\ell,m}^0(r) \varphi_{\xi,\ell,k}(r) \ dr.
\end{align}
Introduce formally the biorthogonal family $(q_{\xi, \ell,j})_{j \geq 0} \subset  L^2(0,T)$ which satisfies 
\begin{align}
    \int_0^T  e^{-\lambda_{\xi,\ell,k}t}q_{\xi, \ell,j}(t) \ dt = \delta_{jk}, \quad \text{for every } j,k \in \mathbb{N}.
\end{align}
Then, a control solving \eqref{eqn: moment equation} formally is given by 
\begin{align}
    \mathsf{u}_{\xi,\ell,m}(t,r) = - \sum_{j \geq 0} e^{-\lambda_{\xi,\ell,j}T} \langle \mathsf{F}_{\xi,\ell,m}^0,\varphi_{\xi,\ell,j} \rangle_{L^2(0,1)} \frac{ \varphi_{\xi,\ell,j}(r)} { \|\varphi_{\xi,\ell,j}\|_{L^2(a,b)}^2}
    q_{\xi, \ell,j}(T-t).
\end{align}

We therefore need two ingredients:
\begin{enumerate}[label = (\roman*)]
    \item Existence of, and estimates for, the biorthogonal family $(q_{\xi,\ell,j})_{j \geq 0}$. This requires in particular a gap between the eigenvalues. Namely, that there exists some $\kappa_\ell > 0$ such that for every $k \in \mathbb{N}$, $\lambda_{\xi,\ell,k+1} - \lambda_{\xi,\ell,k} \geq \kappa_\ell$.
    \item Some lower bound for the $L^2$-norms of the eigenfunctions in the control zone (localization estimates). 
\end{enumerate}

The additional difficulty is that these ingredients must be obtained uniformly in $\ell$. The uniform gap estimate is obtained in Section \ref{section: A uniform spectral gap via the factorization method}, and the localization estimates in Section \ref{section: Uniform localization of eigenfunctions}. We obtain that there exists $C(\xi,T) > 0$ such that, for every $\ell,m$, we have 
\begin{align}
    \|\mathsf{u}_{\xi,\ell,m} \|_{L^2((0,T) \times (a,b))}^2 \leq C(\xi,T) \|\mathsf{F}_{\xi,\ell,m}^0\|_{L^2(0,1)}^2
\end{align}
Taking $\xi = \xi_n$ and indexing by $n$, the control of Theorem \ref{thm: uniform control annulus gamma = 1} is then
\begin{align}
    u_n(t,r\theta) = r^{-(d_x - 1)/2} \sum_{\ell \geq 0}\sum_{m=1}^{m_\ell} \mathsf{u}_{n,\ell,m}(t,r)Y_{\ell,m}(\theta), \quad (t,r\theta) \in (0,T)\times\mathcal{A}.
\end{align}
Indeed, using Parseval's identity, we have
\begin{align*}
    \|u_n\|_{L^2((0,T) \times \mathcal{A})}^2 = \sum_{\ell \geq 0}\sum_{m=1}^{m_\ell}\| \mathsf{u}_{n,\ell,m}\|_{L^2((0,T)\times(a,b))}^2 &\leq C(n,T)\sum_{\ell \geq 0}\sum_{m=1}^{m_\ell} \|\mathsf{F}_{n,\ell,m}^0 \|_{L^2(0,1)}^2 \\
    &= C(n,T) \|f_n^0\|_{L^2(B(0,1))}^2.
\end{align*}

\subsection{Transmutation method}\label{section: strategy proof transmuation method}

We present here the strategy of proof of Theorem \ref{thm: uniform control annulus gamma < 1}. We employ the dual observability point of view. We equivalently prove that, letting $\gamma \in (0,1)$, for every $0 \leq a < b \leq 1$, setting $\mathcal{A} := \{a < |x| < b\}$, there exists $C > 0$ such that, for every $T > 0$, for every $n \in \mathbb{N} \setminus \{0\}$ such that $\xi_n > 0$, for every $f_n^0 \in L^2(B(0,1))$, the solutions $f_n$ of systems \eqref{stm: Fourier heat adjoint system} with $\xi = \xi_n$ satisfy
    \begin{align}\label{eqn: unif obs annulus weak degen 0}
        \|f_n(T)\|_{L^2(B(0,1))}^2 \leq C \exp \left( C T^{\frac{\gamma + 1}{ \gamma - 1}} \right) \|f_n\|_{L^2( (0,T) \times \mathcal{A})}^2.
    \end{align}

To do so, we use the transmutation method (see \textit{e.g.} \cite{ervedoza2011sharp,ervedoza2011observability}) at the level of the radial decomposition of our systems. Namely, given a solution $z$ of the radial heat system \eqref{stm: radial heat adjoint system}, we can construct a solution $w$ of the radial wave system \eqref{stm: radial wave adjoint system} with initial state $(w,w_s)(0,r) = (0,w_1)(r)$. We thus first obtain an observability inequality for such solutions of system \eqref{stm: radial wave adjoint system}.

\begin{theorem}\label{thm: radial wave observability}
    Let $\gamma \geq 0$. Let $0 \leq a < b \leq 1$ and $S \geq 2$. There exists a constant $C > 0$ such that, for every $\xi > 0$, $\ell \in \mathbb{N}$ and $w_1 \in L^2(0,1)$, the associated solution of \eqref{stm: radial wave adjoint system}, with $\mu^2 = \nu_\ell^2$, and with initial state $w_0 = 0$ satisfies
    \begin{align}\label{eqn: wave observability annulus}
        \|w_1\|_{L^2(0,1)}^2 \leq C \exp\left( C \left( \xi + \nu_\ell \ln(2+\xi) \right) \right) \int_a^b \int_{-S}^S w_s(s,r)^2 + w(s,r)^2  \ ds \ dr.
    \end{align}
\end{theorem}

Theorem \ref{thm: radial wave observability} is proved in two steps. The first step is an observability inequality near the singularity $x = 0$, using a cutoff argument and the results of \cite{vancostenoble2009hardy}. The second step is an observability inequality away from $x = 0$, using the lateral propagation of energy as in \cite[Proposition 4.5]{beauchard20152d}. We glue these two inequalities. \\

Then, using the transmutation kernel and its estimates, we transfer the observability inequality \eqref{eqn: wave observability annulus} to an observability inequality for the solutions of systems \eqref{stm: radial heat adjoint system}. We obtain the desired uniform observability inequality for the solutions of the heat systems \eqref{stm: Fourier heat adjoint system} by using the dissipation of the heat and Parseval's identity, leading to \eqref{eqn: unif obs annulus weak degen 0}, and equivalently to Theorem \ref{thm: uniform control annulus gamma < 1}.

\section{Uniform spectral estimates}\label{section: Uniform spectral estimates}

This section is dedicated to the spectral analysis of the operators $G_{\xi,\mu}$ introduced in \eqref{eqn: def G radial}. We first recall in Section \ref{section: some estimates on the eigenvalues} the form of the eigenfunctions of $G_{\xi,\mu}$, and some estimates for its eigenvalues. Then, for $\gamma = 1$, we obtain in Section \ref{section: A uniform spectral gap via the factorization method} a uniform spectral gap for the eigenvalues of $G_{\xi,\mu}$, and some uniform localization estimates for its eigenfunctions in Section \ref{section: Uniform localization of eigenfunctions}. 

\subsection{Some estimates on the eigenvalues}\label{section: some estimates on the eigenvalues}

Let us recall some facts about the eigenfunctions and eigenvalues of $G_{\xi,\mu}$. The Kummer equation is defined by, for $a,b,z \in \mathbb{C}$,
\begin{align}\label{eqn: Kummer equation}
    z w''(z) + \left( b - z \right) w'(z) - a w(z) = 0.
\end{align}
The Kummer function is a solution of \eqref{eqn: Kummer equation}, defined by, as long as $b$ is not a nonpositive integer,
\begin{align}\label{eqn: definition kummer function}
    M (a,b,z) &= \sum_{k \geq 0} \frac{(a)_k}{(b)_k k!}z^k, 
\end{align}
where $(\cdot)_k = \cdot(\cdot+1)...(\cdot+k-1)$, $(\cdot)_0 = 1$, is the Pochhammer symbol. 

\begin{proposition}{\cite[Theorem 1.1]{romankummer2025}}\label{prop: exact form eigenfunctions singular + cond eigenv}
    Let $\gamma = 1$. For every $(\xi,\mu,k) \in (0,+\infty) \times [0,+\infty) \times \mathbb{N}$, the (non-normalized) eigenfunctions $\phi_{\xi,\mu,k}$ of $G_{\xi,\mu}$, associated with $\lambda_{\xi,\mu,k}$, are  
    \begin{align}
    \phi_{\xi,\mu,k} = A e^{- \xi r^2/2} r^{\frac{1}{2} + \mu} M \left(-\frac{\lambda_{\xi,\mu,k} - 2\xi (1 + \mu)}{4\xi}, 1 + \mu, \xi r^2 \right), \quad A \in \mathbb{R} \setminus \{0\},
    \end{align}
where $\lambda_{\xi, \mu ,k} \geq 0$ is an eigenvalue if and only if
\begin{align}\label{eqn : condition eigenvalue zero kummer}
    M \left(-\frac{\lambda_{\xi,\mu,k} - 2\xi (1 + \mu)}{4\xi},1 + \mu, \xi \right) = 0.
\end{align}
\end{proposition}

\begin{proposition}\label{prop: lower bounds eigenvalue}
    Let $\gamma = 1$. For every $(\xi,\mu,k) \in (0,+\infty) \times [0,+\infty) \times \mathbb{N}$, we have 
        \begin{align}\label{eqn: lower bound eigenfunction max}
            \lambda_{\xi,\mu,k} \geq \max \left(2\xi(1+\mu), \mu^2, j_{0,k}^2 \right) \geq 4,
        \end{align}
        where $j_{0,k} > 0$ is the $(k+1)$-th positive zero of the Bessel function $J_0$ of the first kind. Moreover, there exists $c > 0$ such that 
        \begin{align}\label{eqn: lower bound eigenfunction k2}
            \lambda_{\xi,\mu,k} \geq c(k+1)^2.
        \end{align}
\end{proposition}

\begin{proof}
    The bound $\lambda_{\xi,\mu,k} \geq 2\xi(1+\mu)$ follows from \cite[Theorem 1.3]{romankummer2025}. \\
    
    The bound $\lambda_{\xi,\mu,k} \geq j_{0,k}^2$ follows from the min-max theorem. Indeed, we have that $0 \leq \mu_0 \leq \mu_1$ implies $H_{0,\mu_1}^1(0,1) \subset H_{0,\mu_0}^1(0,1)$ and $\mathcal{Q}_{\xi,\mu_1} \geq \mathcal{Q}_{\xi,\mu_0}$, from which it follows that for every $\mu \geq 0$, $\lambda_{\xi,\mu,k} \geq \lambda_{0,0,k}$. The eigenvalues of $G_{0,0}$ are known (see \textit{e.g.} \cite[Proposition 2]{martinez2019cost}) and are given by $j_{0,k}^2$, where $j_{0,k}$ is the $(k+1)$-th zero of the Bessel function $J_0$ of the first kind.\\
    
    We now verify that $\lambda_{\xi,\mu,k} \geq \mu^2$. It follows from a min-max argument using that, for every $f \in H_{0,\mu}^1(0,1)$, we have 
    \begin{align*}
        \mathcal{Q}_{\xi,\mu}(f)  \geq \mu^2 \int_0^1 \frac{g(r)^2}{r} \ dr \geq \mu^2 \|f\|_{L^2(0,1)}^2, \quad g = fr^{-1/2}.
    \end{align*}
    
    To conclude on \eqref{eqn: lower bound eigenfunction max}, we use the numerical approximation $j_{0,0} \approx 2.405$ from \cite[Table 9.5]{abramowitz1948handbook}, so that $j_{0,k}^2 \geq j_{0,0}^2 \geq 4$. \\

    For the second lower bound \eqref{eqn: lower bound eigenfunction k2}, we use that from \cite[Chapter 15.52, page 505]{watson1922treatise}, we have (we index from $k = 0$, while in \cite{watson1922treatise}, the indexing start at $k = 1$)
    \begin{align*}
        j_{0,k} = \left( k + \frac{3}{4}\right)\pi + O \left(\frac{1}{k} \right), \quad \text{as } k \to +\infty. 
    \end{align*}
    Hence, there exists $c > 0$ such that $j_{0,k}^2 \geq c(k+1)^2$, and \eqref{eqn: lower bound eigenfunction k2} follows by substituting this into \eqref{eqn: lower bound eigenfunction max}. 
\end{proof}

We also have the following estimate for any $\gamma > 0$. 

\begin{proposition}\label{prop: lower bound eigenvalues any gamma}
    Let $\gamma > 0$. There exists $C > 0$ such that, for every $(\xi,\mu,k) \in (0,+\infty) \times [0,+\infty) \times \mathbb{N}$, we have
    \begin{align}
        \lambda_{\xi,\mu,k} \geq \max \left( C\xi^{2/(\gamma  +1)},\mu^2 \right).
    \end{align}
\end{proposition}

\begin{proof}
    The bound $\lambda_{\xi,\mu,k} \geq \mu^2$ is proved exactly as in the proof of Proposition \ref{prop: lower bounds eigenvalue}. Let us prove the other lower bound. \\
    
    Let $\tau_\xi := \xi^{1/(1 + \gamma)}$. Making the change of variable $s = \tau_\xi r$, and setting $g(r) = \tau_\xi h(\tau_\xi r)$, $f(r) = \sqrt{r}g(r)$, we have 
    \begin{align*}
        Q_{\xi,\mu}(f) = \tau_\xi^2 \int_0^{\tau_\xi} s h'(s)^2 + \left( s^{2\gamma + 1} + \frac{\mu^2}{s} \right) h(s)^2 \ ds.
    \end{align*}
    Extending $h$ by $0$ on $(\tau_\xi,\infty)$, we deduce using the min-max theorem that
    \begin{align*}
        \lambda_{\xi, \mu, k} \geq \lambda_{\xi,\mu, 0} \geq C_0\xi^{2/(1+\gamma)},
    \end{align*}
    where $C_0$ is the first eigenvalue of the Friedrichs extension of $-\partial_r^2 + r^{2\gamma} + (\mu^2 - 1/4)/r^2$ with domain $C_c^\infty(0,+\infty)$. Using the same monotonicity argument as in the proof of Proposition \ref{prop: lower bounds eigenvalue}, we have $C_0 \geq C$, where $C > 0$ is the first eigenvalue of the Friedrichs extension of $-\partial_r^2 + r^{2\gamma} - 1/(4r^2)$ with domain $C_c^\infty(0,+\infty)$.
\end{proof}

\subsection{A uniform spectral gap via the factorization method}\label{section: A uniform spectral gap via the factorization method}

In this section, we always assume $\gamma = 1$. The aim of this section is to provide a uniform gap for the eigenvalues of $G_{\xi,\mu}$ given in \eqref{eqn: def G radial}, with respect to $\mu,k$. The method we use here stems from the so-called factorization method introduced by E. Schrödinger in \cite{schrodinger1940method} (see also \cite{infeld1951factorization} for a systematic study of the method, in particular to \cite[Section 5.3]{infeld1951factorization} in our case). \\

The gap is given in Theorem \ref{thm: gap uniform} below, and was obtained independently and simultaneously by Helffer, Kachmar and Nicoleau in \cite{helffer2026spectralasymptoticssingularharmonic}. Since our proof differs slightly, we include it to keep the paper self-contained. Moreover, to the best of our knowledge, this method has not been used in connection with the moments method, in particular for the control of Grushin equations (see \textit{e.g.} \cite{allonsius2021analysis, LaurentLeauthaud23, tamekue2026minimal}). We therefore briefly introduce it, which also motivates the introduction of certain operators.

\begin{theorem}\label{thm: gap uniform}
    Let $\gamma = 1$. For every $(\xi,\mu,k) \in (0,+\infty) \times [0,+\infty) \times \mathbb{N}$, we have 
    \begin{align}
        \lambda_{\xi,\mu,k+1} - \lambda_{\xi,\mu,k} \geq 4\xi.
    \end{align}
\end{theorem}

The idea from \cite{schrodinger1940method} is to factorize an operator $H = - \partial_x^2 + V(x)$ into two operators, in order to obtain $H = ST + c$ for some constant $c \in \mathbb{R}$. For $V(x) = x^2$ and $c = 1$, the factors $S$ and $T$ formally map an eigenfunction $\varphi$ of energy $\lambda$ to functions $S\varphi$ and $T\varphi$ of energies $\lambda+2$ and $\lambda-2$ respectively. We wish to do the same for our singular operator in order to obtain an operator $D$ that sends an eigenfunction $\varphi$, with energy $\lambda$, onto $D\varphi$, with energy $\lambda - \kappa$, where $\kappa$ will be an estimated gap (at least formally since we will need to check if $D\varphi$ belongs to the form domain). The following computations are purely algebraic, and similar to what can be found in, for instance, \cite{cabrera2008beyond, quesne2009solvable,infeld1951factorization}. \\

We look for a function $W$ and two operators $G_+ = -\partial_r + W$, $G_- = \partial_r + W$, such that, for some $c \in \mathbb{R}$, we have $G_{\xi,\mu} = G_+G_- + c$. We have
\begin{align*}
    G_+G_- &= (-\partial_r + W)(\partial_r + W) = -\partial_r^2 - W' + W^2, \\
    G_-G_+ &= (\partial_r + W)(-\partial_r + W) = -\partial_r^2 + W' + W^2.
\end{align*}

Hence, we solve for $W$,
\begin{align}\label{eqn: ricatti}
    -W' + W^2 = \xi^2r^2 + \frac{\mu^2 - 1/4}{r^2} - c, \quad c \in \mathbb{R}.
\end{align}

We obtain two possible solutions $W_+,W_-$ of \eqref{eqn: ricatti}, given by 
\begin{align*}
    W_+ &= \xi r - \frac{\mu + \frac{1}{2}}{r}, \quad W_- = \xi r + \frac{\mu - \frac{1}{2}}{r}.
\end{align*}

These give us two possible factorizations using the four factors
\begin{align}\label{eqn: def factors}
    \begin{aligned}
        G_+^{(\mu)} &= -\partial_r + W_+ = -\partial_r + \xi r - \frac{\mu + \frac{1}{2}}{r}, \quad G_-^{(\mu)} = \partial_r + W_+ = \partial_r + \xi r - \frac{\mu + \frac{1}{2}}{r}, \\
    \tilde{G}_+^{(\mu)} &= -\partial_r + W_- = -\partial_r + \xi r + \frac{\mu - \frac{1}{2}}{r}, \quad \tilde{G}_-^{(\mu)} = \partial_r + W_- = \partial_r + \xi r + \frac{\mu - \frac{1}{2}}{r}.
    \end{aligned}
\end{align}
Easy computations lead to 
\begin{align}\label{eqn: factorization G}
    \begin{array}{lllllll}
         G_+^{(\mu)}G_-^{(\mu)} &=& G_{\xi,\mu} - 2\xi(\mu+1),& \quad& G_-^{(\mu)}G_+^{(\mu)} &=& G_{\xi,\mu+1} - 2\xi\mu, \\[6pt]
    \tilde{G}_+^{(\mu)}\tilde{G}_-^{(\mu)} &=& G_{\xi,\mu} + 2\xi(\mu-1),& \quad& \tilde{G}_-^{(\mu)}\tilde{G}_+^{(\mu)} &=&  G_{\xi,\mu-1} + 2\xi\mu.
    \end{array}
\end{align}

We introduced above the notation $(\mu)$ as a superscript indicating the parameter so that it is clear later that, taking the last equation $\tilde{G}_-^{(\mu)}\tilde{G}_+^{(\mu)} =  G_{\xi,\mu-1} + 2\xi\mu$ in \eqref{eqn: factorization G} above for instance, this gives $G_{\xi,\mu} = \tilde{G}_-^{(\mu+1)}\tilde{G}_+^{(\mu+1)} - 2\xi(\mu+1)$.\\

Let now $\varphi \neq 0$ solve $G_{\xi,\mu} \varphi = \lambda \varphi$. Observe, as explained in \cite{schrodinger1940method}, that applying the right-hand side factor of an equation to the whole equation, we have, for instance for the first factorization, 
\begin{align*}
     G_-^{(\mu)}G_+^{(\mu)}G_-^{(\mu)} \varphi = G_-^{(\mu)} G_{\xi,\mu}\varphi - 2\xi(\mu+1)(G_-^{(\mu)}\varphi) \Longrightarrow G_{\xi,\mu+1} (G_-^{(\mu)} \varphi) =  (\lambda - 2\xi) (G_-^{(\mu)}\varphi).
\end{align*}
Thus, $G_-^{(\mu)}$ takes an eigenfunction $\varphi$ at energy $\lambda$ for $G_{\xi,\mu}$, and sends it to an eigenfunction $G_-^{(\mu)}\varphi$ at energy $\lambda - 2\xi$ for $G_{\xi,\mu+1}$, at least formally. Doing this computation with the other equations in \eqref{eqn: factorization G}, we obtain 
\begin{align}
    \begin{array}{llll}
        \text{if }& G_{\xi,\mu}\varphi = \lambda \varphi, \quad &\text{then } &G_{\xi,\mu+1} (G_-^{(\mu)} \varphi) =  (\lambda - 2\xi) (G_-^{(\mu)}\varphi), \\[6pt]
    \text{if }& G_{\xi,\mu+1}\varphi = \lambda \varphi, \quad &\text{then } &G_{\xi,\mu}(G_+^{(\mu)}\varphi) = (\lambda  + 2\xi)(G_+^{(\mu)}\varphi), \\[6pt]
    \text{if }& G_{\xi,\mu}\varphi = \lambda \varphi, \quad &\text{then } &G_{\xi,\mu-1}(\tilde{G}_-^{(\mu)}\varphi) = (\lambda  - 2\xi)(\tilde{G}_-^{(\mu)}\varphi), \\[6pt]
    \text{if }& G_{\xi,\mu-1}\varphi = \lambda \varphi, \quad &\text{then } &G_{\xi,\mu}(\tilde{G}_+^{(\mu)} \varphi) = (\lambda  + 2\xi)(\tilde{G}_+^{(\mu)}\varphi). 
    \end{array}   
\end{align}

Thus, for our problem, if we are given $\varphi \neq 0$ that solves $G_{\xi,\mu} \varphi = \lambda \varphi$, the composition $D := G_-^{(\mu-1)}\tilde{G}_-^{(\mu)}$ sends $\varphi$ to $D\varphi$ that formally solves $G_{\xi,\mu}(D\varphi) = (\lambda - 4\xi)D\varphi$. The operator $D$ is explicitly given by 
\begin{align}\label{eqn: def D}
    D = \partial_r^2 + 2\xi r \partial_r + \xi^2r^2 - \frac{\mu^2 - 1/4}{r^2} + \xi.
\end{align}
Observe that applying $D$ to the eigenfunction $\varphi$, we obtain 
\begin{align}\label{eqn: Du}
    D\varphi = 2\xi r \partial_r \varphi + 2\xi^2r^2 \varphi +\xi \varphi - \lambda \varphi.
\end{align}

We can now proceed to the proof.

\begin{proof}[Proof of Theorem \ref{thm: gap uniform}] 

We use the notation $f' = \partial_r f$. Fix some parameters $(\xi, \mu, k) \in (0,+\infty) \times [0,\infty) \times \mathbb{N}$. Let $\varphi \in D(G_{\xi,\mu}) \setminus \{0\}$ solve $G_{\xi,\mu} \varphi = \lambda_{\xi,\mu, k+1} \varphi$. Then $v := D\varphi$ is given by \eqref{eqn: Du},
\begin{align*}
    v = 2\xi r \varphi' + \left( 2\xi^2r^2  +\xi  - \lambda_{\xi,\mu,k+1} \right) \varphi \in C^\infty((0,1]), 
\end{align*}
and satisfies, 
\begin{align*}
    -v'' + \left( \xi^2 r^2 + \frac{\mu^2 - 1/4}{r^2} \right) v = (\lambda_{\xi,\mu,k+1}-4\xi)v.
\end{align*}

Using Proposition \ref{prop: exact form eigenfunctions singular + cond eigenv}, we have $r\varphi'(r) \sim r^{\mu + 1/2}$ as $r \to 0^+$. Thus, we observe that $v(0) = 0$. However, $v(1) = 2\xi \varphi'(1) \neq 0$, and thus $v$ is not an eigenfunction, but we can build test functions from it. \\

The eigenfunction $\varphi$ has exactly $k+1$ simple zeros in $(0,1)$. Denoting the zeros of $\varphi$ by $z_{j,\varphi}$, they satisfy $0 = z_{0,\varphi} < z_{1,\varphi} < \dots < z_{k+1,\varphi} < z_{k+2,\varphi} =1$. At each $z_{j,\varphi}$, $j \geq 1$, we have $v(z_{j,\varphi}) = 2\xi z_{j,\varphi} \varphi'(z_{j,\varphi}) \neq 0$ which alternates in sign with respect to $j$. Thus, $v$ has at least a zero $z_{j,v}$ in each $(z_{j,\varphi}, z_{j+1,\varphi})$, $j \in \{1, \dots, k + 1\}$. Setting $z_{0,v} = 0$, we thus define
\begin{align}\label{eqn: def vj}
    v_j = v\mathbbm{1}_{(z_{j-1,v}, z_{j,v})}, \quad j \in \{1, \dots, k+1\}.
\end{align}

We claim that 
\begin{align}\label{eqn: vj in form domain}
    v_j \in H_{0,\mu}^1(0,1), \quad \forall j \in \{1, \dots, k+1\}.
\end{align}
We conclude the proof assuming \eqref{eqn: vj in form domain}. Since the $v_j$ have disjoint support, they are linearly independent, and the vector space $E_k := \mathrm{Span}(v_1,\dots,v_{k+1})$ has dimension $\mathrm{Dim}(E_k) =k+1$. Moreover, $\mathcal{Q}_{\xi,\mu}$ introduced in \eqref{eqn: def energy form Q} is well-defined on $E_k$, and for every $f \in E_k$, we have
\begin{align*}
    \mathcal{Q}_{\xi,\mu}(f)  = (\lambda_{\xi,\mu,k+1} - 4\xi) \| f \|_{L^2(0,1)}^2.
\end{align*}
Applying the min-max theorem, we have 
\begin{align*}
    \lambda_{\xi,\mu,k} =  \underset{\underset{\operatorname{dim}(V) = k+1}{V \subset H_{0,\mu}^1(0,1)}}{\min} \underset{f \in V \setminus \{0\}}{\max} \frac{\mathcal{Q}_{\xi,\mu}(f)}{\|f\|_{L^2(0,1)}^2} \leq \underset{f \in E_k \setminus \{0\}}{\max} \frac{\mathcal{Q}_{\xi,\mu}(f)}{\|f\|_{L^2(0,1)}^2} = \lambda_{\xi,\mu,k+1} - 4\xi.
\end{align*}
This concludes the proof assuming \eqref{eqn: vj in form domain}. It remains to prove \eqref{eqn: vj in form domain}.\\

We show that each $v_j$ belongs to the form domain $H_{0,\mu}^1(0,1)$ of $G_{\xi,\mu}$, whose characterization is given in \eqref{eqn: characterization form domain}. First, observe that this is true for any $j \geq 2$, for any $\mu \geq 0$, as the support of $v_j$ is located away from zero where $\varphi$ is smooth. Hence, we only need to verify \eqref{eqn: vj in form domain} for $v_1$.\\

Assume that $\mu > 0$. We have $\varphi \in H_0^1 \cap C^\infty(0,1)$. Thus, on the one hand, $\left( 2\xi^2r^2  +\xi  - \lambda_{\xi,\mu,k+1} \right) \varphi \in H_0^1(0,1)$. On the other hand, it also follows that $r\varphi' \in L^2\cap C^\infty(0,1)$, and
    \begin{align*}
        (r\varphi')' = \varphi' + r\varphi'' = \varphi' + \left(\xi^2r^3 + \frac{\mu^2 - 1/4}{r} - \lambda_{\xi,\mu,k+1}r \right) \varphi \in L^2(0,1),
    \end{align*}
where the inclusion follows from Hardy's inequality. Hence, since $v_1$ satisfies the zero boundary trace conditions, $v_1 \in H_0^1(0,1) = H_{0,\mu}^1(0,1)$.\\

Assume now that $\mu = 0$. Let $\chi \in C^\infty([0,1])$ be such that $\chi(r) = 1$ for every $r \in [0,z_{1,v}/4]$, $\chi(r) = 0$ for every $r \in [z_{1,v}/2,1]$, and $\chi(r) \in [0,1]$. We write 
\begin{align*}
    v_1 = \chi v_1 + (1-\chi)v_1.
\end{align*}
Since $1-\chi$ is supported away from zero, it satisfies $(1-\chi)v_1 \in H_0^1(0,1) \subset H_{0,0}^1(0,1)$. On the other hand, 
\begin{align*}
    \chi v_1 = 2\xi r \chi \varphi' + \chi \left( 2\xi^2r^2  +\xi  - \lambda_{\xi,\mu,k+1} \right) \varphi.
\end{align*}
Since $H_{0,0}^1 (0,1)$ is invariant under multiplication by smooth functions, we have 
\begin{align*}
    \chi  \left( 2\xi^2r^2  +\xi  - \lambda_{\xi,\mu,k+1} \right) \varphi \in H_{0,0}^1(0,1).
\end{align*}
 It remains to treat the term $r \chi \varphi'$. Using the known expression of the eigenfunctions from Proposition \ref{prop: exact form eigenfunctions singular + cond eigenv}, we have $r \chi \varphi' \in L^2 \cap C^\infty(0,1)$ and $\lim_{r \to 0^+} r \chi(r) \varphi'(r) = 0$. Differentiating, we obtain 
\begin{align*}
    (r \chi \varphi')' = \chi' r \varphi ' + \chi r \varphi'' + \chi \varphi' =  \chi' r \varphi ' + \chi r \left(\xi^2r^2 - \frac{1}{4r^2} - \lambda_{\xi,\mu,k+1} \right)\varphi + \chi \varphi'.
\end{align*}
Setting $f = r \chi \varphi' - \chi \varphi/2 \in L^2 \cap C^\infty(0,1)$, which satisfies $f(0)= f(1) = 0$, we obtain
\begin{align*}
    f' = \chi' \left( r \varphi ' - \frac{\varphi}{2} \right) + \chi \left( \frac{\varphi'}{2} - \frac{1}{4r} \varphi \right) + \chi r \left(\xi^2r^2 - \lambda_{\xi,\mu,k+1} \right)\varphi.
\end{align*}
Now, using the known explicit form of the eigenfunctions $\varphi$ given in Proposition \ref{prop: exact form eigenfunctions singular + cond eigenv}, each summand in the above expression belongs to $L^2(0,1)$, and thus $f \in H_0^1(0,1) \subset H_{0,0}^1(0,1)$. Hence, $r\chi \varphi' = f + \chi \varphi/2 \in H_{0,0}^1(0,1)$, and this concludes the proof. 
\end{proof}

\subsection{Uniform localization of eigenfunctions}\label{section: Uniform localization of eigenfunctions}

In this section, we always assume $\gamma = 1$, and we use the notation $\partial_r f = f'$. We present some uniform localization estimates for the eigenfunctions of $G_{\xi,\mu}$ introduced in \eqref{eqn: def G radial}. Our proofs are direct adaptations of \cite{laurent2023uniform} (see also, \textit{e.g.} \cite{allibert1998controle}). We also mention that \cite{dietze2024concentration, dietze2025critical} study the concentration phenomena of eigenfunctions on singular Riemannian manifolds. 

\begin{theorem}\label{thm: localization estimates eigenfunctions}
    Let $\gamma = 1$. Let $(a,b) \subset (0,1)$. For every $\epsilon \in (0, (b-a)/2)$, there exists a constant $C > 0$ such that: for every $(\xi, \mu,k) \in (0,+\infty) \times [0,+\infty) \times \mathbb{N}$, for every pair $(\varphi,\lambda_{\xi,\mu,k}) \in D(G_{\xi,\mu}) \times \mathbb{R}$ such that $G_{\xi,\mu} \varphi = \lambda_{\xi,\mu,k} \varphi$, with $\|\varphi\|_{L^2(0,1)} = 1$, we have
    \begin{align}\label{eqn: uniform localization estimates}
        \begin{aligned}
            &\|\varphi\|_{L^2(a,b)}^2 \geq  \frac{C}{\lambda_{\xi,\mu,k}^2} \exp \left( -2 \sqrt{\lambda_{\xi,\mu,k}}\left(1 + \ln \frac{1}{a+\epsilon} \right) \right) \exp \left( - \xi (a+2\epsilon)^2 \right), \quad \text{if } \mu \geq 1/2,\\[6pt]
            &\|\varphi\|_{L^2(a,b)}^2 \geq \frac{C}{\lambda_{\xi,\mu,k}^2} \exp \left( -2\sqrt{\lambda_{\xi,\mu,k}}\left( 1  + \ln \frac{1}{a+\epsilon} + \frac{\ln \lambda_{\xi,\mu,k}}{4} \right) \right)  \exp \left( - \xi (a+2\epsilon)^2 \right), \quad \text{if } \mu \in [0,1/2).
        \end{aligned}
    \end{align}
\end{theorem}

This theorem is less precise than its counterpart \cite[Section 2.2]{laurent2023uniform}, because we lack strong uniform localization estimates in the allowed region $V - \lambda \leq 0$, where $V$ is given in \eqref{eqn: def E E+ V} below. As a consequence, we have to pay some additional exponential prices to propagate the information there, which are exactly the exponential factors involving the $\sqrt{\lambda}$ in \eqref{eqn: uniform localization estimates} above. \\

Moreover, the presence of the singular term adds technical difficulties. We need to introduce the objects below, which shall also simplify the presentation, and state several preliminary technical lemmas that will be proved afterwards in Section \ref{section: proof of the technical lemmas}.\\

Given $(\xi, \mu,k) \in (0,+\infty) \times [0,+\infty) \times \mathbb{N}$, and a pair $(\varphi,\lambda_{\xi,\mu,k}) \in D(G_{\xi,\mu}) \times \mathbb{R}$ such that $G_{\xi,\mu} \varphi = \lambda_{\xi,\mu,k} \varphi$,  we introduce: 
\begin{align}\label{eqn: def E E+ V}
    \begin{array}{ccl}
        \mathcal{E} (r) &:=&  (\varphi')^2(r) + \lambda_{\xi,\mu,k} \varphi^2(r), \\[6pt]
        \mathcal{E}_+ (r) &:=& (\varphi')^2(r) + \left( \xi^2 r^2 + \displaystyle \frac{\mu^2 - 1/4}{r^2} - \lambda_{\xi,\mu,k} \right)\varphi^2(r), \\[6pt]
        V(r) &=& V_{\xi,\mu}(r) :=  \xi^2 r^2 + \displaystyle \frac{\mu^2 - 1/4}{r^2},
    \end{array}
\end{align}
which are well-defined on $(0,1]$. \\

Let us also introduce some values of interest in $(0,1)$ for the proofs. Given $(\xi, \mu,k) \in (0,+\infty) \times [0,+\infty) \times \mathbb{N}$, a pair $(\varphi,\lambda_{\xi,\mu,k}) \in D(G_{\xi,\mu}) \times \mathbb{R}$ such that $G_{\xi,\mu} \varphi = \lambda_{\xi,\mu,k} \varphi$, and some $\alpha > 0$, we introduce:
\begin{itemize}
    \item \textit{The values $r_-(\alpha),r_+(\alpha)$.} We look for the solutions in $(0,+\infty)$ of the equation $V(r) = \lambda_{\xi,\mu,k} + \alpha^2$. Note that we always have $\lim_{r \to +\infty} V(r) = +\infty$. \\
    
    On the one hand, since $\lim_{r \to 0^+} V(r) = -\infty$ if $\mu < 1/2$ , and $V(0)=0$ if $\mu = 1/2$, this equation has only one solution $R_0 > 0$ if $\mu \in [0,1/2]$. On the other hand, when $\mu > 1/2$, we have $\lim_{r \to 0^+} V(r) = + \infty$ and it has exactly two solutions $0 < R_1 < R_2 $, with $R_1 \in (0,1)$ since $\lambda_{\xi,\mu,k} \geq \min_{r \in (0,1]} V(r)$ by the min-max theorem. We set
    \begin{align}\label{eqn: def r-}
        r_-(\alpha) &:= \left\{
            \begin{array}{lll}
                0 &\text{if}& \mu \in [0,1/2],  \\
                R_1 &\text{if}& \mu > 1/2,
            \end{array}
        \right.\\  
        \label{eqn: def r+}
        r_+(\alpha) &:= \left\{
            \begin{array}{lll}
                \min(R_0,1) &\text{if}& \mu \in [0,1/2],  \\
                \min(R_2,1) &\text{if}& \mu > 1/2.
            \end{array}
        \right.
    \end{align}  
    Note that we always have $V'(r) \leq 0$ for every $r \in (0,r_-(\alpha)]$, $V'(r) \geq 0$ for every $r \in [r_+(\alpha),1)$, and $\{V- \lambda \leq \alpha^2\} \cap (0,1] = (r_-(\alpha), r_+(\alpha)]$.
    \item \textit{The value $r_0$.} We have $(\varphi^2)'' = 2\varphi''\varphi + 2(\varphi')^2 = 2(V-\lambda_{\xi,\mu,k})\varphi^2 + 2(\varphi')^2$. Hence, $\varphi^2$ is convex on $\{ V - \lambda_{\xi,\mu,k} \geq 0\}$. Since $\varphi(0) = \varphi(1) = 0$, and $\varphi$ is continuous on $(0,1)$, then $\varphi^2$ attains its maximum at some point $r_0 \in \{ V - \lambda_{\xi,\mu,k} \leq 0\} \subset [r_-(\alpha),r_+(\alpha)]$ for any $\alpha > 0$. Moreover, 
    \begin{align*}
        \varphi(r_0)^2 = \max_{r \in (0,1)} \varphi(r)^2 \geq \int_0^1 \varphi(r)^2 \ dr,
    \end{align*}
    from which it follows that, for $\|\varphi\|_{L^2(0,1)} = 1$,
    \begin{align}\label{eqn: estimate E(r0)}
        \mathcal{E}(r_0) \geq \lambda_{\xi,\mu,k}.
    \end{align}

\end{itemize}

We will use the following four technical lemmas, which we will prove afterward in Section \ref{section: proof of the technical lemmas}. Lemmas \ref{lemma: estimate E+} and \ref{lemma: estimate E} are respectively analogous to \cite[Lemma 2.5 and Lemma 2.6]{laurent2023uniform}.

\begin{lemma}\label{lemma: estimate E+}
    Let $(\xi, \mu,k) \in (0,+\infty) \times [0,+\infty) \times \mathbb{N}$, and a pair $(\varphi,\lambda_{\xi,\mu,k}) \in D(G_{\xi,\mu}) \times \mathbb{R}$ such that $G_{\xi,\mu} \varphi = \lambda_{\xi,\mu,k} \varphi $. For every $\alpha > 0$, for every $r,s \in (0,1)$ belonging to the same connected component of $\{ V - \lambda_{\xi,\mu,k} \geq \alpha^2\} \cap (0,1)$, we have
    \begin{align}
        \mathcal{E}_+(r) \leq \frac{1}{\alpha^2}\left( \xi^2 + \frac{(\mu^2 - 1/4)_+}{\min(r,s)^2} - \lambda_{\xi,\mu,k} \right) \exp \left( 2 \left| \int_s^r (V(\tau) - \lambda_{\xi,\mu,k})^{1/2} \ d\tau \right| \right) \mathcal{E}_+(s).
    \end{align}
\end{lemma}

\begin{lemma}\label{lemma: estimate E}
    Let $(\xi, \mu,k) \in (0,+\infty) \times [0,+\infty) \times \mathbb{N}$, and a pair $(\varphi,\lambda_{\xi,\mu,k}) \in D(G_{\xi,\mu}) \times \mathbb{R}$ such that $G_{\xi,\mu} \varphi = \lambda_{\xi,\mu,k} \varphi $. For every $r,s \in (0,1)$ we have
    \begin{align}
        \mathcal{E} (r) \leq \exp \left( \frac{1}{\sqrt{\lambda_{\xi,\mu,k}}} \left| \int_s^r \left| V(\tau) \right| \ d\tau \right|\right)\mathcal{E} (s).
    \end{align}
\end{lemma}

\begin{lemma}\label{lemma: estimates r+ r-}
    For every $(\xi, \mu,k) \in (0,+\infty) \times [0,+\infty) \times \mathbb{N}$ and pair $(\varphi,\lambda_{\xi,\mu,k}) \in D(G_{\xi,\mu}) \times \mathbb{R}$ such that $G_{\xi,\mu} \varphi = \lambda_{\xi,\mu,k} \varphi$, we have 
    \begin{align}
           \displaystyle \frac{\mu^2 - 1/4}{r_+\left( \sqrt{\lambda_{\xi,\mu,k}} \right)^2} &\leq \lambda_{\xi,\mu,k}, \label{eqn: esimate c/r+}\\[8pt]
            \xi^2 r_-\left( \sqrt{\lambda_{\xi,\mu,k}} \right)^2 &\leq \lambda_{\xi,\mu,k}, \label{eqn: estimate r-^2}
    \end{align}
    where $r_-(\cdot),r_+(\cdot)$ are defined in \eqref{eqn: def r-} and \eqref{eqn: def r+}.
\end{lemma}

\begin{lemma}\label{lemma: estimate r0}
    For every $(\xi, \mu,k) \in (0,+\infty) \times [0,+\infty) \times \mathbb{N}$ and pair $(\varphi,\lambda_{\xi,\mu,k}) \in D(G_{\xi,\mu}) \times \mathbb{R}$ such that $G_{\xi,\mu} \varphi = \lambda_{\xi,\mu,k} \varphi$ with $\|\varphi\|_{L^2(0,1)} = 1$, let $r_0 \in (0,1)$ be a point where $\varphi(r_0)^2 = \max_{s \in (0,1)}\varphi(s)^2$. Then, if $\mu \in [0,1/2)$, we have 
    \begin{align}
            r_0 \geq \displaystyle \frac{1}{2\lambda_{\xi,\mu,k} \ln(\lambda_{\xi,\mu,k})}.
    \end{align}
\end{lemma}

We can now prove Theorem \ref{thm: localization estimates eigenfunctions}.

\begin{proof}[Proof of Theorem \ref{thm: localization estimates eigenfunctions}]
    Let $(a,b) \subset (0,1)$, $\epsilon \in (0, (b-a)/2)$, $(\xi, \mu,k) \in (0,+\infty) \times [0,+\infty) \times \mathbb{N}$, and a pair $(\varphi,\lambda_{\xi,\mu,k}) \in D(G_{\xi,\mu}) \times \mathbb{R}$ such that $G_{\xi,\mu} \varphi = \lambda_{\xi,\mu,k} \varphi$, with $\|\varphi\|_{L^2(0,1)} = 1$. For simplicity of the presentation, we write $\lambda = \lambda_{\xi,\mu,k}$. We also set $\alpha = \sqrt{\lambda}$, and write $r_\pm = r_\pm(\alpha)$, where $r_\pm(\alpha)$ are defined in \eqref{eqn: def r-} and \eqref{eqn: def r+}. \\

    We divide the proof into two parts. We first prove Theorem \ref{thm: localization estimates eigenfunctions} for $\mu \geq 1/2$, and then we present the necessary modifications for $\mu \in [0,1/2)$.\\
    
    \textbf{First part $\mu \geq 1/2$.} Letting $r \in (0,1)$, we have three possible cases,
    \begin{align}
        \mathrm{(i)} \ r \in [r_-, r_+]\cap(0,1), \quad \mathrm{(ii)} \ r \in [r_+,1) \text{ if } r_+ < 1, \quad \mathrm{(iii)} \ r \in (0,r_-] \text{ if } r_- > 0.
    \end{align}
        
    \textit{First case $r \in [r_-, r_+]\cap(0,1)$}. We apply Lemma \ref{lemma: estimate E} and use \eqref{eqn: estimate E(r0)} to obtain 
    \begin{align*}
         \lambda\leq \mathcal{E} (r_0) \leq \exp \left( \frac{1}{\sqrt{\lambda}} \left| \int_{r_0}^r \left| V(\tau) \right| \ d\tau \right|\right)\mathcal{E} (r),
    \end{align*}
    where we recall that $r_0 \in (r_-,r_+)$. We estimate $|V(\tau)|$ on the segment joining $r_0$ to $r$. Since $\mu \geq 1/2$, we have $V \geq 0$. Hence, for any $r \in [r_-, r_+]$, and $\tau \in [\min(r,r_0), \max(r,r_0)] \subset [r_-,r_+]$, we have
    \begin{align*}
        \begin{array}{lll}
            |V(\tau)| = V(\tau) \leq 2\lambda.
        \end{array}
    \end{align*}
    Since the segment joining $r_0$ to $r$ is of length less than $1$, we thus obtain that 
    \begin{align}\label{eqn: estimate E}
        \mathcal{E} (r) \geq \lambda \exp \left( -2 \sqrt{\lambda} \right), \quad \text{for any } r \in [r_-,r_+] \cap(0,1).
    \end{align}
    
    \textit{Second case $r \in [r_+,1)$ if $r_+ < 1$.} We first apply Lemma \ref{lemma: estimate E+} on the segment $[r_+,r]$, with $\alpha = \sqrt{\lambda}$, to obtain that for every $r \in [r_+,1)$,
    \begin{align*}
        \mathcal{E}_+(r_+) \leq \frac{1}{\lambda}\left( \xi^2 + \frac{(\mu^2 - 1/4)_+}{\min(r_+,r)^2} - \lambda \right) \exp \left( 2 \left| \int_{r_+}^r (V(\tau) - \lambda)^{1/2} \ d\tau \right| \right) \mathcal{E}_+(r).
    \end{align*}
    Now, we have $\mathcal{E}_+(r_+) =  \mathcal{E}(r_+)$, and we can thus bound from below $\mathcal{E}_+(r_+)$ using \eqref{eqn: estimate E} to obtain that, for every $ r \in [r_+,1)$,
    \begin{align}\label{eqn: estimate E_+ brute r+}
        \lambda \exp \left(-2\sqrt{\lambda}  \right) \leq \frac{1}{\lambda}\left( \xi^2 + \frac{(\mu^2 - 1/4)_+}{\min(r_+,r)^2} - \lambda \right) \exp \left( 2 \left| \int_{r_+}^r (V(\tau) - \lambda)^{1/2} \ d\tau \right| \right) \mathcal{E}_+(r).
    \end{align}
    We estimate the term involving the $\min$ in \eqref{eqn: estimate E_+ brute r+}. We use that $r_+ \leq r$ and \eqref{eqn: esimate c/r+} in Lemma \ref{lemma: estimates r+ r-} to obtain that
    \begin{align*}
        \frac{(\mu^2 - 1/4)_+}{\min(r_+,r)^2} - \lambda \leq 0.
    \end{align*}
    Using that $\xi^2/\lambda^2 \leq 1/4$ from Proposition \ref{prop: lower bounds eigenvalue}, and substituting this into \eqref{eqn: estimate E_+ brute r+}, we obtain,
    \begin{align}\label{eqn: estimate E+ i}
         \mathcal{E}_+(r) \geq 4 \exp \left( -2 \sqrt{\lambda}  \right) \exp \left( -2 \left| \int_{r_+}^r (V(\tau) - \lambda)^{1/2} \ d\tau \right| \right), \quad \text{for every } r \in [r_+,1).
    \end{align}

    \textit{Third case $r \in (0,r_-]$ if $r_- > 0$.} We obtain symmetrically to the above reasoning that, for every $ r \in (0,r_-]$,
    \begin{align}\label{eqn: estimate E_+ brute r-}
        \lambda \exp \left( -2 \sqrt{\lambda}  \right) \leq \frac{1}{\lambda}\left( \xi^2 + \frac{(\mu^2 - 1/4)_+}{\min(r_-,r)^2} - \lambda \right) \exp \left( 2 \left| \int_{r_-}^r (V(\tau) - \lambda)^{1/2} \ d\tau \right| \right) \mathcal{E}_+(r).
    \end{align}
    We treat the $\min$ term in \eqref{eqn: estimate E_+ brute r-}. Here we have $r \leq r_-$, and $\mu^2 - 1/4 < \lambda$ by Proposition \ref{prop: lower bounds eigenvalue}, so that 
    \begin{align*}
        \frac{(\mu^2 - 1/4)_+}{\min(r_-,r)^2} - \lambda \leq \frac{\lambda}{r^2}.
    \end{align*}
    Thus, substituting the above inequality into \eqref{eqn: estimate E_+ brute r-}, using again $\xi^2/\lambda^2 \leq 1/4$ and the rough estimate $\lambda \geq 4$ from Proposition \ref{prop: lower bounds eigenvalue}, we obtain
    \begin{align}\label{eqn: estimate E+ ii}
        \mathcal{E}_+(r) \geq 2r^2\exp \left( -2 \sqrt{\lambda}   \right) \exp \left( -2 \left| \int_{r_-}^r (V(\tau) - \lambda)^{1/2} \ d\tau \right| \right), \quad \text{ for every } r \in (0,r_-].
    \end{align}
    We now unify \eqref{eqn: estimate E}, \eqref{eqn: estimate E+ i}, and \eqref{eqn: estimate E+ ii}. Define 
    \begin{align}
        F(r) &:= \mathcal{E}(r) \mathbbm{1}_{[r_-,r_+] \cap (0,1)} + \mathcal{E}_+(r) \mathbbm{1}_{(0,1) \setminus [r_-,r_+]},
    \end{align}
    and
    \begin{align}\label{eqn: def d(r)}
        d(r) := \mathbbm{1}_{r < r_-}(r) \int_r^{r_-} (V(\tau) - \lambda)^{1/2} \ d\tau + \mathbbm{1}_{r > r_+}(r) \int_{r_+}^r (V(\tau) - \lambda)^{1/2} \ d\tau. 
    \end{align}
    With these notations, combining \eqref{eqn: estimate E}, \eqref{eqn: estimate E+ i}, \eqref{eqn: estimate E+ ii}, using that $\lambda \geq 4 \geq 2r^2$ by Proposition \ref{prop: lower bounds eigenvalue}, and $1 \geq \exp \left( -2 d(r)\right)$, we obtain
    \begin{align}\label{eqn: estimate F i}
        F(r) \geq 2r^2 \exp \left( -2 \sqrt{\lambda} \right) \exp \left( -2 d(r)\right), \quad \text{ for every } r \in (0,1].
    \end{align}
    Let now $\epsilon \in (0,(b-a)/2)$, and set $I_\epsilon := (a+\epsilon, a + 2\epsilon)$. The rest of the proof consists of refining this last estimate \eqref{eqn: estimate F i} by estimating $d(r)$ from above, for $r \in I_\epsilon$, and integrating. If $r_+ < 1$, for every $r > r_+$ we have, using \eqref{eqn: esimate c/r+} in Lemma \ref{lemma: estimates r+ r-},
    \begin{align*}
        V(r) - \lambda = \xi^2r^2 + \frac{\mu^2 - 1/4}{r^2} - \lambda \leq  \xi^2r^2 + \frac{\mu^2 - 1/4}{r_+^2} - \lambda \leq \xi^2r^2
    \end{align*}
    On the other hand, if $r_- > 0$, for every $r \leq r_-$ we have, using \eqref{eqn: estimate r-^2} in Lemma \ref{lemma: estimates r+ r-} and Proposition \ref{prop: lower bounds eigenvalue},
    \begin{align*}
       V(r) - \lambda = \xi^2r^2 + \frac{\mu^2 - 1/4}{r^2} - \lambda \leq  \xi^2r_-^2 + \frac{\mu^2 - 1/4}{r^2} - \lambda \leq  \frac{\mu^2 - 1/4}{r^2} \leq \frac{\lambda}{r^2}. 
    \end{align*}
    Hence, it follows from the last two inequalities above that $d(r)$ given in \eqref{eqn: def d(r)} satisfies
    \begin{align}\label{eqn: estimate d}
        d(r) \leq \sqrt{\lambda}\ln \frac{1}{a+\epsilon} + \frac{\xi (a+2\epsilon)^2}{2}, \quad \text{for every } r \in I_\epsilon. 
    \end{align}
    Thus, it follows from \eqref{eqn: estimate F i} and \eqref{eqn: estimate d} that
    \begin{align}\label{eqn: estimate F ii}
        F(r) \geq 2(a+\epsilon)^2 \exp \left( -2 \sqrt{\lambda}\left(1 + \ln \frac{1}{a+\epsilon} \right) \right) \exp \left( - \xi (a+2\epsilon)^2 \right), \quad \text{for every } r \in I_\epsilon.
    \end{align}
    Integrating \eqref{eqn: estimate F ii} on $I_\epsilon$ we obtain, for some $C = C(a,\epsilon) > 0$,
    \begin{align}\label{eqn: estimate F iii}
        \int_{a+\epsilon}^{a+2\epsilon} F(r) \ dr \geq C(a,\epsilon) \exp \left( -2 \sqrt{\lambda}\left(1 + \ln \frac{1}{a+\epsilon} \right) \right) \exp \left( - \xi (a+2\epsilon)^2 \right).
    \end{align}
    Now, on $I_\epsilon$, using Proposition \ref{prop: lower bounds eigenvalue}, for some $C_0 > 0$ that only depends on $a,\epsilon$, we have $\xi^2r^2 \leq \xi^2(a+2\epsilon)^2 \leq C_0 \lambda^2$ and $(\mu^2 - 1/4)/r^2 \leq \lambda/(a+\epsilon)^2 \leq C_0 \lambda^2$. Thus, 
    \begin{align}\label{eqn: estimate F iv}
        \int_{a+\epsilon}^{a+2\epsilon} F(r) \ dr \leq \int_{a+\epsilon}^{a+2\epsilon} |\varphi'|^2(s) \ ds  + C_0 \lambda^2 \int_{a+\epsilon}^{a+2\epsilon} |\varphi|^2(s) \ ds. 
    \end{align}
    We bound from above the integral involving $\varphi'$ in \eqref{eqn: estimate F iv}. Let $\chi \in C_c^\infty(a,b)$, $\chi(x) \in [0,1]$, $\chi = 1$ on $(a+\epsilon, a+ 2\epsilon)$, and $\chi(x) = 0$ on $(a,a+\epsilon/2) \cup (\epsilon + (a+b)/2, b)$. We have 
    \begin{align*}
        \int_{a+\epsilon}^{a+2\epsilon} |\varphi'|^2(s) \ ds \leq \int_{a+\epsilon/2}^{b} \chi(s) |\varphi'|^2(s) \ ds.
    \end{align*}
    Using $(\varphi ')^2 = (\varphi \varphi')' - \varphi \varphi'' = (\varphi \varphi')' - (V - \lambda)\varphi^2$, and integrating by parts, we obtain 
    \begin{align}\label{eqn: estimate F v}
        \begin{aligned}
            \int_{a+\epsilon}^{a+2\epsilon} F(r) \ dr \leq \int_{a+\epsilon/2}^b \frac{|\chi''(s)|}{2}\varphi(s)^2 \ ds + \int_{a+\epsilon/2}^{b} |V(s) - \lambda| \varphi(s)^2 \ ds  +  C_0 \lambda^2 \int_a^b |\varphi|^2(s) \ ds. 
        \end{aligned}
    \end{align}
    Using Proposition \ref{prop: lower bounds eigenvalue}, for $ r \in (a+\epsilon/2,b)$, for some $C_1 > 0$ that does not depend on $\xi, \mu$, we have $|V-\lambda| \leq \xi^2 + \mu^2/(a+\epsilon/2)^2 + \lambda \leq C_1 \lambda^2$. Substituting this into \eqref{eqn: estimate F v}, there exists some $C_2$ that only depends on $\chi, \epsilon, a, b$, such that 
    \begin{align}\label{eqn: estimate F vi}
        \int_{a+\epsilon}^{a+2\epsilon} F(r) &\ dr \leq C_2 \lambda^2 \int_a^b \varphi(s)^2 \ ds.
    \end{align}
    Combining \eqref{eqn: estimate F iii} with \eqref{eqn: estimate F vi} concludes the proof for $\mu \geq 1/2$.\\
    
    \textbf{Second part $\mu \in [0,1/2)$.} We only sketch the proof here, as it is very similar to the above treatment except for the first step due to the sign of $V$. In this case, we necessarily have $r_- = 0$ by definition \eqref{eqn: def r-}. We thus only have two cases for $r \in (0,1)$, given by 
    \begin{align}
        \mathrm{(i)} \ r \in (0, r_+] \cap (0,1), \quad \mathrm{(ii)} \ r \in [r_+,1), \text{ if } r_+ < 1. 
    \end{align}
    
    \textit{First case $r \in (0, r_+] \cap (0,1)$}. We apply Lemma \ref{lemma: estimate E} and use \eqref{eqn: estimate E(r0)} to obtain 
    \begin{align}\label{eqn: estimate part 2}
         \lambda\leq \mathcal{E} (r_0) \leq \exp \left( \frac{1}{\sqrt{\lambda}} \left| \int_{r_0}^r \left| V(\tau) \right| \ d\tau \right| \right)\mathcal{E} (r),
    \end{align}
    where we recall that $r_0 \in (0,r_+)$. We estimate $|V(\tau)|$ on the segment joining $r_0$ to $r$. Since $\mu < 1/2$, we may have $V \leq 0 $. Hence, for any $r \in (0, r_+]$, and $\tau \in [\min(r,r_0), \max(r,r_0)] \subset (0,r_+]$, we have
    \begin{align*}
        \begin{array}{lll}
            |V(\tau)| = V(\tau) &\leq 2\lambda, & \text{if } V \geq 0,\\[6pt]
            |V(\tau)| = \displaystyle \frac{1/4 - \mu^2}{\tau^2} - \xi^2 \tau^2 &\leq \displaystyle \frac{1}{4\tau^2},& \text{if } V(\tau) \leq 0.
        \end{array}
    \end{align*}
    Thus, integrating,  
    \begin{align*}
        \left| \int_{r_0}^r \left| V(\tau) \right| \ d\tau \right|  &\leq  2\lambda + \int_{\min(r,r_0)}^1 \frac{d\tau}{4\tau^2} \leq 2\lambda + \frac{1}{4\min(r,r_0)} \leq 2\lambda + \frac{1}{4r} + \frac{1}{4r_0}.
    \end{align*}
    Using Lemma \ref{lemma: estimate r0}, we thus obtain 
    \begin{align}
        \left| \int_{r_0}^r \left| V(\tau) \right| \ d\tau \right| \leq 2\lambda + \frac{1}{4r} + \frac{\lambda \ln \lambda}{2}.
    \end{align}
    Substituting this into \eqref{eqn: estimate part 2}, we arrive at 
    \begin{align}\label{eqn: estimate part 2 E 0}
        \mathcal{E}(r) \geq  \lambda \exp \left( - \frac{1}{4r\sqrt{\lambda}}\right) \exp \left( -2\sqrt{\lambda}\left( 1 + \frac{\ln \lambda}{4}\right)  \right).
    \end{align}
    Using the rough estimate $\lambda \geq 4$ from Proposition \ref{prop: lower bounds eigenvalue}, we thus finally have 
    \begin{align}\label{eqn: estimate part 2 E i}
        \mathcal{E}(r) \geq  \lambda \exp \left( - \frac{1}{8r}\right) \exp \left( -2\sqrt{\lambda}\left( 1 + \frac{\ln \lambda}{4}\right)  \right).
    \end{align}
    Note that since $r$ will be taken in $(a+\epsilon, a+2\epsilon)$, the exponential involving $r$ will only become a nonnegative constant that depends only on $a,\epsilon$.\\
    
    \textit{Second case $r \in [r_+,1)$, if $r_+ < 1$}. In this case, we follow exactly the steps of the second case of the first part of this proof. We obtain, for every $ r \in [r_+,1)$,
    \begin{align*}
         \mathcal{E}_+(r) \geq  \frac{\lambda}{\xi^2}\exp \left( -2 \left| \int_{r_+}^r (V(\tau) - \lambda)^{1/2} \ d\tau \right| \right) \mathcal{E}_+(r_+) = \frac{\lambda}{\xi^2}\exp \left( -2 \left| \int_{r_+}^r (V(\tau) - \lambda)^{1/2} \ d\tau \right| \right) \mathcal{E}(r_+).
    \end{align*}
    Here we use \eqref{eqn: estimate part 2 E 0} instead of \eqref{eqn: estimate part 2 E i}, to obtain, since $\lambda^2/\xi^2 \geq 4$ by Proposition \ref{prop: lower bounds eigenvalue}, 
    \begin{align}\label{eqn: estimate part 2 E+ i}
         \mathcal{E}_+(r) \geq 4 \exp \left( - \frac{1}{4r_+ \sqrt{\lambda}}\right) \exp \left( -2\sqrt{\lambda}\left( 1 + \frac{\ln \lambda}{4}\right) \right) \exp \left( -2 \left| \int_{r_+}^r (V(\tau) - \lambda)^{1/2} \ d\tau \right| \right).
    \end{align}
    Using that $V(r_+) = 2\lambda = \xi^2 r_+^2 + (\mu^2 - 1/4)/r_+^2 \leq \xi^2 r_+^2$, we obtain that $r_+^2 \geq 2\lambda/\xi^2$. Now, since $\lambda \geq 2\xi$ by Proposition \ref{prop: lower bounds eigenvalue}, we obtain from \eqref{eqn: estimate part 2 E+ i} that 
    \begin{align}\label{eqn: estimate part 2 E+ ii}
         \mathcal{E}_+(r) \geq 4 \exp \left( - \frac{1}{8}\right) \exp \left( -2\sqrt{\lambda}\left( 1 + \frac{\ln \lambda}{4}\right) \right) \exp \left( -2 \left| \int_{r_+}^r (V(\tau) - \lambda)^{1/2} \ d\tau \right| \right).
    \end{align}
    
    Following the end of the proof of the first part, we unify \eqref{eqn: estimate part 2 E i} and \eqref{eqn: estimate part 2 E+ ii}, bound from above $d(r)$, integrate on $(a+\epsilon,a+2\epsilon)$, and use a Caccioppoli-type argument. This concludes the proof. 
\end{proof}

\subsubsection{Proof of the technical lemmas}\label{section: proof of the technical lemmas}

\begin{proof}[\textbf{Proof of Lemma \ref{lemma: estimate E+}}]
    Let us write for simplicity $\lambda := \lambda_{\xi,\mu,k}$. We differentiate $\mathcal{E}_+$, and we obtain, for any $r\in \{ V - \lambda \geq \alpha^2\} \cap (0,1)$,
    \begin{align*}
      \mathcal{E}_+'(r) &= 2\varphi''(r)\varphi'(r) +  \left( V(r) - \lambda\right) 2\varphi'(r)\varphi(r) + V'(r)\varphi^2(r) = 4 (V(r) - \lambda)\varphi'(r)\varphi(r) + V'(r) \varphi^2(r). 
    \end{align*}
    Hence, for any $r\in \{ V - \lambda \geq \alpha^2\} \cap (0,1)$, applying Young's inequality, 
    \begin{align*}
        \left| \mathcal{E}_+'(r) \right| \leq \left( 2 (V(r) - \lambda)^{1/2}  + \frac{|V'|(r)}{V(r)-\lambda} \right) \mathcal{E}_+(r).
    \end{align*}
    Thus, letting $r,s$ belong to the same connected component of $\{V - \lambda \geq \alpha^2\} \cap (0,1)$, and applying Gronwall's lemma, we obtain
    \begin{align}\label{eqn: estimate lemma E+ i}
        \mathcal{E}_+(r) \leq \exp \left( \left| \int_s^r \left( 2 (V(\tau) - \lambda)^{1/2} + \frac{|V'(\tau)|}{V(\tau)-\lambda} \right) \ d\tau \right| \right) \mathcal{E}_+(s).
    \end{align}
    We now compute the integral involving $V'$. Recall $0 \leq r_-(\alpha) < r_+(\alpha) \leq 1$ given in \eqref{eqn: def r-} and \eqref{eqn: def r+}. With these notations we can write $\{V - \lambda \geq \alpha^2\} \cap (0,1) = (0,r_-(\alpha)] \cup [r_+(\alpha),1)$. Recall that we have $V'(r) < 0$ on $(0,r_-(\alpha)]$ and $V'(r) > 0$ on $[r_+(\alpha),1)$. Thus, we compute for $r,s$ in the same connected component of $(0,r_-(\alpha)] \cup [r_+(\alpha),1)$,
    \begin{align}\label{eqn: estimate integral V'}
        \left| \int_s^r \frac{|V'(\tau)|}{V(\tau)-\lambda} \ d\tau \right|  = \left| \int_s^r \frac{V'(\tau)}{V(\tau)-\lambda} \ d\tau \right| = \left| \ln \left( \frac{V(r) - \lambda}{V(s) - \lambda}\right) \right|.
    \end{align}
    
    Therefore, combining \eqref{eqn: estimate lemma E+ i} and \eqref{eqn: estimate integral V'}, we obtain, for every $r,s$ in the same connected component of $\{V - \lambda \geq \alpha^2\} \cap (0,1)$,
    \begin{align}\label{eqn: estimate lemma E+ ii}
        \mathcal{E}_+(r) \leq \max \left( \frac{V(r)-\lambda}{V(s)-\lambda}, \frac{V(s)-\lambda}{V(r)-\lambda} \right) \exp \left( 2 \left| \int_s^r (V(\tau) - \lambda)^{1/2} \ d\tau \right| \right) \mathcal{E}_+(s).
    \end{align}
    We have, 
    \begin{align*}
        \max \left( \frac{V(r)-\lambda}{V(s)-\lambda}, \frac{V(s)-\lambda}{V(r)-\lambda} \right) \leq \frac{\max(V(r)-\lambda, V(s)-\lambda) }{\alpha^2} = \frac{\max(V(r), V(s)) - \lambda}{\alpha^2}. 
    \end{align*}
    Now, for any $r \in (0,1)$, $V(r) \leq \xi^2 + (\mu^2 - 1/4)_+/r^2$, and we thus obtain the sought estimate by substituting this upper bound into \eqref{eqn: estimate lemma E+ ii}.
\end{proof}

\begin{proof}[\textbf{Proof of Lemma \ref{lemma: estimate E}}]
    Let us write for simplicity $\lambda := \lambda_{\xi,\mu,k}$. Differentiating $\mathcal{E}$ we obtain
    \begin{align*}
        \mathcal{E}'(r) = 2\varphi''(r)\varphi'(r) + 2\lambda \varphi'(r)\varphi(r) = 2V(r)\varphi'(r)\varphi(r).
    \end{align*}
    Hence, applying Young's inequality,
    \begin{align*}
        \left| \mathcal{E}'(r) \right| \leq \frac{\left|V(r) \right|}{\sqrt{\lambda}} \mathcal{E}(r).
    \end{align*}
    Applying Gronwall's lemma, it concludes the proof. 
\end{proof}

\begin{proof}[\textbf{Proof of Lemma \ref{lemma: estimates r+ r-}}]
    Let us write for simplicity $\lambda := \lambda_{\xi,\mu,k}$, and $r_\pm = r_\pm\left(\sqrt{\lambda} \right)$. Observe first that if $\mu \in [0,1/2]$, then $r_- = 0$ and $\mu^2-1/4 \leq 0$, and thus \eqref{eqn: esimate c/r+} and \eqref{eqn: estimate r-^2} are true. \\

    Assume now that $\mu > 1/2$. We have $V'(r) = 2\xi^2r - 2 \displaystyle \frac{\mu^2-1/4}{r^3}$, so that
    \begin{align}\label{eqn: comparison harmonic singular V' > 0 i}
        V'(r) \geq 0 \Longrightarrow \xi^2 r^2 \geq \frac{\mu^2 - 1/4}{r^2},\\
        V'(r) \leq 0 \Longrightarrow \xi^2 r^2 \leq \frac{\mu^2 - 1/4}{r^2}. \label{eqn: comparison harmonic singular V' > 0 ii}
    \end{align}
    Suppose first that $0 < r_- < r_+ < 1$ (recall that $r_- > 0$ is forced here since $\mu > 1/2$). In this case, we have $V'(r) \geq 0$ on $[r_+,1)$, and it follows from \eqref{eqn: comparison harmonic singular V' > 0 i} that $2\lambda = V(r_+) \geq 2(\mu^2-1/4)/r_+^2$, which proves \eqref{eqn: esimate c/r+}. On the other hand, we have that $V'(r) \leq 0$ on $(0,r_-]$, and it follows from \eqref{eqn: comparison harmonic singular V' > 0 ii} that $2\lambda = V(r_-) \geq 2\xi^2r_-^2$, which proves \eqref{eqn: estimate r-^2}.  If $r_+ = 1$, then \eqref{eqn: esimate c/r+} holds by Proposition \ref{prop: lower bounds eigenvalue}. This concludes the proof. 
\end{proof}

\begin{proof}[\textbf{Proof of Lemma \ref{lemma: estimate r0}}]
    Let us write for simplicity $\lambda := \lambda_{\xi,\mu,k}$. Recall that $\mathcal{Q}_{\xi,\mu}$ is defined in \eqref{eqn: def energy form Q}. Setting $\varphi = \sqrt{r}\psi$, by definition of $\mathcal{Q}_{\xi,\mu}$ and because $\|\varphi\|_{L^2(0,1)} = 1$, we obtain 
    \begin{align}
        \lambda = \mathcal{Q}_{\xi,\mu}(\varphi) \geq  \int_{r_0}^1 r\psi'(r)^2 \ dr.
    \end{align}
    Using that $\psi(1) = 0$ and the above inequality, an application of Cauchy-Schwarz gives, 
    \begin{align}
        1 \leq \varphi(r_0)^2 = r_0 |\psi(r_0)|^2 \leq r_0  \left| \int_{r_0}^1 \psi'(r) \ dr \right|^2 \leq r_0 \left( \int_{r_0}^1 r\psi'(r)^2 \ dr \right) \left( \int_{r_0}^1 \frac{dr}{r} \right) \leq \lambda r_0 \ln \frac{1}{r_0}.
    \end{align}
    We are therefore interested in the inequality
    \begin{align}
        r_0\ln \frac{1}{r_0} \geq \frac{1}{\lambda}.
    \end{align}
    Set $t = 1/r_0$ and $f(t) = \ln(t)/t$. We have that $f'(t) \geq 0$ for every $t \in (0,e]$, and $f'(t) \leq 0$ for every $t \geq e$. Moreover, $f(e) = e^{-1} > 1/\lambda$. Indeed, Proposition \ref{prop: lower bounds eigenvalue} gives $\lambda \geq 4 > e$. Thus, the set of solutions of 
    \begin{align*}
        \frac{\ln t}{t} \geq \frac{1}{\lambda}
    \end{align*}
    is an interval $[t_0,t_1]$, with $0< t_0 < e < t_1$. To bound from below $r_0$, it thus suffices to bound from above $t_1$. That is, to look for some value $t_2 \geq t_1$ such that $f(t_2) \leq 1/\lambda$. \\
   
    We set $t_2 = 2\lambda \ln(\lambda)$. We have $0 < 2\ln(\lambda) \leq \lambda$ for every $\lambda > 1$. This is equivalent to $2\lambda \ln(\lambda) \leq \lambda^2$ for every $\lambda > 0$. Since $\ln$ is increasing on $(0,+\infty)$, we apply it to both sides and we arrive at $\ln(2\lambda\ln(\lambda)) \leq 2\ln(\lambda)$. Thus, it is clear that
    \begin{align*}
        \frac{\ln(t_2)}{t_2} = \frac{\ln(2\lambda\ln(\lambda))}{2\lambda \ln(\lambda)} \leq \frac{1}{\lambda}.
    \end{align*}
    This concludes the proof since $r_0 \geq 1/t_1 \geq 1/t_2 \geq 1/(2\lambda \ln(\lambda))$.
\end{proof}

\section{Radial observability of waves: proof of Theorem \ref{thm: radial wave observability}}\label{section: radial observability of waves}

The goal of this section is to prove Theorem \ref{thm: radial wave observability}. We shall use the following results on the observability of waves, which we will prove afterward in Section \ref{section: proof technical results waves}. We may indicate the derivative with respect to a variable as a subscript, for instance $\partial_s f = f_s$, $\partial_s^2 f = f_{ss}$, $\partial_r \partial_s f = f_{rs}$, to simplify the presentation.

\begin{proposition}\label{prop: observability near singularity VZ}
    Let $\gamma \geq 0$ and $\nu \geq 0$. There exist $c, C > 0$ such that, for every $\xi \geq 0$, $S > 0$ and $R \in (0,1/2)$ satisfying
    \begin{align*}
        R \leq \frac{S}{5} \quad \text{and} \quad \xi^2 R^{2\gamma+2} \leq c,
    \end{align*}
    and for every $(w_0,w_1) \in H_{0,\nu}^1(B(0,1)) \times L^2(B(0,1))$, the solution $w$ of 
    \begin{align}\label{stm: fourier adjoint wave system}
    \left\{ \begin{array}{lcll}
        \partial_s^2 w - \Delta_x w + \left( \xi^2|x|^{2\gamma} + \displaystyle \frac{\nu^2 - \mathsf{H}}{|x|^2} \right) w &=& 0, & (s,x) \in (-S,S) \times B(0,1), \\
        w(s,x) &=& 0, & (s,x) \in (-S,S) \times \partial B(0,1), \\
        (w,w_s)(0,x) &=& (w_0,w_1)(x), & x \in B(0,1),
    \end{array} \right.
    \end{align}
    satisfies
    \begin{align}\label{eqn: observability near singularity final}
        \int_{B(0,R/2)} |w_1(x)|^2 \ dx \leq \frac{C}{R^3} \int_{-R}^{5R} \int_{\{R/4 < |x| < 2R\}} \left( w_s(s,x)^2 + w(s,x)^2 \right) \ dx \ ds.
    \end{align}
\end{proposition}

As a direct corollary of Proposition \ref{prop: observability near singularity VZ}, we have the following result for the solutions of \eqref{stm: radial wave adjoint system} with initial states $(w,w_s)(0) = (0,w_1)$. 

\begin{corollary}\label{cor: observability near singularity VZ}
    Let $\gamma \geq 0$ and $\nu \geq 0$. There exist $c, C > 0$ such that, for every $\xi \geq 0$, $\ell \in \mathbb{N}$, $S > 0$ and $R \in (0,1/2)$ satisfying
    \begin{align*}
        R \leq \frac{S}{5} \quad \text{and} \quad \xi^2 R^{2\gamma+2} \leq c,
    \end{align*}
    and for every $w_1 \in L^2(0,1)$, the solution $w$ of system \eqref{stm: radial wave adjoint system} with $\mu^2 = \nu_\ell^2$ and initial state $(w,w_s)(0) = (0,w_1)$ satisfies
    \begin{align}\label{eqn: integrated near singularity 0}
        \int_{-4R}^{4R} \int_0^{R/2} w_s(s,r)^2 \ dr \ ds \leq \frac{C}{R^2} \int_{-5R}^{5R} \int_{R/4}^{2R} \left( w_s(s,r)^2 + w(s,r)^2 \right) \ dr \ ds.
    \end{align}
\end{corollary}

\begin{proof}
Let $s_0 \in [-4R,0]$. Since \eqref{stm: fourier adjoint wave system} is invariant by translation in time, we can apply Proposition \ref{prop: observability near singularity VZ} to 
\begin{align*}
    (s,x) \mapsto W(s_0+s,x) := r^{-(d_x-1)/2}w(s_0+s,r)Y_{\ell,1}(\theta), \quad x = r\theta,
\end{align*}
with initial state $(W, W_s)(s_0,\cdot)$. Using that $(s_0-R, s_0+5R) \subset (-5R,5R)$, we thus have
\begin{align*}
    \int_0^{R/2} w_s(s_0,r)^2 \ dr \leq \frac{C}{R^3} \int_{-5R}^{5R} \int_{R/4}^{2R} \left( w_s(s,r)^2 + w(s,r)^2 \right) \ dr \ ds.
\end{align*}
Integrating with respect to $s_0 \in [-4R,0]$, and using that $w_s$ is even with respect to $s$ since $w_0 = 0$, we obtain \eqref{eqn: integrated near singularity 0}.
\end{proof}

\begin{proposition}\label{prop: radial wave lateral propagation}
    Let $\gamma \geq 0$. Let $S > 1$, $0 < a < r_0 < \delta < b < 1$, and $\eta > 0$ such that $\max(1-r_0,\delta) + \eta < S$. There exists a constant $C > 0$ such that, for every $\xi > 0$, $\ell \in \mathbb{N}$, $R \in (0, \min(r_0/2, \eta/5)]$ and $(w_0,w_1) \in H_{0,\nu}^1(B(0,1)) \times L^2(B(0,1))$, the solution of system \eqref{stm: radial wave adjoint system} with $\mu^2 = \nu_\ell^2$ satisfies
    \begin{align}\label{eqn: wave lateral propagation}
        \begin{aligned}
            &\int_{-5R}^{5R} \int_{R/4}^{2R} w_s(s,r)^2 + w(s,r)^2  \ dr \ ds +  \int_{-\eta}^\eta \int_{R/4}^1 w_s(s,r)^2 \ dr \ ds \\
            &\leq \frac{C}{\delta - r_0} \left( 1 + \xi^2\delta^{2\gamma} + \frac{|\nu_\ell^2 - 1/4|}{r_0^2} \right) \exp \left( \int_{R/4}^1 \mathcal{M}(r) \ dr \right) \int_{-S}^S \int_a^{b} w_s(s,r)^2 + w(s,r)^2  \ dr \ ds,
        \end{aligned}
    \end{align}
    where
    \begin{align}
        V(r) = \xi^2r^{2\gamma} + \frac{\nu_\ell^2 - 1/4}{r^2}, \quad \tilde{V}(r) = \sqrt{1 + V(r)^2} > 0, \quad \mathcal{M}(r) = \frac{|\tilde{V}'(r)|}{\tilde{V}(r)} + 2\sqrt{\tilde{V}(r)}.
    \end{align}
\end{proposition}

\subsection{Proof of Theorem \ref{thm: radial wave observability}}

We can now proceed to the proof of Theorem \ref{thm: radial wave observability}.

\begin{proof}[Proof of Theorem \ref{thm: radial wave observability}]

Let $S \geq 2$, $0 < a < r_0 < \delta < b < 1$, and $\eta > 0$ such that $\max(1-r_0,\delta) + \eta < S$. Let also $R= R_*/2$, where $R_*$ is the largest admissible radius in both Corollary \ref{cor: observability near singularity VZ} and Proposition \ref{prop: radial wave lateral propagation}, that is 
\begin{align}\label{eqn: def R*}
    R_* := \min \left( \frac{1}{2}, \frac{S}{5}, \frac{c^{1/(2\gamma + 2)}}{\xi^{1/(\gamma + 1)}}, \frac{r_0}{2}, \frac{\eta}{5} \right). 
\end{align}

where $c > 0$ is given by Corollary \ref{cor: observability near singularity VZ}. Let $w$ be a solution of system \eqref{stm: radial wave adjoint system}, with initial state $(0,w_1)$. From Corollary \ref{cor: observability near singularity VZ}, we have
\begin{align}\label{eqn: integrated near singularity}
    \int_{-4R}^{4R} \int_0^{R/2} w_s(s,r)^2 \ dr \ ds \leq \frac{C_1}{R^2} \int_{-5R}^{5R} \int_{R/4}^{2R} \left( w_s(s,r)^2 + w(s,r)^2 \right) \ dr \ ds.
\end{align}

Using now Proposition \ref{prop: radial wave lateral propagation}, substituting \eqref{eqn: integrated near singularity} into \eqref{eqn: wave lateral propagation}, we thus have 
    \begin{align}\label{eqn: lower bound gluing i}
        \begin{aligned}
            &\frac{R^2}{C_1} \int_{-4R}^{4R} \int_0^{R/2} w_s^2(s,r) \ dr \ ds  + \int_{R/4}^1 \int_{-\eta}^\eta w_s(s,r)^2 \ ds \ dr \\
            &\leq \frac{C_2}{\delta - r_0} \left( 1 + \xi^2\delta^{2\gamma} + \frac{|\nu_\ell^2 - 1/4|}{r_0^2} \right)
            \exp \left( \int_{R/4}^1 \mathcal{M}(r) \ dr \right)  \int_{-S}^S \int_a^{b} w_s(s,r)^2 + w(s,r)^2  \ dr \ ds,
        \end{aligned}
    \end{align}

Recalling that $\eta \geq 5R$, and since we may assume $R^2/C \leq 1$, we have
\begin{align}\label{eqn: lower bound gluing ii}
    \frac{R^2}{C_1} \int_{-4R}^{4R} \int_0^{R/2} w_s^2(s,r) \ dr \ ds  + \int_{-\eta}^\eta \int_{R/4}^1 w_s(s,r)^2 \ dr \ ds \geq \frac{R^2}{C_1} \int_{-4R}^{4R} \int_0^1 w_s^2(s,r) \ dr \ ds.
\end{align}

Expanding $w_1$ in an orthonormal basis of eigenfunctions of $G_{\xi,\nu_\ell}$ (given in \eqref{eqn: def G radial} and Section \ref{section: Uniform spectral estimates}), $w_1 = \sum_k w_k \varphi_{\xi,\ell,k}$, we have 
\begin{align*}
    \int_{-4R}^{4R} \int_0^1 w_s^2(s,r) \ dr \ ds &=  \sum_{k \geq 0} |w_k|^2  \int_{-4R}^{4R} \cos^2 \left( s \sqrt{\lambda_{\xi,\ell,k}} \right) \ ds \\
    &= \sum_{k \geq 0} |w_k|^2  \left( 4R + \frac{  \sin \left( 8R \sqrt{\lambda_{\xi,\ell,k}}\right) }{2\sqrt{\lambda_{\xi,\ell,k}}} \right) \\
    &\geq 2R \sum_{k \geq 0} |w_k|^2  =  2R \|w_1\|_{L^2(0,1)}^2.
\end{align*}

Substituting this into \eqref{eqn: lower bound gluing ii}, we obtain 
\begin{align}\label{eqn: lower bound gluing ii bis}
    \frac{R^2}{C_1} \int_{-4R}^{4R} \int_0^{R/2} w_s^2(s,r) \ dr \ ds  + \int_{R/4}^1 \int_{-\eta}^\eta w_s(s,r)^2 \ ds \ dr \geq \frac{2R^3}{C_1}\|w_1\|_{L^2(0,1)}^2.
\end{align}

Substituting now \eqref{eqn: lower bound gluing ii bis} into \eqref{eqn: lower bound gluing i}, we arrive at, for some new constant $C > 0$ that does not depend on $R,\xi,\ell$,
\begin{align}\label{eqn: lower bound gluing iii}
        \begin{aligned}
            \|w_1\|_{L^2(0,1)}^2\leq \frac{C}{R^3} \left( 1 + \xi^2\delta^{2\gamma} + \frac{|\nu_\ell^2 - 1/4|}{r_0^2} \right) \exp \left( \int_{R/4}^1 \mathcal{M}(r) \ dr \right) \int_{-S}^S \int_a^{b}  w_s(s,r)^2 + w(s,r)^2  \ dr \ ds.
        \end{aligned}
    \end{align}

Let us now estimate the exponential term involving $\mathcal{M}$ in \eqref{eqn: lower bound gluing iii} above. We have
\begin{align*}
    \mathcal{M}(r) = \frac{|\tilde{V}'(r)|}{\tilde{V}(r)} + 2\sqrt{\tilde{V}(r)} &= \frac{|V'(r)V(r)|}{1+V(r)^2} + 2\sqrt{\tilde{V}(r)} \\
    &= \frac{1}{2} \left|\frac{d}{dr}\ln(1+ V(r)^2) \right| + 2(1+V(r)^2)^{1/4} \\
    &\leq \frac{1}{2}\left|\frac{d}{dr}\ln(1+ V(r)^2) \right| + 2 + 2\sqrt{|V(r)|} \\
    &\leq \frac{1}{2}\left|\frac{d}{dr}\ln(1+ V(r)^2) \right| + 2 + 2\xi r^\gamma + 2\frac{\sqrt{|\nu_\ell^2 - 1/4|}}{r}.
\end{align*}

Integrating and taking the exponential, we obtain
\begin{align}\label{upper bound integral M i}
    \exp \left( \int_{R/4}^1 \mathcal{M}(r) \ dr \right) \leq \left( 1 + \left( \xi^2 + \frac{16|\nu_\ell^2 - 1/4|}{R^2} \right)^2 \right) \exp \left( 2 + \frac{2\xi}{1+\gamma} + 2\sqrt{|\nu_\ell^2 - 1/4|} \ln\left( \frac{4}{R}\right) \right).
\end{align}

Using that $4/R = 8/R_*$, where $R_*$ is given in \eqref{eqn: def R*}, and $1/(\gamma + 1) \leq 1$, we have, for some $C_0 > 0$ that does not depend on $\xi > 0$, that $\ln(4/R) \leq C_0 \ln(2+\xi)$. Thus, 
\begin{align*}
    \exp \left( 2 + \frac{2\xi}{1+\gamma} + 2\sqrt{|\nu_\ell^2 - 1/4|} \ln\left( \frac{4}{R}\right) \right) &\leq \exp \left( 2 (1 + \xi) + 2C_0(\nu_\ell + 1/2) \ln(2+\xi) \right) \\
    &\leq \exp \left( (2+C_0) (1 + \xi) + 2C_0\nu_\ell\ln(2+\xi) \right).
\end{align*}

Substituting this into \eqref{upper bound integral M i}, we obtain 
\begin{align}\label{upper bound integral M ii}
    \exp \left( \int_{R/4}^1 \mathcal{M}(r) \ dr \right) \leq \left( 1 + \left( \xi^2 + \frac{16|\nu_\ell^2 - 1/4|}{R^2} \right)^2 \right) \exp \left( (2+C_0) (1 + \xi) + 2C_0\nu_\ell\ln(2+\xi) \right).
\end{align}

Substituting \eqref{upper bound integral M ii} into \eqref{eqn: lower bound gluing iii}, using that $R = R_*/2$, and absorbing the polynomial term into the exponential, we obtain the sought estimate. This concludes the proof. 
\end{proof}

\subsection{Proof of the technical intermediate results}\label{section: proof technical results waves}

We start by recalling some results from \cite{vancostenoble2009hardy} that we shall use as a black box. We then prove Proposition \ref{prop: observability near singularity VZ}, and then Proposition \ref{prop: radial wave lateral propagation}.

\subsubsection{Preliminaries on the exact boundary observability of waves}\label{section: preliminaries on the exact boundary observability of waves}

We shall use the following results from \cite{vancostenoble2009hardy}. For $T > 0$, $(v_0,v_1) \in H_{0,\nu}^1(B(0,1)) \times L^2(B(0,1))$ and $f \in L^1(0,T;L^2(B(0,1)))$, we consider
\begin{align}\label{stm: wave system zuazua vancostenoble}
    \left\{ \begin{array}{lcll}
        \partial_\tau^2 v - \Delta_z v + \displaystyle \frac{\nu^2 - \mathsf{H}}{|z|^2} v &=& f, & (\tau,z) \in (0,T) \times B(0,1), \\
        v(\tau,z) &=& 0, & (\tau,z) \in (0,T) \times \partial B(0,1), \\
        (v,v_\tau)(0,z) &=& (v_0,v_1)(z), & z \in B(0,1),
    \end{array} \right.
\end{align}
and the associated energy
\begin{align*}
    E_v(\tau) := \frac{1}{2} \left( \|v(\tau)\|_{\nu}^2 + \|v_\tau(\tau)\|_{L^2(B(0,1))}^2 \right), 
\end{align*}
where $\|\cdot\|_\nu$ is the norm given in \eqref{eqn: def norm nu} with $\xi = 0$. The following proposition is from \cite[Proposition 2.1 (i)]{vancostenoble2009hardy}.

\begin{proposition}\label{prop: VZ well-posedness}
    Let $T > 0$ and $\nu \geq 0$. For every $(v_0,v_1) \in H_{0,\nu}^1(B(0,1)) \times L^2(B(0,1))$ and $f \in L^1(0,T;L^2(B(0,1)))$, there exists a unique solution $v \in C([0,T];H_{0,\nu}^1(B(0,1))) \cap C^1([0,T];L^2(B(0,1)))$ of \eqref{stm: wave system zuazua vancostenoble}. Moreover, there exists $C_1 = C_1(T,\nu) > 0$ such that, for every $(v_0,v_1) \in H_{0,\nu}^1(B(0,1)) \times L^2(B(0,1))$ and $f \in L^1(0,T;L^2(B(0,1)))$, we have
    \begin{align*}
        \|(v,v_\tau)\|_{L^\infty(0,T;H_{0,\nu}^1(B(0,1)) \times L^2(B(0,1)))} \leq C_1 \|(v_0,v_1)\|_{H_{0,\nu}^1(B(0,1)) \times L^2(B(0,1))} + C_1 \|f\|_{L^1(0,T;L^2(B(0,1)))}.
    \end{align*}
\end{proposition}

The following proposition is from \cite[Proposition 2.2 (i)]{vancostenoble2009hardy}.

\begin{proposition}\label{prop: VZ hidden regularity}
    Let $T > 0$ and $\nu \geq 0$. There exists $C_2 = C_2(T,\nu) > 0$ such that, for every $(v_0,v_1) \in H_{0,\nu}^1(B(0,1)) \times L^2(B(0,1))$ and $f \in L^1(0,T;L^2(B(0,1)))$, the solution $v$ of \eqref{stm: wave system zuazua vancostenoble} satisfies
    \begin{align*}
        \int_0^T \int_{\partial B(0,1)} \left| \partial_\eta v \right|^2 \ d\sigma \ d\tau \leq C_2 \left( E_v(0) + \|f\|_{L^1(0,T;L^2(B(0,1)))}^2 \right).
    \end{align*}
\end{proposition}

The following theorem is from \cite[Theorem 3.1]{vancostenoble2009hardy}.

\begin{theorem}\label{thm: VZ observability}
    Let $\nu \geq 0$ and $T > 2$. Then, there exists $C > 0$ such that, for every $(v_0,v_1) \in H_{0,\nu}^1(B(0,1)) \times L^2(B(0,1))$, the solution $v$ of \eqref{stm: wave system zuazua vancostenoble} with $f = 0$ satisfies
    \begin{align*}
        E_v(0) \leq C \int_0^T \int_{\partial B(0,1)} \left| \partial_\eta v \right|^2 \ d\sigma \ d\tau,
    \end{align*}
    where $\partial_\eta$ stands for the outward normal derivative.
\end{theorem}

\subsubsection{Observability of waves near the singularity: proof of Proposition \ref{prop: observability near singularity VZ}}\label{section: Observability of waves near the singularity}

We employ a cutoff argument to obtain an observability inequality near the singularity.

\begin{proof}[Proof of Proposition \ref{prop: observability near singularity VZ}]
Let $(w_0,w_1) \in H_{0,\nu}^1(B(0,1)) \times L^2(B(0,1))$, and $w$ be the solution of system \eqref{stm: fourier adjoint wave system}. \\

For any $R \in (0,1)$, let $\chi \in C_c^\infty(B(0,1))$, $\chi = 1$ on $B(0,1/2)$. Setting $\chi_R(x) = \chi(x/R)$ and $w_R := \chi_R w$, $w_R$ solves 
\begin{align}\label{stm: fourier adjoint wave system cutoff}
    \left\{ \begin{array}{lcll}
        \partial_s^2 w_R - \Delta_x w_R + \left( \xi^2|x|^{2\gamma} + \displaystyle \frac{\nu^2 - \mathsf{H}}{|x|^2} \right) w_R &=& - w\Delta\chi_R - 2\nabla \chi_R \cdot \nabla w, & (-S,S) \times B(0,R), \\
        w_R(s,x) &=& 0, & (-S,S) \times \partial B(0,R), \\
        (w_R, w_{R,s})(0,x) &=& (\chi_R w_0, \chi_R w_1)(x), & B(0,R),
    \end{array} \right.
\end{align}

We make the change of scale $\mathcal{W}(\tau, z) = w_R(R\tau, Rz)$, $(\tau, z) \in (-S/R,S/R) \times B(0,1)$. Setting $\tilde{w}(\tau,z) = w(R\tau, Rz)$, $\tilde{w}_0(z) = w_0(Rz)$ and $\tilde{w}_1(z) = w_1(Rz)$, $\mathcal{W}$ solves
\begin{align}\label{stm: fourier adjoint wave system rescaled}
    \left\{ \begin{array}{lcl}
        \partial_\tau^2 \mathcal{W} - \Delta_z \mathcal{W} +  \displaystyle \frac{\nu^2 - \mathsf{H}}{|z|^2} \mathcal{W} &=& -\tilde{w}\Delta\chi - 2\nabla \chi \cdot \nabla \tilde{w} - \xi^2R^{2\gamma+2}|z|^{2\gamma} \mathcal{W},  \\
        \mathcal{W}(\tau,z) &=& 0, \quad (\tau,z) \in (-S/R,S/R) \times \partial B(0,1), \\
        (\mathcal{W}, \mathcal{W}_\tau)(0,z) &=& (\chi \tilde{w}_0, R\chi \tilde{w}_1)(z).
    \end{array} \right.
\end{align}

We denote by $\mathcal{R}$ the solution of system \eqref{stm: fourier adjoint wave system rescaled} with initial state $(\mathcal{R}, \mathcal{R}_\tau)(0,z) = (0, 0)$, and by $\mathcal{W}^h$ the solution of system \eqref{stm: fourier adjoint wave system rescaled} with initial state $(\mathcal{W}^h, \mathcal{W}^h_\tau)(0,z) = (\chi \tilde{w}_0, R\chi \tilde{w}_1)(z)$, and with no source term. In particular $\mathcal{W} = \mathcal{W}^h + \mathcal{R}$. Applying Theorem \ref{thm: VZ observability} to $\mathcal{W}^h$ on $(0,4)$, we obtain 
\begin{align*}
    \frac{1}{2}E_{\mathcal{W}^h}(0) \leq C \int_0^4 \int_{\partial B(0,1)} \left| \partial_\eta \mathcal{W}^h \right|^2 \ d\sigma \ d\tau.
\end{align*}
Since $\partial_\eta \mathcal{W} = 0$, we have $\partial_\eta \mathcal{W}^h = - \partial_\eta \mathcal{R}$, and we obtain 
\begin{align}\label{eqn: upper bound energy Wh}
    E_{\mathcal{W}^h}(0) \leq 2C \int_0^4 \int_{\partial B(0,1)} \left| \partial_\eta \mathcal{R} \right|^2 \ d\sigma \ d\tau.
\end{align}

What follows aims at bounding from above the right-hand side of \eqref{eqn: upper bound energy Wh} using the results of \cite{vancostenoble2009hardy} that we recalled in Section \ref{section: preliminaries on the exact boundary observability of waves}. We use the notation $C > 0$ for a constant that does not depend on $R$ and $\xi$, and that may change from line to line. \\

Using Proposition \ref{prop: VZ hidden regularity}, we have, using that $E_\mathcal{R}(0) = 0$,
    \begin{align}\label{eqn: upper bound normal derivative R}
        \begin{aligned}
            &\int_0^4 \int_{\partial B(0,1)} \left| \partial_\eta \mathcal{R} \right|^2 \ d\sigma \ d\tau \leq C \left\|\tilde{w}\Delta\chi + 2\nabla \chi \cdot \nabla \tilde{w} + \xi^2R^{2\gamma+2}|z|^{2\gamma} \mathcal{W} \right\|_{L^1(0,4;L^2(B(0,1)))}^2.
        \end{aligned}
    \end{align}

We first have, by the triangle inequality, and using $\mathcal{W} = \mathcal{W}^h + \mathcal{R}$, 
\begin{align}\label{eqn: upper bound source term wave}
    \begin{aligned}
        &\left\|\tilde{w}\Delta\chi + 2\nabla \chi \cdot \nabla \tilde{w} + \xi^2R^{2\gamma+2}|z|^{2\gamma} \mathcal{W} \right\|_{L^1(0,4;L^2(B(0,1)))}^2  \\
        &\leq 2\left\| \tilde{w}\Delta\chi + 2\nabla \chi \cdot \nabla \tilde{w} \right\|_{L^1(0,4;L^2(B(0,1)))}^2 + 4\xi^4R^{4\gamma+4} \left( \left\|  \mathcal{W}^h \right\|_{L^1(0,4;L^2(B(0,1)))}^2 + \left\| \mathcal{R} \right\|_{L^1(0,4;L^2(B(0,1)))}^2 \right) 
    \end{aligned}
\end{align}

On the one hand we have, by conservation of the energy $E_{\mathcal{W}^h}(\tau) = E_{\mathcal{W}^h}(0)$, 
\begin{align}\label{eqn: upper bound Wh}
    \left\|  \mathcal{W}^h \right\|_{L^1(0,4;L^2(B(0,1)))}^2 = \left(\int_0^4 \left\|\mathcal{W}^h (\tau) \right\|_{L^2(B(0,1))} \ d\tau \right)^2 \leq C E_{\mathcal{W}^h}(0).
\end{align}

On the other hand, we have
\begin{align*}
    \left\| \mathcal{R} \right\|_{L^1(0,4;L^2(B(0,1)))} &= \int_0^4 \|\mathcal{R}(\tau)\|_{L^2(B(0,1))} \ d\tau \leq C \int_0^4 \|\mathcal{R}(\tau)\|_{H_{0,\nu}^1(B(0,1))} d \tau \\
    &\leq 4 C \sup_{\tau \in [0,4]} \left( \|\mathcal{R}(\tau)\|_{H_{0,\nu}^1(B(0,1))}^2 + \|\mathcal{R}_\tau(\tau)\|_{L^2(B(0,1))}^2 \right)^{1/2}. 
\end{align*}
Applying Proposition \ref{prop: VZ well-posedness}, we thus deduce,
\begin{align}\label{eqn: upper bound R}
    \left\| \mathcal{R} \right\|_{L^1(0,4;L^2(B(0,1)))}^2 \leq C \left\|\tilde{w}\Delta\chi + 2\nabla \chi \cdot \nabla \tilde{w} + \xi^2R^{2\gamma+2}|z|^{2\gamma} \mathcal{W} \right\|_{L^1(0,4;L^2(B(0,1)))}^2.
\end{align}

Hence, substituting \eqref{eqn: upper bound Wh} and \eqref{eqn: upper bound R} into \eqref{eqn: upper bound source term wave}, we obtain
\begin{align*}
    \begin{aligned}
        &\left\|\tilde{w}\Delta\chi + 2\nabla \chi \cdot \nabla \tilde{w} + \xi^2R^{2\gamma+2}|z|^{2\gamma} \mathcal{W} \right\|_{L^1(0,4;L^2(B(0,1)))}^2 \leq \left\| \tilde{w}\Delta\chi + 2\nabla \chi \cdot \nabla \tilde{w} \right\|_{L^1(0,4;L^2(B(0,1)))}^2 \\
        &+ C\xi^4R^{4\gamma+4} E_{\mathcal{W}^h}(0) + C\xi^4R^{4\gamma+4} \left\|\tilde{w}\Delta\chi + 2\nabla \chi \cdot \nabla \tilde{w} + \xi^2R^{2\gamma+2}|z|^{2\gamma} \mathcal{W} \right\|_{L^1(0,4;L^2(B(0,1)))}^2.
    \end{aligned}
\end{align*}
Thus, taking $R > 0$ such that $C \xi^4R^{4\gamma+4} \leq 1/2$ and $R \leq S/4$, it follows that
\begin{align}\label{eqn: upper bound source term wave ii}
    \begin{aligned}
        &\left\|\tilde{w}\Delta\chi + 2\nabla \chi \cdot \nabla \tilde{w} + \xi^2R^{2\gamma+2}|z|^{2\gamma} \mathcal{W} \right\|_{L^1(0,4;L^2(B(0,1)))}^2 \leq 2\left\| \tilde{w}\Delta\chi + 2\nabla \chi \cdot \nabla \tilde{w} \right\|_{L^1(0,4;L^2(B(0,1)))}^2 \\
        &+ 2C\xi^4R^{4\gamma+4} E_{\mathcal{W}^h}(0).
    \end{aligned}
\end{align}
Substituting \eqref{eqn: upper bound source term wave ii} into \eqref{eqn: upper bound normal derivative R} we obtain 
\begin{align}\label{eqn: upper bound source term wave iii}
    \int_0^4 \int_{\partial B(0,1)} \left| \partial_\eta \mathcal{R} \right|^2 \ d\sigma \ d\tau \leq C\left\| \tilde{w}\Delta\chi + 2\nabla \chi \cdot \nabla \tilde{w} \right\|_{L^1(0,4;L^2(B(0,1)))}^2 + C\xi^4R^{4\gamma+4} E_{\mathcal{W}^h}(0).
\end{align}

Substituting now \eqref{eqn: upper bound source term wave iii} into \eqref{eqn: upper bound energy Wh}, we have 
\begin{align}
     E_{\mathcal{W}^h}(0) \leq C\left\| \tilde{w}\Delta\chi + 2\nabla \chi \cdot \nabla \tilde{w} \right\|_{L^1(0,4;L^2(B(0,1)))}^2 + C\xi^4R^{4\gamma+4} E_{\mathcal{W}^h}(0).
\end{align}

Up to taking $R$ smaller than before, since the constants do not depend on $R$, we have arrived at, for some new constant $C > 0$ independent of $R$,  
\begin{align}\label{eqn: estimate energy Wh}
    E_{\mathcal{W}^h}(0) \leq C\left\| \tilde{w}\Delta\chi + 2\nabla \chi \cdot \nabla \tilde{w} \right\|_{L^1(0,4;L^2(B(0,1)))}^2.
\end{align}

The left-hand side of \eqref{eqn: estimate energy Wh} satisfies
\begin{align}
    E_{\mathcal{W}^h}(0) \geq \frac{R^2}{2}\|\chi \tilde{w}_1\|_{L^2(B(0,1))}^2 \geq \frac{R^{2-d_x}}{2} \int_{B(0,R/2)} |w_1|^2 \ dx.
\end{align}

On the other hand, we bound from above the right-hand side of \eqref{eqn: estimate energy Wh}, that is supported in $\{1/2 \leq |z| \leq 1\}$. By Cauchy-Schwarz in time, since $\Delta\chi$ and $\nabla\chi$ are supported in $\{1/2 \leq |z| \leq 1\}$, and changing variables back to $(s,x) = (R\tau, Rz)$, we have, setting $A_R := \{R/2 \leq |x| \leq R\}$,
\begin{align*}
    &\left\| \tilde{w}\Delta\chi + 2\nabla \chi \cdot \nabla \tilde{w} \right\|_{L^1(0,4;L^2(B(0,1)))}^2 \leq C \int_0^4 \int_{\{1/2 \leq |z| \leq 1\}} |\tilde{w}|^2 + |\nabla_z \tilde{w}|^2 \ dz \ d\tau \\
    &= C R^{-1-d_x} \int_0^{4R} \int_{A_R} |w|^2 \ dx \ ds + C R^{1-d_x} \int_0^{4R} \int_{A_R} |\nabla_x w|^2 \ dx \ ds.
\end{align*}
Hence,
\begin{align}\label{eqn: observability near singularity i}
    \int_{B(0,R/2)} |w_1|^2 \ dx \leq C \left( \frac{1}{R^3} \int_0^{4R} \int_{A_R} |w|^2 \ dx \ ds + \frac{1}{R} \int_0^{4R} \int_{A_R} |\nabla_x w|^2 \ dx \ ds \right).
\end{align}

We get rid of the gradient term in the right-hand side of \eqref{eqn: observability near singularity i} using a Caccioppoli-type argument. Set $\tilde{A}_R := \{R/4 < |x| < 2R\}$, and let $\zeta \in C_c^\infty((-R,5R) \times \tilde{A}_R)$ that satisfies
\begin{align*}
    0 \leq \zeta \leq 1, \quad \zeta = 1 \text{ on } [0,4R] \times A_R, \quad |\zeta_{ss}| + |\Delta_x \zeta| \leq C/R^2.
\end{align*}
We have
\begin{align*}
    \int_0^{4R} \int_{A_R} |\nabla_x w|^2 \ dx \ ds \leq \int_{-R}^{5R} \int_{\tilde{A}_R} \zeta |\nabla_x w|^2 \ dx \ ds.
\end{align*}

Using \eqref{stm: fourier adjoint wave system}, we have the identity
\begin{align}\label{eqn: identity wave multi d}
    |\nabla_x w|^2 = \operatorname{div}_x(w\nabla_x w) - w\Delta_x w = w_s^2 - (ww_s)_s + \operatorname{div}_x(w\nabla_x w) - \left( \xi^2|x|^{2\gamma} + \frac{\nu^2 - \mathsf{H}}{|x|^2} \right) w^2.
\end{align}
Hence, integrating by parts, 
\begin{align*}
    \int_{-R}^{5R} \int_{\tilde{A}_R} \zeta |\nabla_x w|^2 \ dx \ ds
    &= \int_{-R}^{5R} \int_{\tilde{A}_R} \zeta w_s^2 + \zeta_s ww_s - w\nabla_x\zeta \cdot \nabla_x w - \zeta\left( \xi^2|x|^{2\gamma} + \frac{\nu^2 - \mathsf{H}}{|x|^2} \right) w^2 \ dx \ ds \\
    &= \int_{-R}^{5R} \int_{\tilde{A}_R} \zeta w_s^2 + \frac{1}{2}\left( \Delta_x \zeta - \zeta_{ss} \right) w^2 - \zeta \left( \xi^2|x|^{2\gamma} + \frac{\nu^2 - \mathsf{H}}{|x|^2} \right) w^2 \ dx \ ds.
\end{align*}

It thus follows, since $R < 1$,
\begin{align}
    \int_0^{4R} \int_{A_R} |\nabla_x w|^2 \ dx \ ds \leq \frac{C}{R^2} \int_{-R}^{5R} \int_{\tilde{A}_R} w_s^2 + w^2 \ dx \ ds.
\end{align}
Substituting this into \eqref{eqn: observability near singularity i}, we obtain the sought result
\begin{align}
     \int_{B(0,R/2)} |w_1|^2 \ dx  \leq \frac{C}{R^3} \int_{-R}^{5R} \int_{\tilde{A}_R} w_s^2 + w^2 \ dx \ ds.
\end{align}
\end{proof}

\subsubsection{Observability of waves away from the singularity: proof of Proposition \ref{prop: radial wave lateral propagation}}\label{section: Lateral propagation of energy of waves}

We now give an observability estimate for the solutions of \eqref{stm: radial wave adjoint system} with initial states $w_0 = 0$, away from the singularity, using lateral propagation of energy. The following proof is an adaptation of \cite[Proposition 4.5]{beauchard20152d}.

\begin{proof}[Proof of Proposition \ref{prop: radial wave lateral propagation}]
Introduce, for $R \in (0,1)$, 
\begin{align}
    \mathcal{F}(r) &= \int_{-1+r-\eta}^{1-r+\eta} w_s(s,r)^2 + w_r(s,r)^2 + \tilde{V}(r)w(s,r)^2 \ ds, \\
    \mathcal{F}_R(r) &= \int_{-r + \frac{R}{4}-\eta}^{r - \frac{R}{4} +\eta} w_s(s,r)^2 + w_r(s,r)^2 + \tilde{V}(r)w(s,r)^2 \ ds.
\end{align}

We divide the proof into five steps. The first and second steps aim at obtaining respectively the following estimates, for $0 < R/4 < r_0 < \delta < 1$,
\begin{align}
    \int_{r_0}^1  \mathcal{F}(r) \ dr &\leq \frac{1-r_0}{\delta - r_0}\exp \left(\int_{r_0}^1 \mathcal{M}(s) \ ds \right) \int_{r_0}^\delta \mathcal{F}(\rho) \ d\rho, \label{eqn: upper bound int F(r) i}\\
    \int_{R/4}^{r_0} \mathcal{F}_R(r) \ dr &\leq \frac{r_0 - R/4}{\delta - r_0}\exp \left(\int_{R/4}^1 \mathcal{M}(s) \ ds \right) \int_{r_0}^\delta \mathcal{F}_R(\rho) \ d\rho. \label{eqn: upper bound int FR(r) i}
\end{align}

The third step simply sums the above two estimates. The fourth and fifth steps then respectively bound from below and from above the obtained sum. The last step concludes. \\

\textit{First step: proof of \eqref{eqn: upper bound int F(r) i}.} Differentiating $\mathcal{F}$, we obtain
\begin{align*}
    \mathcal{F}'(r) &=  - \left( w_s^2 + w_r^2 + \tilde{V}(r)w^2 \right)(1-r+\eta, r) - \left( w_s^2 + w_r^2 + \tilde{V}(r)w^2 \right)(-1+r-\eta, r) \\
    &+\int_{-1+r-\eta}^{1-r+\eta} \tilde{V}'(r)w(s,r)^2 \ ds + 2\int_{-1+r-\eta}^{1-r+\eta} \left( w_sw_{sr} + w_rw_{rr} + \tilde{V}(r)ww_r \right) (s,r) \ ds. 
\end{align*}

Substituting $w_{rr}$ using \eqref{stm: radial wave adjoint system} and integrating by parts, we have 
\begin{align*}
    &2\int_{-1+r-\eta}^{1-r+\eta} w_rw_{rr} (s,r) \ ds = 2\int_{-1+r-\eta}^{1-r+\eta} w_rw_{ss} (s,r) + V(r)w_rw(s,r) \ ds \\
    &= 2w_rw_s(1-r+\eta,r) - 2w_rw_s(-1+r-\eta,r) + 2 \int_{-1+r-\eta}^{1-r+\eta} V(r)w_rw(s,r) - w_{rs}w_s (s,r) \ ds.
\end{align*}
Hence, using that $V \leq \tilde{V}$ and applying Young's inequality, 
\begin{align*}
    \mathcal{F}'(r) &=  - \left( (w_s - w_r)^2 + \tilde{V}(r)w^2 \right)(1-r+\eta, r) - \left( (w_s + w_r)^2 + \tilde{V}(r)w^2 \right)(-1+r-\eta, r) \\
    &+ \int_{-1+r-\eta}^{1-r+\eta} \tilde{V}'(r)w(s,r)^2 \ ds + 2\int_{-1+r-\eta}^{1-r+\eta} \left( (\tilde{V}(r) + V(r))ww_r \right) (s,r) \ ds \\
    &\leq \left( M(r) + 2\sqrt{\tilde{V}(r)} \right) \mathcal{F}(r),
\end{align*}

where $M(r) = |\tilde{V}'(r)|/\tilde{V}(r)$. Recall that $a < r_0 < \delta < b$. Setting $\mathcal{M}(r) := M(r) + 2\sqrt{\tilde{V}(r)}$, it follows that the mapping 
\begin{align*}
    [r_0, 1)\ni r \mapsto \mathcal{F}(r) \exp \left( - \int_{r_0}^r\mathcal{M}(s) \ ds\right)
\end{align*}
is nonincreasing. We thus have, for any $r_0 \leq \rho \leq \delta  \leq  r \leq 1$,
\begin{align}\label{eqn: upper bound F(r)}
    \mathcal{F}(r) \leq \mathcal{F}(\rho)\exp \left(\int_\rho^r \mathcal{M}(s) \ ds \right) \leq \mathcal{F}(\rho) \exp \left(\int_{r_0}^1 \mathcal{M}(s) \ ds \right).
\end{align}

Integrating \eqref{eqn: upper bound F(r)} with respect to $\rho$ on $(r_0,\delta)$, and then to $r$ on $(\delta, 1)$, we obtain 
\begin{align*}
    \int_\delta^1  \mathcal{F}(r) \ dr \leq \frac{1-\delta}{\delta - r_0}\exp \left(\int_{r_0}^1 \mathcal{M}(s) \ ds \right) \int_{r_0}^\delta \mathcal{F}(\rho) \ d\rho.
\end{align*}
Adding the integral of $\mathcal{F}$ on $(r_0,\delta)$ on both sides, we arrive at the sought estimate
\begin{align*}
    \int_{r_0}^1  \mathcal{F}(r) \ dr \leq \frac{1-r_0}{\delta - r_0}\exp \left(\int_{r_0}^1 \mathcal{M}(s) \ ds \right) \int_{r_0}^\delta \mathcal{F}(\rho) \ d\rho.
\end{align*}

\textit{Second step: proof of \eqref{eqn: upper bound int FR(r) i}.} Following the same arguments as for the computations of $\mathcal{F}$, taking care of the sign of the derivatives, we arrive at 
\begin{align*}
    \mathcal{F}'_R(r) \geq - \mathcal{M}(r) \mathcal{F}_R(r),
\end{align*}

It follows that the mapping 
\begin{align*}
    \left[ \frac{R}{4}, \delta \right] \ni r \mapsto \mathcal{F}_R(r) \exp \left( \int_{R/4}^r\mathcal{M}(s) \ ds\right)
\end{align*}
is nondecreasing. We thus have, for any $R/4 \leq r \leq r_0 \leq \rho \leq \delta$,
\begin{align}\label{eqn: upper bound FR(r)}
    \mathcal{F}_R(r) \leq \mathcal{F}_R(\rho)\exp \left(\int_{R/4}^\delta \mathcal{M}(s) \ ds \right).
\end{align}
Integrating \eqref{eqn: upper bound FR(r)} with respect to $r$ on $(R/4,r_0)$, and to $\rho$ on $(r_0,\delta)$, we obtain the sought estimate 
\begin{align*}
    \int_{R/4}^{r_0} \mathcal{F}_R(r) \ dr \leq \frac{r_0 - R/4}{\delta - r_0}\exp \left(\int_{R/4}^1 \mathcal{M}(s) \ ds \right) \int_{r_0}^\delta \mathcal{F}_R(\rho) \ d\rho.
\end{align*}

\textit{Third step: summing \eqref{eqn: upper bound int F(r) i} and \eqref{eqn: upper bound int FR(r) i}.} Summing \eqref{eqn: upper bound int F(r) i} and \eqref{eqn: upper bound int FR(r) i}, we obtain 
\begin{align}\label{eqn: estimate sum F FR i}
      \int_{R/4}^{r_0} \mathcal{F}_R(r) \ dr + \int_{r_0}^1  \mathcal{F}(r) \ dr &\leq \frac{2}{\delta - r_0}\exp \left(\int_{R/4}^1 \mathcal{M}(s) \ ds \right) \left( \int_{r_0}^\delta \mathcal{F}_R(\rho) + \mathcal{F}(\rho) \ d\rho \right).
\end{align}

\textit{Fourth step: bounding \eqref{eqn: estimate sum F FR i} from below.}
On the one hand, since $r - \frac{R}{4} + \eta \geq \eta$ for $r \geq R/4$, and since $\tilde{V} \geq 1$, we have
\begin{align*}
    \int_{R/4}^{r_0} \mathcal{F}_R(r) \ dr \geq \int_{R/4}^{r_0} \int_{-\eta}^{\eta}  w_s(s,r)^2 + w(s,r)^2 \ ds \ dr.
\end{align*}
Assuming moreover $\eta \geq 5R$ and $2R \leq r_0$, it also follows that
\begin{align*}
    \int_{R/4}^{r_0} \mathcal{F}_R(r) \ dr \geq \int_{-5R}^{5R} \int_{R/4}^{2R} w_s(s,r)^2 + w(s,r)^2  \ dr \ ds.
\end{align*}
In particular,
\begin{align}\label{eqn: lower bound FR i}
    \int_{R/4}^{r_0} \mathcal{F}_R(r) \ dr \geq \frac{1}{2} \int_{R/4}^{r_0} \int_{-\eta}^{\eta} w_s(s,r)^2 \ ds \ dr + \frac{1}{2} \int_{-5R}^{5R} \int_{R/4}^{2R} w_s(s,r)^2 + w(s,r)^2 \ dr \ ds.
\end{align}

On the other hand, using $1-r+\eta \geq \eta$, $-1+r-\eta \leq -\eta$, and $\tilde{V} \geq 1$, we similarly obtain  
\begin{align}\label{eqn: lower bound F i}
    \int_{r_0}^1  \mathcal{F}(r) \ dr \geq \frac{1}{2} \int_{r_0}^1 \int_{-\eta}^\eta w_s(s,r)^2 \ ds \ dr. 
\end{align}

Substituting \eqref{eqn: lower bound FR i} and \eqref{eqn: lower bound F i} into \eqref{eqn: estimate sum F FR i}, we thus obtain
\begin{align}\label{eqn: estimate sum F FR ii}
       \begin{aligned}
           \int_{-5R}^{5R} \int_{R/4}^{2R} w_s(s,r)^2 + w(s,r)^2  \ dr \ ds &+ \int_{R/4}^1 \int_{-\eta}^\eta w_s(s,r)^2 \ ds \ dr \\
           &\leq \frac{4}{\delta - r_0}\exp \left(\int_{R/4}^1 \mathcal{M}(s) \ ds \right) \left( \int_{r_0}^\delta \mathcal{F}_R(\rho) + \mathcal{F}(\rho) \ d\rho \right).
       \end{aligned}
\end{align}

\textit{Fifth step: bounding \eqref{eqn: estimate sum F FR ii} from above.} We first treat the term involving $\mathcal{F}(\rho)$ in the right-hand side of \eqref{eqn: estimate sum F FR ii}. The term involving $\mathcal{F}_R(\rho)$ is estimated similarly. \\

We start with the term in $w_r$ appearing in $\mathcal{F}(\rho)$ in the right-hand side of \eqref{eqn: estimate sum F FR ii}, using a Caccioppoli-type argument. Multiplying \eqref{stm: radial wave adjoint system} by $w$, we have the identity
\begin{align}\label{eqn: identity wave}
    w_r^2 = w_s^2 - (ww_s)_s + (ww_r)_r - Vw^2.
\end{align}

Introducing a cutoff $\chi \in C_c^\infty((-S,S) \times (a,b))$, $0 \leq \chi \leq 1$, $\chi = 1$ on $[-1+ r_0 - \eta, 1- r_0 +\eta] \times [r_0, \delta]$, using \eqref{eqn: identity wave} and integrating by parts, we obtain,
\begin{align*}
    \int_{r_0}^\delta \int_{-1+ \rho - \eta}^{1-\rho+\eta} w_r(s,\rho)^2 \ ds \ d\rho &\leq \int_a^{b} \int_{-S}^S  \chi(s,\rho) w_r(s,\rho)^2 \ ds \ d\rho \\
    &= \int_a^{b} \int_{-S}^S \chi(s,\rho) \left( w_s^2 - (ww_s)_s + (ww_r)_r - Vw^2 \right)(s,\rho) \ ds \ d\rho  \\
    &= \int_a^{b} \int_{-S}^S \chi w_s^2 + \left( \frac{1}{2}\chi_{rr} - \frac{1}{2}\chi_{ss}  - \chi V \right)w^2 \ ds \ d\rho.
\end{align*}
Since $-V(r) \leq 1/(4r^2) \leq 1/(4a^2)$ on $(a,b)$, it follows that, for some constant $C(\chi,a) > 0$ that only depends on $\chi$ and $a$, we have
\begin{align}\label{eqn: Caccioppoli type ineq}
     \int_{r_0}^\delta \int_{-1+ \rho - \eta}^{1-\rho+\eta} w_r(s,\rho)^2 \ ds \ d\rho \leq C(\chi,a) \int_a^{b} \int_{-S}^S \left( w_s^2 + w^2 \right) \ ds \ d\rho.
\end{align}

On the other hand, we have 
\begin{align}\label{eqn: estimate tilde V}
    \int_{r_0}^\delta \int_{-1+ \rho - \eta}^{1-\rho+\eta} \tilde{V}(\rho)w(s,\rho)^2 \ ds \ d\rho \leq \left( 1 + \xi^2 \delta^{2\gamma} + \frac{|\nu_\ell^2 - 1/4|}{r_0^2} \right) \int_a^{b} \int_{-S}^S w(s,\rho)^2 \ ds \ d\rho
\end{align}

From \eqref{eqn: Caccioppoli type ineq} and \eqref{eqn: estimate tilde V}, we obtain 
\begin{align}\label{eqn: upper bound int F(rho) ii}
    \begin{aligned}
        \int_{r_0}^\delta  \mathcal{F}(\rho) \ d\rho \leq \left( C(\chi,a) + 1 + \xi^2 \delta^{2\gamma} + \frac{|\nu_\ell^2 - 1/4|}{r_0^2} \right)\int_a^b \int_{-S}^S w_s^2 + w^2 \ ds \ d\rho.
    \end{aligned}
\end{align}

Similarly, the same cutoff $\chi$ arguments yield
\begin{align}\label{eqn: upper bound int FR(rho) ii}
    \int_{r_0}^\delta  \mathcal{F}_R(\rho) \ d\rho \leq \left( C(\chi) + 2 + \xi^2 \delta^{2\gamma} + \frac{|\nu_\ell^2 - 1/4|}{r_0^2} \right) \int_a^{b} \int_{-S}^S \left( w_s^2 + w^2 \right) \ ds \ d\rho.
\end{align}

Substituting \eqref{eqn: upper bound int F(rho) ii} and \eqref{eqn: upper bound int FR(rho) ii} into \eqref{eqn: estimate sum F FR ii}, we obtain 
\begin{align}\label{eqn: estimate sum F FR iii}
       \begin{aligned}
           &\int_{-5R}^{5R} \int_{R/4}^{2R} w_s(s,r)^2 + w(s,r)^2  \ dr \ ds + \int_{R/4}^1 \int_{-\eta}^\eta w_s(s,r)^2 \ ds \ dr \\
           &\leq \frac{8}{\delta - r_0}\exp \left(\int_{R/4}^1 \mathcal{M}(s) \ ds \right) \left( C(\chi) + 2 + \xi^2 \delta^{2\gamma} + \frac{|\nu_\ell^2 - 1/4|}{r_0^2} \right) \int_a^{b} \int_{-S}^S \left( w_s^2 + w^2 \right) \ ds \ d\rho.
       \end{aligned}
\end{align}

This concludes the proof.
\end{proof}

\section{Proof of the controllability results in the critical degenerate case: Theorem \ref{thm: main thm minimal time gamma = 1}}\label{section: proof of the controllability results in the critical degenerate case}

We prove Theorem \ref{thm: uniform control annulus gamma = 1} in Section \ref{section: Proof of Theorem 2}, and the positive part of Theorem \ref{thm: main thm minimal time gamma = 1} in Section \ref{section: Proof of Theorem 1 positive}, and the negative part of Theorem \ref{thm: main thm minimal time gamma = 1} in Section \ref{section: Proof of Theorem 1 negative}.

\subsection{Proof of Theorem \ref{thm: uniform control annulus gamma = 1}: control on annuli}\label{section: Proof of Theorem 2}

We will use the following Theorem \ref{thm: uniform biorthogonal families} on existence of, and estimates for, biorthogonal families, which is a consequence of \cite[Theorems IV.1.8, IV.1.10 and IV.1.18]{boyer:hal-02470625}. For the reader's convenience, we give a short proof of this theorem, with precise references to \cite{boyer:hal-02470625}, after the proof of Theorem \ref{thm: uniform control annulus gamma = 1}.

\begin{theorem}\label{thm: uniform biorthogonal families}
    Let $\Lambda = (\Lambda_j)_{j \geq 0} \subset (0,+\infty)$ be an increasing sequence of real numbers. For every $E > 0$, denote the counting function by $\mathcal{N}(E) := \# \{j \geq 0: \Lambda_j \leq E\}$. Assume that there exist $\rho, N> 0$ and $\alpha \in (0,1)$ such that
    \begin{align}\label{hyp: condition eigenv biorth}
        \begin{array}{clll}
            \left| \Lambda_{i} - \Lambda_j \right| &\geq& \rho, & \forall i \neq j ,\\[6pt]
            \mathcal{N}(E) &\leq& NE^\alpha, & \forall E > 0.
        \end{array}
    \end{align}
    Then, there exists $C = C(N,\alpha,\rho)$ such that, for every $T > 0$, there exists a real biorthogonal family $(q_j)_{j \geq 0} \subset L^2(0,T)$ satisfying
    \begin{align}
        \int_0^T q_j(t)e^{-\Lambda_i t} \ dt = \delta_{ij}, \qquad \forall i,j \geq 0,
    \end{align}
    and
    \begin{align}\label{eqn: estimate biorth family boyer}
        \|q_j\|_{L^2(0,T)} \leq  C e^{CT^{-\frac{\alpha}{1-\alpha}}}e^{\Lambda_j^\alpha(C+\ln \Lambda_j)},
        \qquad \forall j \geq 0.
    \end{align}
\end{theorem}

\begin{proof}[Proof of Theorem \ref{thm: uniform control annulus gamma = 1}]
  
Let $0 \leq a < b \leq 1$, $\mathcal{A} = \{x\in B(0,1):a<|x|<b\}$, and $T > 0$. For every $n \geq 1$ such that $\xi_n > 0$, let $f_n^0\in L^2(B(0,1))$, and given a control $u_n$, denote by $f_n$ the associated solution of system \eqref{stm : Fourier heat control system} with $\xi = \xi_n$ and $\mathcal{U} = \mathcal{A}$. Recall from Section \ref{section: Spherical harmonic decomposition} that the state and the control are decomposed as
\begin{align*}
    f_n(t,r\theta) =r^{-(d_x-1)/2} \sum_{\ell\geq 0}\sum_{m=1}^{m_\ell} \mathsf{F}_{n,\ell,m} (t,r)Y_{\ell,m}(\theta), \quad u_n(t,r\theta) =r^{-(d_x-1)/2} \sum_{\ell\geq 0}\sum_{m=1}^{m_\ell}\mathsf{u}_{n,\ell,m}(t,r)Y_{\ell,m}(\theta).
\end{align*}
For each $(\ell,m)$, $\mathsf{F}_{n,\ell,m}$ solves system \eqref{stm: radial heat control system} with $\xi=\xi_n$, control $\mathsf{u}_{n,\ell,m}$, and initial state
\begin{align*}
    \mathsf{F}_{n,\ell,m}^0(r) = r^{(d_x-1)/2} \int_{\mathbb{S}^{d_x-1}}f_n^0(r\theta)Y_{\ell,m}(\theta)\ d\theta.
\end{align*}
We now construct these radial controls $\mathsf{u}_{n,\ell,m}$ that, we recall, are formally given by 
\begin{align}\label{eqn: def control radial in proof}
    \mathsf{u}_{n,\ell,m}(t,r) = -  \sum_{k \geq 0} e^{-\lambda_{n,\ell,k}T} \mathsf{F}_{n,\ell,m,k}^0 \frac{ \varphi_{n,\ell,k}(r)} { \|\varphi_{n,\ell,k}\|_{L^2(a,b)}^2} q_{n, \ell,k}(T-t),
\end{align}

where we set $\mathsf{F}_{n,\ell,m,k}^0 = \langle \mathsf{F}_{n,\ell,m}^0,\varphi_{n,\ell,k} \rangle_{L^2(0,1)}$.\\

First, we apply Theorem \ref{thm: uniform biorthogonal families} with the eigenvalues of $G_{n,\nu_\ell}$. We verify the conditions outlined in \eqref{hyp: condition eigenv biorth}. \\

Using Proposition \ref{prop: lower bounds eigenvalue}, for every $n,\ell$ we have $\lambda_{n,\ell,k} \geq c(k+1)^2$, from which it follows that, for every $E > 0$,
\begin{align*}
    \mathcal{N}(E) \leq \frac{1}{\sqrt{c}}\sqrt{E}.
\end{align*}
On the other hand, from Theorem \ref{thm: gap uniform}, for every $n,\ell,k$ we have
\begin{align*}
   \lambda_{n,\ell,k+1} - \lambda_{n,\ell,k} \geq 4\min \{\xi_n > 0\} > 0.
\end{align*}
Hence, we can apply Theorem \ref{thm: uniform biorthogonal families} with $\alpha = 1/2$, $\rho = 4\min \{\xi_n > 0\}$, $N = 1/\sqrt{c}$, and constants that are independent of $n,\ell,k$. We thus obtain, substituting \eqref{eqn: estimate biorth family boyer} into \eqref{eqn: def control radial in proof}, that
\begin{align}\label{eqn: estimate cost control i}
    \|\mathsf{u}_{n,\ell,m} \|_{L^2((0,T) \times (a,b))} \leq C \exp\left( \frac{C}{T} \right) \sum_{k \geq 0} \exp \left( - \left( T-\frac{C + \ln \lambda_{n,\ell,k}}{\sqrt{\lambda_{n,\ell,k}}} \right)\lambda_{n,\ell,k}\right) \frac{\left| \mathsf{F}_{n,\ell,m,k}^0 \right|}{\|\varphi_{n,\ell,k}\|_{L^2(a,b)}}.
\end{align}

Second, we apply Theorem \ref{thm: localization estimates eigenfunctions}. Let $\epsilon \in (0, (b-a)/2)$. There exists a constant $C > 0$ such that, for every $n, \ell, k$, 
\begin{align}\label{eqn: estimate inverse localization eigenfunctions}
    \|\varphi_{n,\ell,k}\|_{L^2(a,b)}^{-1} \leq \frac{\lambda_{n,\ell,k}}{C} \exp \left( \sqrt{\lambda_{n,\ell,k}}\left( 1  + \ln \frac{1}{a+\epsilon} + \frac{\ln \lambda_{n,\ell,k}}{4} \right) \right)  \exp \left( \frac{\xi_n}{2} (a+2\epsilon)^2 \right).
\end{align}

Substituting \eqref{eqn: estimate inverse localization eigenfunctions} into \eqref{eqn: estimate cost control i}, we obtain that there exists $C > 0$, that may differ from the preceding ones, such that, for every $n,\ell,m,k$, we have
\begin{align*}
     \begin{aligned}
         &\|\mathsf{u}_{n,\ell,m} \|_{L^2((0,T) \times (a,b))} \leq  C \exp\left( \frac{C}{T} \right) \exp \left( \frac{\xi_n}{2} (a+2\epsilon)^2 \right) \\
         &\times \sum_{k \geq 0} \lambda_{n,\ell,k}  |\mathsf{F}_{n,\ell,m,k}^0| \exp \left( - \left( T - \frac{C + \ln \lambda_{n,\ell,k}}{\sqrt{\lambda_{n,\ell,k}}} - \frac{1}{\sqrt{\lambda_{n,\ell,k}}}\left( 1  + \ln \frac{1}{a+\epsilon} + \frac{\ln \lambda_{n,\ell,k}}{4} \right) \right)\lambda_{n,\ell,k}\right),
     \end{aligned}
\end{align*}
that we further simplify into, for some $C(\epsilon) > 0$,
\begin{align}\label{eqn: estimate cost control iii}
     \begin{aligned}
         \|\mathsf{u}_{n,\ell,m} \|_{L^2((0,T) \times (a,b))} &\leq  C \exp\left( \frac{C}{T} \right) \exp \left( \frac{\xi_n}{2} (a+2\epsilon)^2 \right) \\
         &\times \sum_{k \geq 0} \lambda_{n,\ell,k}  |\mathsf{F}_{n,\ell,m,k}^0| \exp \left( -T\lambda_{n,\ell,k}
        + C(\epsilon)\sqrt{\lambda_{n,\ell,k}} \ln\lambda_{n,\ell,k} \right).
     \end{aligned}
\end{align}
We have that for every $\delta \in (0,1/2)$ and $\beta > 1$, there exists $C_0 = C(\epsilon,\delta,\beta) > 0$ such that 
\begin{align*}
    \frac{C}{T} + C(\epsilon)\sqrt{\lambda}\ln\lambda \leq \delta T\lambda + \frac{C_0}{T^\beta} = 2\delta T\lambda - \delta T\lambda + \frac{C_0}{T^\beta}.
\end{align*}
Moreover, recalling that $\lambda_{n,\ell,k} \geq 2\xi_n(1+\nu_\ell)$ from Proposition \ref{prop: lower bounds eigenvalue}, we have
\begin{align*}
    \exp \left( -(1-2\delta)T\lambda_{n,\ell,k} \right) \leq \exp \left( -2(1-2\delta)T\xi_n(1+\nu_\ell)\right).
\end{align*}
Substituting these estimates into \eqref{eqn: estimate cost control iii}, and then using
Cauchy-Schwarz and Parseval, we obtain, for some new $C > 0$ that differs from the preceding ones, and that does not depend on $n,\ell,m,k$,
\begin{align}\label{eqn: estimate cost control iv}
    \begin{aligned}
        \|\mathsf{u}_{n,\ell,m}\|_{L^2((0,T)\times(a,b))} &\leq C \exp\left( \frac{C}{T^\beta} \right) \exp\left( \frac{\xi_n}{2}(a+2\epsilon)^2 - 2(1 - 2\delta)T\xi_n(1+\nu_\ell) \right) \\
        &\times  \left( \sum_{k\geq 0}\lambda_{n,\ell,k}^2 e^{-2\delta T\lambda_{n,\ell,k}} \right)^{1/2} \|\mathsf{F}_{n,\ell,m}^0\|_{L^2(0,1)}.
    \end{aligned}
\end{align}
Since we have $\lambda_{n,\ell,k}\geq c(k+1)^2$ by Proposition \ref{prop: lower bounds eigenvalue}, we obtain, uniformly in $n,\ell$,
\begin{align*}
    \sum_{k\geq 0}\lambda_{n,\ell,k}^2 e^{- 2\delta T\lambda_{n,\ell,k}} \leq \frac{16}{\delta^2 T^2} \sum_{k \geq 0}e^{-\frac{\delta T}{2}\lambda_{n,\ell,k}} \leq \frac{16}{\delta^2 T^2} \sum_{k \geq 0}e^{-\frac{c\delta T}{2}(k+1)^2} \leq \frac{16}{\delta^2 T^2} \sqrt{\frac{\pi}{2c\delta T}}.
\end{align*}
Substituting this into \eqref{eqn: estimate cost control iv}, and absorbing the resulting $T^{-5/4}$ polynomial term into the first exponential, we arrive at, for some $C > 0$ that differs from the preceding ones, and that does not depend on $n,\ell,m,k$,
\begin{align}\label{eqn: control cost minimal time from annulus}
    \|\mathsf{u}_{n,\ell,m}\|_{L^2((0,T)\times(a,b))} \leq C \exp\left( \frac{C}{T^\beta} \right) \exp\left( \frac{\xi_n}{2}(a+2\epsilon)^2 - 2(1 - 2\delta)T\xi_n (1+\nu_\ell) \right) \|\mathsf{F}_{n,\ell,m}^0\|_{L^2(0,1)}.
\end{align}

Applying Parseval's identity, summing, and using $1+\nu_\ell \geq 1 + \nu$, we finally obtain 
\begin{align}\label{eqn: final control cost annulus proof}
    \|u_n\|_{L^2((0,T)\times\mathcal{A})} \leq C \exp\left(\frac{C}{T^\beta}\right) \exp\left( \frac{\xi_n}{2}(a+2\epsilon)^2 -2(1-2\delta)(1+\nu)T\xi_n \right) \|f_n^0\|_{L^2(B(0,1))},
\end{align}
which concludes the proof. 
\end{proof}

Let us give a brief proof of Theorem \ref{thm: uniform biorthogonal families}.

\begin{proof}[Proof of Theorem \ref{thm: uniform biorthogonal families}]
    Under the assumptions of Theorem \ref{thm: uniform biorthogonal families}, it follows from \cite[Theorems IV.1.8 and IV.1.10]{boyer:hal-02470625} that there exists a biorthogonal family $(q_j^\infty)_j \subset L^2(0,+\infty)$ satisfying
    \begin{align}
        \int_0^{+ \infty} q_j^\infty(t)e^{-\Lambda_i t} \ dt = \delta_{ij}, \qquad \forall i,j \geq 0,
    \end{align}
    and
    \begin{align}\label{eqn: estimate biorth half line}
        \|q_j^\infty\|_{L^2(0,+\infty)} \leq  \exp \left( \Lambda_j^\alpha (C + \ln \Lambda_j) \right),
        \qquad \forall j \geq 0.
    \end{align}
    where $C > 0$ depends only on $N,\alpha,\rho$ appearing in \eqref{hyp: condition eigenv biorth}. Introduce the spaces (see \cite[Sections IV.1.2.2 and IV.1.2.3]{boyer:hal-02470625})
    \begin{align*}
        F_{\Lambda,\infty} = \overline{\mathrm{Span}\left( t \mapsto e^{-\Lambda_j t}, \Lambda_j \in \Lambda \right)}^{L^2(0,+\infty)}, \quad F_{\Lambda,T} = \overline{\mathrm{Span}\left( t \mapsto e^{-\Lambda_j t}, \Lambda_j \in \Lambda \right)}^{L^2(0,T)}.
    \end{align*}

    We denote by $\Gamma_{\Lambda,T} : F_{\Lambda,\infty} \to F_{\Lambda,T}$ the restriction operator to $(0,T)$. From \cite[Theorem IV.1.18]{boyer:hal-02470625}, it is invertible, and its inverse satisfies, for some $C_1 > 0$ that depends only on $N,\alpha$ appearing in \eqref{hyp: condition eigenv biorth}, 
    \begin{align}\label{eqn: estimate norm restriction operator biorth fam}
        \left\| \Gamma_{\Lambda,T}^{-1} \right\| \leq C_1 \exp \left( C_1 T^{-\frac{\alpha}{1 - \alpha}} \right).
    \end{align}
    We set 
    \begin{align*}
        q_j := \left(\Gamma_{\Lambda,T}^{-1} \right)^* q_j^\infty .
    \end{align*}
    By construction in \cite[Page 61]{boyer:hal-02470625}, we have $q_j^\infty \in F_{\Lambda,\infty}$, and thus,
    \begin{align*}
        \delta_{ij} = \left(q_j^\infty, \exp(-\Lambda_i \cdot)\right)_{L^2(0,+\infty)} &= (q_j^\infty, \Gamma_{\Lambda,T}^{-1} \Gamma_{\Lambda,T} \exp(-\Lambda_i \cdot))_{L^2(0,+\infty)} \\
        &= \left( \left(\Gamma_{\Lambda,T}^{-1} \right)^* q_j^\infty, \Gamma_{\Lambda,T} \exp(-\Lambda_i \cdot)\right)_{L^2(0,T)}.
    \end{align*}
    The estimate \eqref{eqn: estimate biorth family boyer} simply follows from 
    \begin{align*}
        \left\|q_j \right\|_{L^2(0,T)} = \left\| \left(\Gamma_{\Lambda,T}^{-1} \right)^* q_j^\infty \right\|_{L^2(0,T)}
        \leq \left\|\left(\Gamma_{\Lambda,T}^{-1} \right)^* \right\| \left\| q_j^\infty \right\|_{L^2(0,+\infty)}  = \left\| \Gamma_{\Lambda,T}^{-1} \right\| \left\| q_j^\infty \right\|_{L^2(0,+\infty)},
    \end{align*}
    combined with \eqref{eqn: estimate biorth half line} and \eqref{eqn: estimate norm restriction operator biorth fam}.
\end{proof}

\subsection{Proof of the upper bound of Theorem \ref{thm: main thm minimal time gamma = 1}: angular reduction and Lebeau-Robbiano construction}\label{section: Proof of Theorem 1 positive}

We can now prove the upper bound part of Theorem \ref{thm: main thm minimal time gamma = 1}.

\begin{proof}[Proof of the upper bound in Theorem \ref{thm: main thm minimal time gamma = 1}]

Let $f_0 \in L^2(\Omega)$, and given a control $u$, denote by $f$ the corresponding solution of system \eqref{stm : control system}. Recall that we denote by $f_n$ its Fourier coefficients, and those of $f_0$ by $f_n^0$, and that the $f_n$ are solutions of system \eqref{stm : Fourier heat control system} with $\xi = \xi_n$ and initial state $f_n^0$. Recall also, from Section \ref{section: Fourier decomposition and the control cost from annuli}, that it suffices to show that for every $T > \operatorname{dist}(\omega_x, 0_{\mathbb{R}^{d_x}})^2/(4(1+\nu))$, there exists $C(T) > 0$ such that, for every $n$, there exists a control $u_n \in L^2((0,T) \times \omega_x)$ that steers $f_n^0$ to $0$ in time $T > 0$, and satisfies
\begin{align}
    \|u_n\|_{L^2((0,T) \times \omega_x)} \leq C(T) \|f_n^0\|_{L^2(B(0,1))}.
\end{align}

Let $\Gamma \subset \mathbb{S}^{d_x-1}$ be a nonempty open subset, and set
\begin{align*}
    \mathcal{A}_\Gamma
    := \{x=r\theta\in B(0,1):a<r<b,\ \theta\in\Gamma\},
\end{align*}
which will be chosen so that $\mathcal{A}_\Gamma \subset \omega_x$. \\

Fix $\epsilon\in(0,(b-a)/2)$ and $\beta>1$. Recall that we denote by $\mathcal{G}_n$ the Friedrichs operator of the Fourier system \eqref{stm : Fourier heat control system} with $\xi = \xi_n$, and by $S_n(t)=e^{-t\mathcal{G}_n}$ the continuous semigroup of contractions that it generates. For $L \geq 1$, we denote by $\Pi_L$ the spectral projections onto 
\begin{align*}
    \mathcal{V}_L:=\operatorname{Span} \{Y_{\ell,m}: \ell\leq L,\ 1\leq m\leq m_\ell\},
\end{align*}
and its action onto elements of $L^2(B(0,1))$ is understood as, setting $x=r\theta$, for almost every $r \in (0,1)$, the mapping $\theta \mapsto \Pi_L f(r\theta)$ belongs to $\mathcal{V}_L$. We slightly abuse notation here to simplify the presentation, and we shall use $x  =r\theta$ in what follows when necessary. \\

Moreover, by Proposition \ref{prop: lower bounds eigenvalue}, $\lambda_{n,\ell,k} \geq \nu_\ell^2 \geq \ell^2$, and we have the dissipation estimate
\begin{align}\label{eqn: angular high frequency dissipation}
    \|S_n(t)(\operatorname{Id}-\Pi_L) g \|_{L^2(B(0,1))} \leq e^{-L^2t}\|g\|_{L^2(B(0,1))}, \qquad t > 0.
\end{align}

\textit{First step: }Let $\tau > 0$. We first show that for every $L > 0$, there exists a control $v_n$ supported in $\mathcal{A}_\Gamma$ such that the corresponding solution $f_n$, with initial state $f_n^0$, satisfies $\Pi_Lf_n(\tau) = 0$, with a quantitative estimate on the control cost. \\

Applying Theorem \ref{thm: uniform control annulus gamma = 1} with $\delta=1/4$ for instance, for every $ \tau > 0$ and initial state $f_n^0 \in L^2(B(0,1))$, there exists a control $w_n$ supported in $\mathcal{A}$ satisfying
\begin{align}\label{eqn: uniform cost fourier annulus gamma = 1 in proof i}
        \|w_n\|_{L^2((0,\tau)\times\mathcal{A})} \leq C \exp\left(\frac{C}{\tau^\beta}\right) \exp\left( \frac{\xi_n}{2}(a+2\epsilon)^2 \right) \|\Pi_L f_n^0\|_{L^2(B(0,1))},
    \end{align}
that steers $\Pi_L f_n^0$ to $0$ in time $\tau$ in $L^2(B(0,1))$. If we denote by $g_n$ the associated solution with initial state $\Pi_L f_n^0$, then $g_n(\tau) = 0$. The solution of system \eqref{stm : Fourier heat control system} starting from $(\mathrm{Id}-\Pi_L)f_n^0$, and without control, is $S_n(t)(\mathrm{Id}-\Pi_L)f_n^0$. Thus, denoting by $h_n$ the solution of system \eqref{stm : Fourier heat control system} with initial state $f_n^0$ and control $w_n$, it is given by
\begin{align*}
    h_n(t) = g_n(t) + S_n(t)(\mathrm{Id}-\Pi_L)f_n^0.
\end{align*}
Using that $\Pi_L$ commutes with $S_n(t)$, and with $\mathbbm{1}_\mathcal{A}$ since $\mathbbm{1}_\mathcal{A}$ does not depend on $\theta$, we thus have that $\Pi_L(h_n(t)) = \Pi_L(g_n(t))$, and that  
\begin{align}
    (\partial_t + \mathcal{G}_n)\Pi_L(h_n(t)) = \Pi_L\mathbbm{1}_\mathcal{A} w_n = \mathbbm{1}_\mathcal{A} \Pi_L w_n, \quad \Pi_L(h_n(0)) = \Pi_L f_n^0, \quad \Pi_L(h_n(\tau)) = 0.
\end{align}
In particular, the solution $f_n$ of system \eqref{stm : Fourier heat control system} with initial state $f_n^0$ and control $\tilde{w}_n := \Pi_L  w_n$ satisfies $\Pi_Lf_n(\tau) = 0$. Moreover, $\tilde{w}_n(t,r\cdot) \in \mathcal{V}_L$, and it satisfies
    \begin{align}\label{eqn: uniform cost fourier annulus gamma = 1 in proof ii}
        \|\tilde{w}_n\|_{L^2((0,\tau)\times\mathcal{A})} \leq C \exp\left(\frac{C}{\tau^\beta}\right) \exp\left( \frac{\xi_n}{2}(a+2\epsilon)^2 \right) \|\Pi_L f_n^0\|_{L^2(B(0,1))}.
    \end{align}
    
We now construct $v_n$ from $\tilde{w}_n$. Recall that we want $v_n$ to be supported in $\mathcal{A}_\Gamma$ and such that the solution $f_n$ of system \eqref{stm : Fourier heat control system} with initial state $f_n^0$ and control $v_n$ satisfies $\Pi_Lf_n(\tau) = 0$. For that, it suffices that $\Pi_L \mathbbm{1}_{\mathcal{A}_\Gamma} v_n = \Pi_L \mathbbm{1}_{\mathcal{A}} \tilde{w}_n = \mathbbm{1}_{\mathcal{A}} \tilde{w}_n$. Namely, denoting by $R_\Gamma$ the restriction operator $R_\Gamma : g \in \mathcal{V}_L \mapsto R_\Gamma g = g_{|\Gamma}$, it suffices that $v_n$ satisfies, for every $t \in (0,\tau)$ and almost every $r \in (0,1)$, 
\begin{align}\label{eqn: Riesz equation}
    \langle v_n(t,r\cdot), R_\Gamma g \rangle_{L^2(\mathbb{S}^{d_x - 1})} = \langle \tilde{w}_n(t,r\cdot), g \rangle_{L^2(\mathbb{S}^{d_x - 1})} , \quad \text{for every } g \in \mathcal{V}_L.
\end{align}
Recall the spectral inequality on the sphere (see \textit{e.g.} \cite{jerison1999nodal}). There exists $C(\Gamma) > 0$ such that, for every $g \in \mathcal{V}_L$,
\begin{align*}
    \|g\|_{L^2(\mathbb{S}^{d_x-1})}^2 \leq C(\Gamma) e^{C(\Gamma) L} \|g\|_{L^2(\Gamma)}^2,
\end{align*}
since the spherical eigenvalues satisfy $\sqrt{\Lambda_L} = \sqrt{L(L+d_x - 2)} \leq L + (d_x-2)/2 \leq CL$. In particular, $R_\Gamma$ is injective, and thus invertible on its image $R_\Gamma(\mathcal{V}_L)$ that we endow with the scalar product induced by $L^2(\mathbb{S}^{d_x - 1})$, with operator norm
\begin{align}\label{eqn: norm inverse of R}
    \|R_\Gamma^{-1} \| \leq C(\Gamma)^{1/2} e^{C(\Gamma) L/2} .
\end{align}
Hence, we can solve \eqref{eqn: Riesz equation}, and $v_n$ is given by $v_n(t,r\cdot) := (R_\Gamma^{-1})^* \tilde{w}_n(t,r\cdot)$, where $(R_\Gamma^{-1})^* : \mathcal{V}_L \to R_\Gamma(\mathcal{V}_L)$ is the adjoint of $R_\Gamma^{-1}$.
Indeed, for every $g \in \mathcal{V}_L$, for every $t > 0$ and $a < r < b$, we have, dropping the $(t,r\cdot)$ notation, 
\begin{align*}
    \langle \tilde{w}_n,g \rangle_{L^2(\mathbb{S}^{d_x - 1})} = \langle \tilde{w}_n, R_\Gamma^{-1} R_\Gamma g \rangle_{L^2(\mathbb{S}^{d_x - 1})}  =  \langle (R_\Gamma^{-1})^* \tilde{w}_n,  R_\Gamma g \rangle_{L^2(\mathbb{S}^{d_x - 1})} =  \langle (R_\Gamma^{-1})^* \tilde{w}_n,  g \rangle_{L^2(\Gamma)}.
\end{align*}
Moreover, using \eqref{eqn: uniform cost fourier annulus gamma = 1 in proof ii} and \eqref{eqn: norm inverse of R},
\begin{align}\label{eqn: uniform cost fourier annulus gamma = 1 in proof iii}
    \begin{aligned}
        \|v_n\|_{L^2((0,\tau) \times \mathcal{A}_\Gamma)} &\leq \|(R_\Gamma^{-1})^*\| \| \tilde{w}_n \|_{L^2((0,\tau) \times \mathcal{A})} =  \|R_\Gamma^{-1}\| \| \tilde{w}_n \|_{L^2((0,\tau) \times \mathcal{A})} \\
        &\leq C(\Gamma)^{1/2} e^{C(\Gamma) L/2}  C \exp\left(\frac{C}{\tau^\beta}\right) \exp\left( \frac{\xi_n}{2}(a+2\epsilon)^2 \right) \|\Pi_L f_n^0\|_{L^2(B(0,1))}.
    \end{aligned}
\end{align}
This concludes the first step.\\

\textit{Second step:} We now apply the Lebeau-Robbiano strategy (see \textit{e.g.} \cite[Section 6.2]{le2012carleman} and \cite{lebeau1995controle}), with controls constructed as in first step, and with costs as prescribed by \eqref{eqn: uniform cost fourier annulus gamma = 1 in proof iii}. \\

First, let $\tau \in (0,\min(1,T))$, and let $f_n^0$ be an initial state. We let the solution of system \eqref{stm : Fourier heat control system} dissipate freely on $(0,T-\tau)$, \textit{i.e.} $f_n(t) = S_n(t)f_n^0$. We have, using Proposition \ref{prop: lower bounds eigenvalue}, 
\begin{align}\label{eqn: dissipation estimate LR proof}
    \|S_n(T-\tau)f_n^0\| \leq e^{-2\xi_n(1+\nu)(T-\tau)}\|f_n^0\|.
\end{align}

The control we construct below using the Lebeau-Robbiano method drives $S_n(T-\tau)f_n^0$ to $0$ in time $\tau$. The construction is classical, and we follow \cite[Section 6.2]{le2012carleman}. Let us set, for simplicity of the presentation, 
\begin{align}\label{eqn: def Mn}
    \mathcal{C}_n := C(\Gamma)^{1/2}\exp\left(\frac{\xi_n}{2}(a+2\epsilon)^2\right),
\end{align}
where we assume $C(\Gamma) > 1$ without loss of generality. Set 
\begin{align}
    \tau_0 = 0, \quad \tau_{j+1} = \tau_j + 2T_j, \quad T_j = K2^{-j\rho},
\end{align}
where $\rho \in (0,1)$, and $K := \tau(1-2^{-\rho})/2$ so that $\tau = 2 \sum_{j \geq 0}T_j$ and $T_j \leq 1$ for $\tau \leq 2$. We also let $L_j$ be a sequence of real numbers to be chosen later. \\

For each $j \geq 0$, we apply on $(T-\tau + \tau_j,T-\tau + \tau_j + T_j)$ the control $v_n$ provided by the first step of the proof. We have, for some $C > 0$ that does not depend on $n$,
\begin{align}\label{eqn: cost control LR i}
    \begin{aligned}
        &\Pi_{L_j} f_n( T-\tau + \tau_j + T_j ) = 0, \\
        &\|v_n\|_{L^2((T-\tau + \tau_j,T-\tau + \tau_j + T_j) \times\mathcal{A}_\Gamma)} \leq \mathcal{C}_n \exp \left(\frac{C}{T_j^\beta} + C L_j \right) \|f_n(T-\tau + \tau_j)\|_{L^2(B(0,1))}.
    \end{aligned}
\end{align}

By Duhamel's formula and the contraction property of $S_n(t)$, it follows from \eqref{eqn: cost control LR i} that
\begin{align}\label{eqn: estimate solution controlled}
    \|f_n(T - \tau + \tau_j + T_j)\|_{L^2(B(0,1))} \leq \left( 1+\sqrt{T_j} \ \mathcal{C}_n \exp\left(\frac{C}{T_j^\beta} + CL_j \right) \right)\|f_n(T-\tau + \tau_j)\|_{L^2(B(0,1))}.
\end{align}
We then leave the system uncontrolled on $(T - \tau + \tau_j + T_j,  T - \tau + \tau_{j+1})$. Since $\Pi_{L_j} f_n(T - \tau + \tau_j + T_j )=0$, estimate \eqref{eqn: angular high frequency dissipation} yields
\begin{align}\label{eqn: estimate solution at rest}
    \|f_n(T - \tau + \tau_{j+1})\|_{L^2(B(0,1))} \leq e^{-L_j^2 T_j}\|f_n(T - \tau + \tau_j + T_j) \|_{L^2(B(0,1))}.
\end{align}
Combining \eqref{eqn: estimate solution controlled} and \eqref{eqn: estimate solution at rest}, we obtain
\begin{align}
    \begin{aligned}
        &\|f_n(T - \tau + \tau_{j+1})\|_{L^2(B(0,1))} \\
        &\leq \exp \left(- L_j^2 T_j \right)\left( 1+\sqrt{T_j} \ \mathcal{C}_n \exp\left(\frac{C}{T_j^\beta} + C L_j \right) \right)\|f_n(T-\tau + \tau_j)\|_{L^2(B(0,1))}.
    \end{aligned}
\end{align}
Hence, iterating the above estimate, we have that, using $T_j \leq 1$, for every $j \geq 0$, 
\begin{align}\label{eqn: estimate controlled solution LR ii}
    \begin{aligned}
        &\|f_n(T - \tau + \tau_{j+1})\|_{L^2(B(0,1))} \\
        &\leq \exp \left( \sum_{k=0}^j -L_k^2 T_k + \ln \left( 1 + \ \mathcal{C}_n \exp\left(\frac{C}{T_k^\beta} + C L_k \right) \right) \right)\|S_n(T-\tau)f_n^0\|_{L^2(B(0,1))}.
    \end{aligned}
\end{align}
As $j \to + \infty$, we have $T-\tau + \tau_j \to T$. Hence, to show that $f_n(T) = 0$, it suffices to show that the exponential factor in \eqref{eqn: estimate controlled solution LR ii} goes to zero as $j \to +\infty$. We have, for some $C_0 > 0$, using first that $\ln(1+x) \leq \ln 2 + \ln x$ for every $x \geq 1$, and then Young's inequality on the term $CL_k$,
\begin{align*}
    - L_k^2 T_k + \ln \left( 1 + \ \mathcal{C}_n \exp\left(\frac{C}{T_k^\beta} + C L_k \right) \right) &\leq -L_k^2 T_k + \ln 2 + \ln \mathcal{C}_n + \frac{C}{T_k^\beta} + C L_k \\
    &\leq -\frac{L_k^2}{2}T_k + \ln \mathcal{C}_n + \frac{C_0}{T_k^\beta}.
\end{align*}
Hence, setting for every $k$
\begin{align}\label{eqn: def L_j}
    L_k = 2 \sqrt{\frac{\ln \mathcal{C}_n}{T_k} + \frac{C_0}{T_k^{\beta + 1}}} \to + \infty, \text{ as } k \to + \infty,
\end{align}
we have that 
\begin{align}
    - L_k^2 T_k + \ln \left( 1 + \ \mathcal{C}_n \exp\left(\frac{C}{T_k^\beta} + C L_k \right) \right) \leq - \ln \mathcal{C}_n - C_0T_k^{-\beta}. 
\end{align}
Substituting this into \eqref{eqn: estimate controlled solution LR ii}, we obtain
\begin{align}\label{eqn: estimate controlled solution LR iii}
    \begin{aligned}
        \|f_n(T - \tau + \tau_{j+1})\|_{L^2(B(0,1))}\leq \mathcal{C}_n^{-j-1} \exp\left( - C_0 \sum_{k=0}^j T_k^{-\beta} \right)\|S_n(T-\tau)f_n^0\|_{L^2(B(0,1))}.
    \end{aligned}
\end{align}
The right-hand side of \eqref{eqn: estimate controlled solution LR iii} goes to zero as $j \to + \infty$, and it follows that $f_n(T) = 0$.\\

We now show that the constructed control is in $L^2((0,T) \times \mathcal{A}_\Gamma)$. Denote it by $u_n$, which is zero on $(0,T-\tau)$, equal to the active control constructed on each $(T-\tau + \tau_j, T - \tau + \tau_j + T_j)$, and is equal to zero on each $(T-\tau + \tau_j + T_j, T - \tau + \tau_{j+1})$. Combining \eqref{eqn: cost control LR i} and \eqref{eqn: estimate controlled solution LR iii}, it satisfies, for every $j \geq 1$,
\begin{align}\label{eqn: control cost LR ii}
    \|u_n\|_{L^2((T-\tau,T-\tau + \tau_{j+1}) \times \mathcal{A}_\Gamma)} \leq \sum_{i = 0}^j \mathcal{C}_n \exp \left( \frac{C}{T_i^\beta} + C L_i \right)  \exp\left( - C_0 \sum_{k=0}^{i-1} T_k^{-\beta} \right)\|S_n(T-\tau)f_n^0\|_{L^2(B(0,1))},
\end{align}
where we used $\mathcal{C}_n^{-i} \leq 1$. We bound from above the factor in \eqref{eqn: control cost LR ii}. We have, by definition \eqref{eqn: def L_j} of $L_i$, applying Young's inequality and using that $T_i \leq 1$, for every $\kappa > 0$, for some $C(\kappa) > 0$,
\begin{align*}
    \frac{C}{T_i^\beta} + C L_i \leq \kappa \ln \mathcal{C}_n +  \frac{C(\kappa)}{T_i^{\beta}}.
\end{align*}
Thus, substituting this into \eqref{eqn: control cost LR ii}, we have 
\begin{align}\label{eqn: control cost LR iii}
    \|u_n\|_{L^2((T-\tau,T-\tau + \tau_{j+1}) \times \mathcal{A}_\Gamma)} \leq \mathcal{C}_n^{\kappa + 1} \sum_{i =0}^j \exp \left( C(\kappa)T_i^{-\beta} - C_0 \sum_{k=0}^{i-1} T_k^{-\beta} \right) \|S_n(T-\tau)f_n^0\|_{L^2(B(0,1))}.
\end{align}
Computing the geometric sum in the exponential above using the definition of $T_j$, we have 
\begin{align*}
    C(\kappa)T_i^{-\beta} - C_0 \sum_{k=0}^{i-1} T_k^{-\beta} = \left( C(\kappa) - \frac{C_0}{2^{\rho \beta} - 1} \right)T_i^{-\beta} + \frac{C_0 T_0^{-\beta}}{2^{\rho \beta} - 1}.
\end{align*}
Choosing $\rho$ small enough such that $C(\kappa) - C_0/(2^{\rho \beta} - 1) < 0$, the sum appearing in \eqref{eqn: control cost LR iii} is thus bounded by, for some $C > 0$ that does not depend on $n$,  
\begin{align*}
    \sum_{i =0}^j \exp \left( C(\kappa)T_i^{-\beta} - C_0 \sum_{k=0}^{i-1} T_k^{-\beta} \right) &\leq \exp \left( \frac{C_0 T_0^{-\beta}}{2^{\rho \beta} - 1}\right)\sum_{i \geq 0} \exp\left( C(\kappa) - \frac{C_0}{2^{\rho \beta} - 1} \right)T_i^{-\beta} \\
    &\leq C\exp \left(C\tau^{-\beta} \right),
\end{align*}
where we used $T_0 \leq C(\rho)\tau$ by definition. Hence, combining \eqref{eqn: control cost LR iii} with \eqref{eqn: dissipation estimate LR proof}, replacing $\mathcal{C}_n$ by its definition \eqref{eqn: def Mn}, there exists some $C > 0$, that does not depend on $n$, such that, for every $T > 0$, for every $j \geq 1$, we have 
\begin{align*}
    &\|u_n\|_{L^2((T-\tau,T-\tau + \tau_{j+1}) \times \mathcal{A}_\Gamma)} \\
    &\leq C\exp \left(C\tau^{-\beta} \right) \exp\left(\frac{ \xi_n}{2}(\kappa +1)(a+2\epsilon)^2 - 2\xi_n(1+\nu)(T-\tau)\right)\|f_n^0\|_{L^2(B(0,1))}.
\end{align*}

By passing to the limit $j \to + \infty$, we obtain
\begin{align}\label{eqn: control cost LR iv}
    \|u_n\|_{L^2((0,T) \times \mathcal{A}_\Gamma)} \leq C\exp \left(C\tau^{-\beta} \right) \exp\left(\frac{ \xi_n}{2}(\kappa +1)(a+2\epsilon)^2 - 2\xi_n(1+\nu)(T-\tau)\right)\|f_n^0\|_{L^2(B(0,1))}.
\end{align}

\textit{Third step: Conclusion.} Observe that the cost appearing in \eqref{eqn: control cost LR iv} is uniform in $n$ if 
\begin{align}
    T > \frac{(a+2\epsilon)^2}{4(1+\nu)}(\kappa + 1) + \tau.
\end{align}
Since all the parameters can be chosen arbitrarily close to $0$, we deduce that, for every $\mathcal{A}_\Gamma$, we have 
\begin{align}
    T(\mathcal{A}_\Gamma \times \Omega_y) \leq \frac{a^2}{4(1+\nu)}.
\end{align}
To conclude, let $\eta > 0$ be such that there exists $x_\eta \in \omega_x$ with $|x_\eta| = \operatorname{dist}(\omega_x, 0_{\mathbb{R}^{d_x}}) + \eta$, and let $r_\eta > 0$ be such that $B(x_\eta,r_\eta) \subset \omega_x$. Such $\eta$ can be taken arbitrarily small. Choosing $\mathcal{A}_\Gamma \subset B(x_\eta,r_\eta)$ with $a < |x_\eta|$, and using that $\omega_0 \subset \omega_1$ implies $T(\omega_1) \leq T(\omega_0)$, we obtain, for every such $\eta > 0$,
\begin{align}
    T(\omega_x \times \Omega_y) \leq \frac{\left( \operatorname{dist} \left(\omega_x, 0_{\mathbb{R}^{d_x}} \right) + \eta \right)^2}{4(1+\nu)}.
\end{align}
This concludes the proof by passing to the limit $\eta \to 0$. 
\end{proof}

\subsection{Proof of the lower bound of Theorem \ref{thm: main thm minimal time gamma = 1}}\label{section: Proof of Theorem 1 negative}

We prove here that $T(\omega_x \times \Omega_y) \geq \operatorname{dist} \left(\omega_x, 0_{\mathbb{R}^{d_x}} \right)^2/(4(1+\nu))$. Namely, that for any
\begin{align*}
    T < \frac{\operatorname{dist} \left(\omega_x, 0_{\mathbb{R}^{d_x}} \right)^2}{4(1+\nu)},
\end{align*}
system \eqref{stm : control system} is not null-controllable from $\omega_x \times \Omega_y$. We already proved this statement in \cite[Theorem 1.3]{vanlaere2026observability} for $\nu > 0$, by means of Agmon estimates (see \cite[Proposition 5.1]{vanlaere2026observability}) whose proof holds only for $\nu > 0$. We propose here to take advantage of the explicit knowledge of the eigenfunctions from the geometry of $B(0,1)$, rather than extending the Agmon estimates to $\nu = 0$. We shall use the below Lemma \ref{lemma: monotonicity kummer function}, from \cite[Lemma 5.9]{vanlaere2026observability}. This Lemma was stated for $\nu > 0$ for consistency with the results of the paper, but holds verbatim for $\nu = 0$. 

\begin{lemma}\label{lemma: monotonicity kummer function}
    Let $\nu \geq 0$. There exists $\xi_0> 0$ such that, for every $\xi \geq \xi_0$, the function 
    \begin{align}\label{eqn: definition of kummer in lemma monotonicity}
        r \in [0,+\infty) \mapsto M\left(-\frac{\lambda_{\xi,\nu,0} - 2\xi(1+\nu)}{4\xi}, 1+\nu,\xi r^2 \right),
    \end{align}
    where $M(\cdot,\cdot,\cdot)$ is the Kummer function introduced in \eqref{eqn: definition kummer function}, is decreasing and strictly concave.
\end{lemma}

\begin{proof}[Proof of the lower bound of Theorem \ref{thm: main thm minimal time gamma = 1}]
Let $T < \operatorname{dist} \left(\omega_x, 0_{\mathbb{R}^{d_x}} \right)^2/(4(1+\nu))$. We prove that there exists no constant $C > 0$ independent of $n$ such that \eqref{eqn: uniform cost fourier} holds. We employ the dual observability point of view (see \textit{e.g.} \cite[Theorem 2.44]{coron2007control}). Namely, we equivalently prove that there exists no constant $C > 0$ independent of $n$ such that, for any $f_n^0 \in L^2(B(0,1))$, the solution $f_n$ of system \eqref{stm: Fourier heat adjoint system} with $\xi = \xi_n$ satisfies 
\begin{align}\label{eqn: obs ineq}
    \int_{B(0,1)} |f_n(T,x)|^2 \ dx \leq C \int_0^T \int_{\omega_x} |f_n(t,x)|^2 \ dx \ dt.
\end{align}
Assume that \eqref{eqn: obs ineq} holds. Let $a = \mathrm{dist}\left(\omega_x, 0_{\mathbb{R}^{d_x}} \right)$. Then, it holds on $\mathcal{A} = \{ x \in B(0,1): \ a \leq |x| \leq 1\}$ since $ \omega_x \subset \mathcal{A}$. In particular, there must exist $C > 0$ such that, for every $n \geq 1$, 
\begin{align}\label{eqn: obs ineq annulus}
    \int_{B(0,1)} |f_n(T,x)|^2 \ dx \leq C \int_0^T \int_{\mathcal{A}} |f_n(t,x)|^2 \ dx \ dt.
\end{align}
Setting $x = r\theta$, $r \in (0,1)$, $\theta \in \mathbb{S}^{d_x - 1}$, \eqref{eqn: obs ineq annulus} reads
\begin{align}\label{eqn: obs ineq polar}
    \int_0^1 \int_{\mathbb{S}^{d_x - 1}} |f_n(T,r\theta)|^2 r^{d_x - 1} \ d\theta \ dr \leq C \int_0^T \int_a^1 \int_{\mathbb{S}^{d_x - 1}} |f_n(t,r\theta)|^2 r^{d_x - 1} \ d\theta \ dr \ dt.
\end{align}

Recall, from Proposition \ref{prop: exact form eigenfunctions singular + cond eigenv} and the spherical decomposition of Section \ref{section: Spherical harmonic decomposition}, that setting
    \begin{align}
        \phi_{n,\nu,0}(r) = e^{- \xi_n r^2/2} r^{\frac{1}{2} + \nu} M \left(-\frac{\lambda_{n,\nu,0} - 2\xi_n (1 + \nu)}{4\xi_n}, 1 + \nu, \xi_n r^2 \right), 
    \end{align}
the functions $(r,\theta) \mapsto r^{-(d_x - 1)/2}\phi_{n,\nu,0}(r)Y_{0,1}(\theta)$ are eigenfunctions of $\mathcal{G}_n$. We test \eqref{eqn: obs ineq polar} against the solutions
\begin{align}
    (t,x) \in (0,T) \times B(0,1) \mapsto e^{-\lambda_{n,\nu,0} t }r^{-(d_x - 1)/2}\phi_{n,\nu,0}(r)Y_{0,1}(\theta).
\end{align}
For simplicity of the presentation, let us write
\begin{align}
    M_{n,\nu}(r) = M \left(-\frac{\lambda_{n,\nu,0} - 2\xi_n (1 + \nu)}{4\xi_n}, 1 + \nu, \xi_n r^2 \right).
\end{align}

We let $\xi_0 > 0$ be given by Lemma \ref{lemma: monotonicity kummer function}, and let $n_0$ be such that $\xi_n \geq \xi_0$ for $n \geq n_0$. In what follows we take $n \geq n_0$.\\

Using Lemma \ref{lemma: monotonicity kummer function}, since $M_{n,\nu}(0) = 1$ and $M_{n,\nu}(1) = 0$ by definition \eqref{eqn: definition kummer function} of $M$ and Proposition \ref{prop: exact form eigenfunctions singular + cond eigenv}, we have 
\begin{align}\label{eqn: bounds M n nu}
    1-r \leq M_{n,\nu}(r) \leq 1, \quad \text{for every } r \in (0,1).
\end{align}

On the one hand, using the lower bound in \eqref{eqn: bounds M n nu}, and that from \cite[Theorem 1.4]{romankummer2025} we have $\lambda_{n,\nu,0} = 2\xi_n(1+\nu) + \epsilon(n)$ with $\lim_n \epsilon(n) = 0$, we obtain
\begin{align*}
    \int_0^1 \int_{\mathbb{S}^{d_x - 1}} |f_n(T,r\theta)|^2 r^{d_x - 1} \ d\theta \ dr &= \int_0^1 e^{- 2\lambda_{n,\nu,0} T }e^{- \xi_n r^2} r^{1 + 2\nu} M_{n,\nu}(r)^2 \ dr \\
    &\geq \int_0^1 e^{-4\xi_n(1+\nu)T - 2\epsilon(n)T }e^{- \xi_n r^2} r^{1 + 2\nu} (1-r)^2 \ dr \\
    &\geq e^{-4\xi_n(1+\nu)T - 2\epsilon(n)T } e^{-1}\frac{1}{2(1+\nu)}\xi_n^{-(1+\nu)} \left(1- \xi_n^{-1/2} \right)^2.
\end{align*}
Thus, up to taking $n_0$ even larger than before, for every $n \geq n_0$ we have, for some $C_- > 0$ that does not depend on $n$, 
\begin{align}\label{eqn: lower bound obs annulus}
    \int_0^1 \int_{\mathbb{S}^{d_x - 1}} |f_n(T,r\theta)|^2 r^{d_x - 1} \ d\theta \ dr \geq C_- e^{-4\xi_n(1+\nu)T } \xi_n^{-(1+\nu)}.
\end{align}
On the other hand, using the upper bound in \eqref{eqn: bounds M n nu}, we have 
\begin{align*}
    \int_0^T \int_a^1 \int_{\mathbb{S}^{d_x - 1}} |f_n(t,r\theta)|^2 r^{d_x - 1} \ d\theta \ dr \ dt &\leq \int_0^T \int_a^1 e^{- 2\lambda_{n,\nu,0} t }e^{- \xi_n r^2} r^{1 + 2\nu} M_{n,\nu}(r)^2  \ dr \ dt \\
    &\leq T \int_a^1 e^{- \xi_n r^2}  \ dr \\
    &\leq  T (1-a) e^{-\xi_n a^2}.
\end{align*}

Hence, assuming that \eqref{eqn: obs ineq annulus} holds, it would imply, combining the above inequality with \eqref{eqn: lower bound obs annulus}, that there exists some $C > 0$, that may differ from the preceding ones, such that for any $n \geq n_0$ we have 
\begin{align}
    C \leq \xi_n^{1+\nu}e^{4\xi_n(1+\nu)T }e^{-\xi_n a^2}.
\end{align}
This is impossible for any $T < a^2/(4(1+\nu))$, as in this case
\begin{align}
    \xi_n^{1+\nu} e^{4\xi_n(1+\nu)T }e^{-\xi_n a^2} \to 0, \quad n \to +\infty.
\end{align}
This concludes the proof. 
\end{proof}

\section{Proof of the null-controllability results in the weakly degenerate case: Theorem \ref{thm: main thm minimal time gamma < 1}}\label{section: Proof of null-controllability results in the weakly degenerate case}

In this section, we prove Theorem \ref{thm: uniform control annulus gamma < 1}, and then Theorem \ref{thm: main thm minimal time gamma < 1}. \\

\subsection{Preliminaries on the transmutation method: from heat to wave}

We recall here the transmutation method, to construct from a solution of the heat equation \eqref{stm: radial heat adjoint system}, a solution of the wave equation \eqref{stm: radial wave adjoint system}. The following proposition is from \cite[Proposition 3.1]{ervedoza2011sharp} (see also \cite[Proposition 1]{ervedoza2011observability}, and \cite[Proposition 4.3]{beauchard20152d}).

\begin{proposition}\label{prop: estimate transmutation kernel}
    Let $T > 0$ and $L > 0$. For every $\alpha > 3L^2$, there exists a unique solution of
    \begin{align}
        \left\{ \begin{array}{lcll}
            \partial_t k_T(t,s) + \partial_s^2 k_T(t,s) &=& 0, & (t,s) \in (0,T) \times (-L,L), \\[6pt]
            k_T(0,s) = k_T(T,s) &=& 0, & s \in (-L,L), \\[6pt]
            k_T(t,0) &=& 0, & t \in (0,T), \\[6pt]
            \partial_s k_T(t,0) &=& \displaystyle \exp\left( -\frac{\alpha T}{t(T-t)} \right), & t \in (0,T).
        \end{array} \right.
    \end{align}
    Moreover, for every $(t,s) \in (0,T) \times (-L,L)$, this solution satisfies
    \begin{align*}
        |k_T(t,s)| &\leq |s|\exp\left( \frac{2}{\min(t,T-t)}\left( s^2 - \frac{\alpha}{3} \right) \right), \\
        |\partial_s k_T(t,s)| &\leq \exp\left( \frac{2}{\min(t,T-t)}\left( s^2 - \frac{\alpha}{3} \right) \right).
    \end{align*}
\end{proposition}

The following proposition is from \cite[Theorem 2.1]{ervedoza2011sharp} (see also \cite[Proposition 2]{ervedoza2011observability}, and \cite[Proposition 4.4]{beauchard20152d}). 

\begin{proposition}\label{prop: transmutation}
    Let $z_0 \in L^2(0,1)$, and let $z$ be the associated solution of \eqref{stm: radial heat adjoint system}. We define
    \begin{align}
        w(s,r) := \int_0^T k_T(t,s)z(t,r) \ dt, \quad (s,r) \in (-L,L) \times (0,1).
    \end{align}
    Then, $w$ is the unique solution of \eqref{stm: radial wave adjoint system} with initial conditions
    \begin{align}
        w_0(r) := 0, \quad w_1(r) := \int_0^T \exp\left( -\frac{\alpha T}{t(T-t)} \right)z(t,r) \ dt.
    \end{align}
\end{proposition}

\subsection{From radial uniform observability of solutions of the heat equation to uniform null-controllability from annuli: proof of Theorem \ref{thm: uniform control annulus gamma < 1}}\label{section: From radial uniform observability of solutions of the heat equation to uniform null-controllability from annuli}

Building on Theorem \ref{thm: radial wave observability} and the transmutation method, we start by proving a small time observability estimate for the solutions of \eqref{stm: radial heat adjoint system} (see Proposition \ref{prop: radial small time cost observability}), from which we then easily deduce the sought uniform observability with respect to $\nu_\ell$ (see Proposition \ref{prop: radial uniform small time cost observability}). Then, we prove the dual uniform observability inequality equivalent of Theorem \ref{thm: uniform control annulus gamma < 1}.

\subsubsection{Radial uniform observability for the solutions of the heat equation}

\begin{proposition}\label{prop: radial small time cost observability}
    For every $0 \leq a < b \leq 1$, there exists $C > 0$ such that, for every $T,\xi > 0$, $\ell \in \mathbb{N}$, and every $z_0 \in L^2(0,1)$, the solution $z$ of system \eqref{stm: radial heat adjoint system} with $\mu^2 = \nu_\ell^2$ satisfies
    \begin{align}\label{eqn: cost of obs heat}
        \begin{aligned}
            \| z(T)\|_{L^2(0,1)}^2 \leq \frac{C}{T} \exp\left( \frac{C}{T} \right) \exp\left( C \left( \xi + \nu_\ell \ln(2+\xi) \right) \right) \|z\|_{L^2( (0,T) \times (a,b))}^2,
        \end{aligned}
    \end{align}
\end{proposition}

\begin{proof}
Let $z_0 \in L^2(0,1)$, and denote by $z$ the associated solution of \eqref{stm: radial heat adjoint system}. Using Proposition \ref{prop: transmutation}, $w$ defined by 
    \begin{align}
        w(s,r) = \int_0^T k_T(t,s)z(t,r) \ dt, \quad (s,r) \in (-S,S) \times (0,1),
    \end{align}
solves \eqref{stm: radial wave adjoint system} with initial states
\begin{align}\label{eqn: initial states heat to wave}
    w_0(r) = 0, \quad w_1(r) = \int_0^T \exp\left( -\frac{\alpha T}{t(T-t)} \right)z(t,r) \ dt, 
\end{align}
and satisfies, thanks to Theorem \ref{thm: radial wave observability},
\begin{align}\label{eqn: wave observability annulus i}
        \|w_1\|_{L^2(0,1)}^2 \leq C \exp\left( C \left( \xi + \nu_\ell \ln(2+\xi) \right) \right) \int_a^b \int_{-S}^S w_s(s,r)^2 + w(s,r)^2  \ ds \ dr.
    \end{align}
First, we estimate from above 
\begin{align*}
     \int_a^b \int_{-S}^S w_s(s,r)^2 &+ w(s,r)^2  \ ds \ dr \\
    &=  \int_a^b \int_{-S}^S \left( \int_0^T \partial_s k_T(t,s)z(t,r) \ dt\right)^2 + \left( \int_0^T k_T(t,s)z(t,r) \ dt\right)^2  \ ds \ dr.
\end{align*}

Using Cauchy-Schwarz, we have 
\begin{align*}
    \begin{aligned}
        & \int_a^b \int_{-S}^S w_s(s,r)^2 + w(s,r)^2  \ ds \ dr \\
        &\leq \left( \|\partial_s k_T\|_{L^2((0,T) \times (-S,S))}^2 + \|k_T\|_{L^2((0,T) \times (-S,S))}^2 \right)  \int_a^b \int_0^T |z(t,r)|^2 \ dt \ dr
    \end{aligned}
\end{align*}
and thus, applying the estimates on $k_T$ of Proposition \ref{prop: estimate transmutation kernel}, it follows that 
\begin{align}\label{eqn: upper bound energy}
     \int_a^b \int_{-S}^S w_s(s,r)^2 + w(s,r)^2  \ ds \ dr \leq 2ST\left( 1 + S^2\right)  \int_a^b \int_0^T |z(t,r)|^2 \ dt \ dr.
\end{align}

Second, we bound from below $\|w_1\|_{L^2(0,1)}^2$ given in \eqref{eqn: initial states heat to wave}. To do so, we decompose the solution $z$ of system \eqref{stm: radial heat adjoint system} into an orthonormal basis of $L^2(0,1)$ consisting of eigenfunctions of $G_{\xi,\nu_\ell}$ (see \eqref{eqn: def G radial} and Section \ref{section: Uniform spectral estimates}), denoting $z_0 = \sum_k z_{0,k} \varphi_{\xi,\ell,k}$,
\begin{align*}
    \|w_1 \|_{L^2(0,1)}^2 &= \int_0^1 \left| \int_0^T \exp\left( -\frac{\alpha T}{t(T-t)} \right)z(t,r) \ dt \right|^2 \ dr \\
    &= \int_0^1 \left| \int_0^T \exp\left( -\frac{\alpha T}{t(T-t)} \right) \sum_{k \geq 0}  \exp \left(-\lambda_{\xi,\ell,k} t \right) z_{0,k} \varphi_{\xi,\ell,k}(r)  \ dt \right|^2 \ dr \\
    &= \sum_{k \geq 0} |z_{0,k}|^2 \left| \int_0^T \exp\left( -\frac{\alpha T}{t(T-t)} \right) \exp \left(-\lambda_{\xi,\ell,k} t \right)  \ dt \right|^2 \\
    &\geq \left| \int_0^T \exp\left( -\frac{\alpha T}{t(T-t)} \right) \ dt \right|^2 \sum_{k \geq 0} |z_{0,k}|^2  \exp \left(-2\lambda_{\xi,\ell,k} T \right). 
\end{align*}

Hence, 
\begin{align}\label{eqn: lower bound initial state wave}
    \|w_1 \|_{L^2(0,1)}^2 \geq \left| \int_0^T \exp\left( -\frac{\alpha T}{t(T-t)} \right) \ dt \right|^2 \| z(T)\|_{L^2(0,1)}^2. 
\end{align}

Thus, substituting \eqref{eqn: upper bound energy} and \eqref{eqn: lower bound initial state wave} into \eqref{eqn: wave observability annulus i}, we obtain, for some new constant $ C > 0$,
\begin{align}\label{eqn: cost of obs heat i}
        \begin{aligned}
            \left| \int_0^T \exp\left( -\frac{\alpha T}{t(T-t)} \right) \ dt \right|^2 \| z(T)\|_{L^2(0,1)}^2 
            &\leq CT\exp\left( C \left( \xi + \nu_\ell \ln(2+\xi) \right) \right) \|z\|_{L^2( (0,T) \times (a,b))}^2,
        \end{aligned}
    \end{align}

which is the sought estimate since
\begin{align*}
    \int_0^T \exp\left( -\frac{\alpha T}{t(T-t)} \right) \ dt \geq \frac{T}{2} \exp \left( -\frac{16\alpha}{3T} \right). 
\end{align*}
\end{proof}

We now deduce from Proposition \ref{prop: radial small time cost observability} a uniform cost of observability with respect to $\ell$ for the solutions of \eqref{stm: radial heat adjoint system}.

\begin{proposition}\label{prop: radial uniform small time cost observability}
    For every $0 \leq a < b \leq 1$, there exists $C > 0$ such that, for every $T,\xi > 0$, $\ell \in \mathbb{N}$, and every $z_0 \in L^2(0,1)$, the solution $z$ of system \eqref{stm: radial heat adjoint system} with $\mu^2 = \nu_\ell^2$ satisfies
    \begin{align}
         \|z(T)\|_{L^2(0,1)}^2 \leq \frac{C}{T} \exp\left( \frac{C}{T} \right) \exp\left( C\xi + \frac{C^2\ln(2 + \xi)^2}{4T} \right) \|z\|_{L^2( (0,T) \times (a,b))}^2.
    \end{align}
\end{proposition}

\begin{proof}
    Let $z_0 \in L^2(0,1)$, and $z$ the associated solution of system \eqref{stm: radial heat adjoint system}. By Proposition \ref{prop: lower bound eigenvalues any gamma}, we have 
    \begin{align*}
        \|z(T)\|_{L^2(0,1)}^2 \leq \exp\left(-\nu_\ell^2 T\right)\|z(T/2)\|_{L^2(0,1)}^2.
    \end{align*}
    Applying Proposition \ref{prop: radial small time cost observability} to the right-hand side above on $(0,T/2)$, we obtain 
    \begin{align}\label{eqn: unif radial observability i}
        \|z(T)\|_{L^2(0,1)}^2 \leq  \exp\left(-\nu_\ell^2 T \right) \frac{2C}{T} \exp\left( \frac{2C}{T} \right) \exp\left( C \left( \xi + \nu_\ell \ln(2+\xi) \right) \right) \|z\|_{L^2( (0,T) \times (a,b))}^2.
    \end{align}
    We have that, for any $C,\xi,T > 0$, 
    \begin{align*}
         \sup_{\mu \geq 0}  C\ln(2+\xi)\mu - \mu^2T = \frac{C^2\ln(2 + \xi)^2}{4T},
    \end{align*}
    and thus, substituting this into \eqref{eqn: unif radial observability i}, 
    \begin{align*}
        \|z(T)\|_{L^2(0,1)}^2 \leq \frac{2C}{T} \exp\left( \frac{2C}{T} \right) \exp\left( C\xi + \frac{C^2\ln(2 + \xi)^2}{4T} \right) \|z\|_{L^2( (0,T) \times (a,b))}^2,
    \end{align*}
    which concludes the proof.
\end{proof}

\subsubsection{Proof of Theorem \ref{thm: uniform control annulus gamma < 1}}\label{section: uniform observability from annuli}

By duality, Theorem \ref{thm: uniform control annulus gamma < 1} is equivalent to the following theorem. 

\begin{theorem}
    Let $\gamma \in (0,1)$. For every $0 \leq a < b \leq 1$, setting $\mathcal{A} := \{x \in B(0,1), \ a < |x| < b\}$, there exists $C > 0$ such that, for every $T > 0$, for every $n\geq 1$ such that $\xi_n > 0$, for every $f_n^0 \in L^2(B(0,1))$, the solutions $f_n$ of systems \eqref{stm: Fourier heat adjoint system} with $\xi = \xi_n$ satisfy
    \begin{align}
        \|f_n(T)\|_{L^2(B(0,1))}^2 \leq C \exp \left( C T^{\frac{\gamma + 1}{ \gamma - 1}} \right) \|f_n\|_{L^2( (0,T) \times \mathcal{A})}^2.
    \end{align}
\end{theorem}

\begin{proof}
Let $f_n^0 \in L^2(B(0,1))$. Recall that by the decomposition in spherical harmonics, we have that the associated solution $f_n$ satisfies, 
\begin{align*}
    f_n(t,x) = r^{-(d_x - 1)/2} \sum_{\ell \geq 0} \sum_{m = 1}^{m_\ell} z_{n,\ell,m}(t,r) Y_{\ell,m}(\theta), \quad x = r\theta,
\end{align*}
with $z_{n,\ell,m}$ solution of system \eqref{stm: radial heat adjoint system} with $\xi = \xi_n$ and $\mu^2 = \nu_\ell^2$. Using Proposition \ref{prop: lower bound eigenvalues any gamma}, we have for some $C_0 > 0$, 
\begin{align}\label{eqn: dissip gamma 0,1 proof thm}
    \|f_n(T)\|_{L^2(B(0,1))}^2 \leq e^{-C_0|\xi_n|^{2/(\gamma + 1)}T} \|f_n(T/2)\|_{L^2(B(0,1))}^2.
\end{align}

By Parseval's identity, we have
\begin{align*}
    \|f_n(T/2)\|_{L^2(B(0,1))}^2 = \sum_{\ell \geq 0} \sum_{m = 1}^{m_\ell} \|z_{n,\ell,m}(T/2)\|_{L^2(0,1)}^2.
\end{align*}
Applying Proposition \ref{prop: radial uniform small time cost observability}, summing back and using that $(2C_1/T) \exp(2C_1/T) \leq \exp(4C_1/T)$, we obtain
\begin{align}\label{eqn: unif observability heat annulus gamma < 1 0}
    \|f_n(T/2)\|_{L^2(B(0,1))}^2 \leq  \exp\left( C_1\xi_n + \frac{8C_1 + C_1^2\ln(2 + \xi_n)^2}{2T} \right)  \|f_n\|_{L^2( (0,T) \times \mathcal{A})}^2.
\end{align}

Hence, combining \eqref{eqn: dissip gamma 0,1 proof thm} and \eqref{eqn: unif observability heat annulus gamma < 1 0}, we have
\begin{align}\label{eqn: unif observability heat annulus gamma < 1 i}
    \|f_n(T)\|_{L^2(B(0,1))}^2 \leq \exp\left( C_1\xi_n + \frac{8C_1 + C_1^2\ln(2 + \xi_n)^2}{2T} - C_0T|\xi_n|^{2/(\gamma + 1)} \right)  \|f_n\|_{L^2( (0,T) \times \mathcal{A})}^2.
\end{align}
Now, for some $C_2,C_3 > 0$, for every $T > 0$, 
\begin{align*}
    \sup_{\xi_n \geq 0} C_1\xi_n - \frac{C_0T}{2}\xi_n^{2/(\gamma + 1)} = C_2 T^{(\gamma + 1)/(\gamma - 1)},
\end{align*}
and 
\begin{align*}
    \sup_{\xi_n \geq 0} \frac{8C_1 + C_1^2\ln(2 + \xi_n)^2}{2T} - \frac{C_0T}{2}|\xi_n|^{2/(\gamma + 1)} \leq \frac{C_3}{T} \left( 1 + \ln \left( 1 + \frac{1}{T} \right) \right)^2.
\end{align*}
Hence, for some $C_4 > 0$, for every $T > 0$, since $\gamma < 1$,
\begin{align*}
    \sup_{\xi_n \geq 0} C_1\xi_n + \frac{8C_1 + C_1^2\ln(2 + \xi_n)^2}{2T} - C_0T|\xi_n|^{2/(\gamma + 1)} &\leq C_2 T^{(\gamma + 1)/(\gamma - 1)} + \frac{C_3}{T} \left( 1 + \ln \left( 1 + \frac{1}{T} \right) \right)^2 \\
    &\leq C_4 \left( 1 + T^{(\gamma + 1)/(\gamma - 1)} \right).
\end{align*}

Substituting this into \eqref{eqn: unif observability heat annulus gamma < 1 i} gives the sought result. 
\end{proof}

\subsection{Small time null-controllability: Proof of Theorem \ref{thm: main thm minimal time gamma < 1}}

We are now in position to prove Theorem \ref{thm: main thm minimal time gamma < 1}. Its proof relies on reducing the control zones using spectral inequalities and the Lebeau-Robbiano method, in the almost exact same manner as for the proof of Theorem \ref{thm: main thm minimal time gamma = 1}. For this reason, we only sketch the proof, and we emphasize where are the differences. In what follows, we set
\begin{align}\label{def: beta}
    \beta := \frac{1+\gamma}{1-\gamma} > 1.
\end{align}

\begin{proof}[Proof of Theorem \ref{thm: main thm minimal time gamma < 1}]
Let $\omega \subset B(0,1) \times \Omega_y$ be a nonempty open set. Since it contains a set of the form $\omega_x \times \omega_y$, with $\omega_x \subset B(0,1)$ and $\omega_y \subset \Omega_y$ nonempty open sets, it suffices to prove the null-controllability result in time $T > 0$ from $\omega_x \times \omega_y$. Moreover, it suffices to consider $T \in (0,1]$, since for $T > 1$ one can take the control given on $(0,1)$ and extend it by zero. \\

Let $0 < a < b < 1$, and $\Gamma \subset \mathbb{S}^{d_x-1}$ a nonempty open set, such that $\mathcal{A}_\Gamma := \{x = r\theta \in B(0,1): \ a < r < b, \ \theta \in \Gamma\} \subset \omega_x$. Set also $\mathcal{A} := \{x = r\theta \in B(0,1): \ a < r < b\}$. Let $f_0 \in L^2(B(0,1) \times \Omega_y)$, and let $f_n^0 \in L^2(B(0,1))$ be its Fourier coefficients. \\

\textit{First step: angular reduction and control from $\mathcal{A}_\Gamma \times \Omega_y \subset \omega_x \times \Omega_y$.} 
Applying Theorem \ref{thm: uniform control annulus gamma < 1}, there exists $C > 0$ such that, for every $\tau > 0$, $n \geq 1$ such that $\xi_n > 0$, and $f_n^0 \in L^2(B(0,1))$, there exists a control $u_n \in L^2((0,\tau) \times \mathcal{A})$ satisfying
\begin{align}
        \begin{aligned}
            &\|u_n\|_{L^2((0,\tau) \times \mathcal{A})} \leq C \exp\left( C \tau^{-\beta} \right) \|f_n^0\|_{L^2(B(0,1))}, \\
            &f_n(\tau) = 0,
        \end{aligned}
    \end{align}
where $f_n$ is the solution of system \eqref{stm : Fourier heat control system} with $\xi = \xi_n$, $\mathcal{U} = \mathcal{A}$, and initial state $f_n^0$.\\

We apply almost exactly the first and second steps of the proof of the upper bound of Theorem \ref{thm: main thm minimal time gamma = 1}, with $\beta$ given in \eqref{def: beta}, with the difference that there is no use of the dissipation estimate \eqref{eqn: dissipation estimate LR proof}, and that $\mathcal{C}_n = C(\Gamma)^{1/2}$, where $C(\Gamma)$ is the constant appearing in the spectral inequality on the sphere. The application of the Lebeau-Robbiano method does not change the power of $\tau$ in the control cost, but only the constants. We thus obtain that there exists $C_0 > 0$ such that, for every $\tau \in (0,1]$, $n \geq 1$ such that $\xi_n > 0$, and $f_n^0 \in L^2(B(0,1))$, there exists a control $v_n \in L^2((0,\tau) \times \mathcal{A}_\Gamma)$ satisfying
\begin{align}\label{eqn: cost control strip gamma < 1 proof}
        \begin{aligned}
            &\|v_n\|_{L^2((0,\tau) \times \mathcal{A}_\Gamma)} \leq C_0 \exp\left( C_0 \tau^{-\beta} \right) \|f_n^0\|_{L^2(B(0,1))}, \\
            &f_n(\tau) = 0,
        \end{aligned}
\end{align}
where now $\mathcal{U} = \mathcal{A}_\Gamma$. If $\partial \Omega_y = \emptyset$, then $\xi_1 = 0$, and the corresponding Fourier component is controlled from $\mathcal{A}_\Gamma$ with a cost of the same form thanks to \cite{ervedoza2008control}. Hence, up to increasing $C_0$, \eqref{eqn: cost control strip gamma < 1 proof} holds for every $n \geq 1$. \\

\textit{Second step: vertical reduction in the $y$-direction, and control from $\mathcal{A}_\Gamma \times \omega_y \subset \omega_x \times \omega_y$.} We apply once again the strategy of the first step, now in the $y$-direction, with the following differences. \\

The spectral projector, for $N \geq 1$, is now $\Pi_N$ onto the eigenspaces $\mathcal{V}_N := \operatorname{Span} \left\{ \phi_n, \ n \leq N \right\}$, where the $\phi_n$ are the eigenfunctions of $-\Delta_y$. By \cite{jerison1999nodal}, there exists $C(\omega_y) > 0$ such that, for every $N \geq 1$ and $g \in \mathcal{V}_N$,
    \begin{align*}
        \|g\|_{L^2(\Omega_y)}^2 \leq C(\omega_y) e^{C(\omega_y)\xi_N} \|g\|_{L^2(\omega_y)}^2.
    \end{align*}
Hence, proceeding as for the angular reduction of the first step, for every $\tau \in (0,1]$, $N \geq 1$ and $f_0 \in L^2(\Omega)$, there exists a control $u \in L^2((0,\tau) \times \Omega)$ supported in $(0,\tau) \times \mathcal{A}_\Gamma \times \omega_y$ such that the associated solution $f$ of system \eqref{stm : control system} satisfies
\begin{align*}
    \Pi_N f(\tau) = 0, \qquad \|u\|_{L^2((0,\tau) \times \omega)} \leq C \exp\left( C \tau^{-\beta} + C \xi_N \right) \|f_0\|_{L^2(\Omega)},
\end{align*}
where $C > 0$ does not depend on $N$. Also, by Proposition \ref{prop: lower bound eigenvalues any gamma}, there exists $C_* > 0$ such that, for every $N \geq 1$, the solutions $f$ of system \eqref{stm : control system} without control, and initial state satisfying $\Pi_N f(s) = 0$, satisfy
\begin{align}\label{eqn: high frequency dissipation y}
    \|f(t)\|_{L^2(\Omega)} \leq \exp\left(-C_* \xi_{N+1}^{2/(\gamma+1)}(t-s) \right) \|f(s)\|_{L^2(\Omega)}, \qquad 0 \leq s \leq t.
\end{align}

The Lebeau-Robbiano construction of \cite[Section 6.2]{le2012carleman} (see also \cite[proof of Theorem 1 (1)]{beauchard2014null}) then applies, as the cost of control of low-frequencies is absorbed by the dissipation of the high-frequencies since $\gamma \in (0,1)$ (as $2/(\gamma + 1) > 1$). Notice that when $\gamma = 1$, the dissipation of the high frequencies no longer absorbs the spectral cost in any time, and the reduction in the $y$-direction fails. This is consistent with the geometric condition results of \cite{LISSYprolate2025, duprez2020control, darde2023null, CRMATH_2017__355_12_1215_0}.\\

\textit{Conclusion.} We thus obtain from the second step that system \eqref{stm : control system} is null-controllable in any time $T \in (0,1]$, and thus in any time $T > 0$, from $\mathcal{A}_\Gamma \times \omega_y \subset \omega_x \times \omega_y \subset \omega$. This concludes the proof. 
\end{proof}

\subsection*{Acknowledgments}

The author would like to thank S. Ervedoza for suggesting the use of the transmutation method, and Mathilda Trabut for several discussions and comments. The author was supported by the project of Consolidación Investigadora CNS2024-154725 funded by MICIU/AEI/10.13039/501100011033.

\subsection*{Tool and computational resource disclosure}

The author acknowledges the use of AI: 
\begin{itemize}
    \item During the investigation phase, AI was used only for bibliographical research, and to assist in understanding the literature. No proof or idea presented in this work has been AI-generated. Every proofs and ideas were obtained and independently checked by the author. This work was entirely written by the author, without any use of AI.
    \item The author also acknowledges the use of AI at the final stage to verify typos and certain computations in a preliminary version. In particular, it spotted a mistake at the end of the proof of Theorem \ref{thm: gap uniform}, and computational errors in the application of the Lebeau-Robbiano strategy in the proof of Theorem \ref{thm: main thm minimal time gamma = 1} and in the tracking of parameters and constants in the computations of Section \ref{section: radial observability of waves}. In each case, it suggested corrections, which were independently verified by the author.
\end{itemize}

The author bears sole responsibility for the content of the present work.

\printbibliography

\end{document}

%% file: structurepaper.tex
\begin{tikzpicture}[tikzfig]
	\begin{pgfonlayer}{nodelayer}
		\node [style=none] (5) at (-19.5, 1.5) {};
		\node [style=none] (6) at (-13.5, 1.5) {};
		\node [style=none] (7) at (-13.5, -1.5) {};
		\node [style=none] (8) at (-19.5, -1.5) {};
		\node [style=none] (37) at (0, 14) {};
		\node [style=none] (38) at (0, -2.75) {};
		\node [style=none] (39) at (12.5, 14) {};
		\node [style=none] (40) at (12.5, -2.75) {};
		\node [style=none] (41) at (-8.5, 1.5) {};
		\node [style=none] (42) at (-2.5, 1.5) {};
		\node [style=none] (43) at (-2.5, -1.5) {};
		\node [style=none] (44) at (-8.5, -1.5) {};
		\node [style=none] (45) at (-19.5, 8.5) {};
		\node [style=none] (46) at (-8.5, 8.5) {};
		\node [style=none] (47) at (-8.5, 6.5) {};
		\node [style=none] (48) at (-19.5, 6.5) {};
		\node [style=none] (49) at (-19.5, 5) {};
		\node [style=none] (50) at (-8.5, 5) {};
		\node [style=none] (51) at (-8.5, 3) {};
		\node [style=none] (52) at (-19.5, 3) {};
		\node [style=none] (53) at (15, 1.75) {};
		\node [style=none] (54) at (21, 1.75) {};
		\node [style=none] (55) at (21, -1.75) {};
		\node [style=none] (56) at (15, -1.75) {};
		\node [style=none] (57) at (2.5, 1.75) {};
		\node [style=none] (58) at (10, 1.75) {};
		\node [style=none] (59) at (10, -1.75) {};
		\node [style=none] (60) at (2.5, -1.75) {};
		\node [style=none] (61) at (2.5, 7.5) {};
		\node [style=none] (62) at (10, 7.5) {};
		\node [style=none] (63) at (10, 4) {};
		\node [style=none] (64) at (2.5, 4) {};
		\node [style=none] (65) at (15, 7.5) {};
		\node [style=none] (66) at (21, 7.5) {};
		\node [style=none] (67) at (21, 4) {};
		\node [style=none] (68) at (15, 4) {};
		\node [style=none] (69) at (-12.75, 0) {};
		\node [style=none] (70) at (-9.25, 0) {};
		\node [style=none] (71) at (-1.75, 0) {};
		\node [style=none] (72) at (1.75, 0) {};
		\node [style=none] (73) at (-7.75, 4) {};
		\node [style=none] (74) at (-7.75, 7.5) {};
		\node [style=none] (75) at (-5.5, 7.5) {};
		\node [style=none] (76) at (10.75, 0) {};
		\node [style=none] (77) at (14.25, 0) {};
		\node [style=none] (78) at (10.75, 5.75) {};
		\node [style=none] (79) at (14.25, 5.75) {};
		\node [style=none] (80) at (-5.5, 4) {};
		\node [style=none] (81) at (-5.5, 5.75) {};
		\node [style=none] (82) at (1.75, 5.75) {};
        \node [style=none] (83) at (-19.5, 2.375) {};
		\node [style=none] (84) at (21, 2.375) {};
		\node [style=none] (85) at (-19.5, 12.5) {};
		\node [style=none] (86) at (-8.5, 12.5) {};
		\node [style=none] (87) at (-8.5, 10.5) {};
		\node [style=none] (88) at (-19.5, 10.5) {};
		\node [style=none] (89) at (15, 13.25) {};
		\node [style=none] (90) at (21, 13.25) {};
		\node [style=none] (91) at (21, 9.75) {};
		\node [style=none] (92) at (15, 9.75) {};
		\node [style=none] (93) at (-7.75, 11.5) {};
		\node [style=none] (94) at (14.25, 11.5) {};
		\node [style=none, font=\small\scshape] at (-11, 14) {Radial systems};
		\node [style=none, font=\small\scshape] at (6.25, 14) {Fourier systems};
		\node [style=none, font=\small\scshape] at (18, 14) {Grushin systems};
		\node [style=none, align=center, font=\scriptsize] at (-14, 11.5) {Explicit eigenfunctions of 1D singular\\Grushin operators (Kummer functions)\\{\color{black!55}Proposition \ref{prop: exact form eigenfunctions singular + cond eigenv}, Lemma \ref{lemma: monotonicity kummer function}}};
		\node [style=none, align=center, font=\scriptsize] at (18, 11.5) {Lower bound on\\ the minimal time\\ of null-controllability\\ from vertical strips\\when $\gamma = 1$\\{\color{black!55}Theorem \ref{thm: main thm minimal time gamma = 1}}};
		\node [style=none, align=center, font=\scriptsize\itshape, text=black!70, fill=white, inner sep=1pt, anchor=south] at (3.25, 11.75) {Testing the dual observability inequality\\on a sequence of solutions consisting of eigenfunctions};
		\node [style=none, align=center, font=\scriptsize] at (-14, 7.5) {Uniform spectral gap of 1D singular\\Grushin operators\\{\color{black!55}Theorem \ref{thm: gap uniform}}};
		\node [style=none, align=center, font=\scriptsize] at (-14, 4) {Localization estimates of eigenfunctions\\of 1D singular Grushin operators\\{\color{black!55}Theorem \ref{thm: localization estimates eigenfunctions}}};
		\node [style=none, align=center, font=\scriptsize] at (6.25, 5.75) {Cost of small-time\\null-controllability of the\\Fourier heat systems\\from annuli when $\gamma = 1$\\{\color{black!55}Theorem \ref{thm: uniform control annulus gamma = 1}}};
		\node [style=none, align=center, font=\scriptsize] at (18, 5.75) {Upper bound on\\ the minimal time\\ of null-controllability\\ from vertical strips\\when $\gamma = 1$\\{\color{black!55}Theorem \ref{thm: main thm minimal time gamma = 1}}};
		\node [style=none, align=center, font=\scriptsize] at (-16.5, 0) {Observability\\of radial waves\\{\color{black!55}Theorem \ref{thm: radial wave observability}}};
		\node [style=none, align=center, font=\scriptsize] at (-5.5, 0) {Uniform observability\\of the radial heat\\{\color{black!55}Proposition \ref{prop: radial uniform small time cost observability}}};
		\node [style=none, align=center, font=\scriptsize] at (6.25, 0) {Cost of small-time\\null-controllability of the\\Fourier heat systems\\from annuli when $\gamma \in (0,1)$\\{\color{black!55}Theorem \ref{thm: uniform control annulus gamma < 1}}};
		\node [style=none, align=center, font=\scriptsize] at (18, 0) {Null-controllability in \\ arbitrarily small times\\ from open sets\\when $\gamma \in (0,1)$\\{\color{black!55}Theorem \ref{thm: main thm minimal time gamma < 1}}};
		\node [style=none, align=center, font=\scriptsize\itshape, text=black!70, fill=white, inner sep=0pt, anchor=south] at (-2, 6) {Moments method and\\Parseval's identity};
		\node [style=none, align=center, font=\scriptsize\itshape, text=black!70, fill=white, inner sep=0pt, anchor=south] at (12.5, 6) {Angular\\reduction and};
		\node [style=none, align=center, font=\scriptsize\itshape, text=black!70, fill=white, inner sep=0pt, anchor=north] at (12.5, 5.5) {Lebeau-Robbiano\\method};
		\node [style=none, align=center, font=\scriptsize\itshape, text=black!70, anchor=south] at (-11, 0.25) {Transmutation\\method};
		\node [style=none, align=center, font=\scriptsize\itshape, text=black!70, anchor=north] at (-11, -0.25) {and dissipation\\estimates};
		\node [style=none, align=center, font=\scriptsize\itshape, text=black!70, fill=white, inner sep=0pt, anchor=south] at (0, 0.25) {Parseval's\\identity and};
		\node [style=none, align=center, font=\scriptsize\itshape, text=black!70, fill=white, inner sep=0pt, anchor=north] at (0, -0.25) {duality};
		\node [style=none, align=center, font=\scriptsize\itshape, text=black!70, fill=white, inner sep=0pt, anchor=south] at (12.5, 0.25) {Angular and\\ vertical\\ reductions and};
		\node [style=none, align=center, font=\scriptsize\itshape, text=black!70, fill=white, inner sep=0pt, anchor=north] at (12.5, -0.25) {Lebeau-Robbiano\\method};
	\end{pgfonlayer}
	\begin{pgfonlayer}{edgelayer}
		\draw [dotted, black!45] (37.center) to (38.center);
		\draw [dotted, black!45] (39.center) to (40.center);
		\draw [bend left, looseness=0.25] (45.center) to (46.center);
		\draw [bend left=90, looseness=0.25] (46.center) to (47.center);
		\draw [bend left, looseness=0.25] (47.center) to (48.center);
		\draw [bend left=90, looseness=0.25] (48.center) to (45.center);
		\draw [bend left, looseness=0.25] (49.center) to (50.center);
		\draw [bend left=90, looseness=0.25] (50.center) to (51.center);
		\draw [bend left, looseness=0.25] (51.center) to (52.center);
		\draw [bend left=90, looseness=0.25] (52.center) to (49.center);
		\draw [bend left, looseness=0.25] (61.center) to (62.center);
		\draw [bend left=90, looseness=0.25] (62.center) to (63.center);
		\draw [bend left, looseness=0.25] (63.center) to (64.center);
		\draw [bend left=90, looseness=0.25] (64.center) to (61.center);
		\draw [bend left, looseness=0.25] (65.center) to (66.center);
		\draw [bend left=90, looseness=0.25] (66.center) to (67.center);
		\draw [bend left, looseness=0.25] (67.center) to (68.center);
		\draw [bend left=90, looseness=0.25] (68.center) to (65.center);
		\draw [bend left, looseness=0.25] (5.center) to (6.center);
		\draw [bend left=90, looseness=0.25] (6.center) to (7.center);
		\draw [bend left, looseness=0.25] (7.center) to (8.center);
		\draw [bend left=90, looseness=0.25] (8.center) to (5.center);
		\draw [bend left, looseness=0.25] (41.center) to (42.center);
		\draw [bend left=90, looseness=0.25] (42.center) to (43.center);
		\draw [bend left, looseness=0.25] (43.center) to (44.center);
		\draw [bend left=90, looseness=0.25] (44.center) to (41.center);
		\draw [bend left, looseness=0.25] (57.center) to (58.center);
		\draw [bend left=90, looseness=0.25] (58.center) to (59.center);
		\draw [bend left, looseness=0.25] (59.center) to (60.center);
		\draw [bend left=90, looseness=0.25] (60.center) to (57.center);
		\draw [bend left, looseness=0.25] (53.center) to (54.center);
		\draw [bend left=90, looseness=0.25] (54.center) to (55.center);
		\draw [bend left, looseness=0.25] (55.center) to (56.center);
		\draw [bend left=90, looseness=0.25] (56.center) to (53.center);
		\draw [bend left, looseness=0.25] (85.center) to (86.center);
		\draw [bend left=90, looseness=0.25] (86.center) to (87.center);
		\draw [bend left, looseness=0.25] (87.center) to (88.center);
		\draw [bend left=90, looseness=0.25] (88.center) to (85.center);
		\draw [bend left, looseness=0.25] (89.center) to (90.center);
		\draw [bend left=90, looseness=0.25] (90.center) to (91.center);
		\draw [bend left, looseness=0.25] (91.center) to (92.center);
		\draw [bend left=90, looseness=0.25] (92.center) to (89.center);
		\draw [-stealth] (93.center) to (94.center);
		\draw [-stealth] (81.center) to (82.center);
		\draw [-stealth] (78.center) to (79.center);
		\draw [-stealth] (69.center) to (70.center);
		\draw [-stealth] (71.center) to (72.center);
		\draw [-stealth] (76.center) to (77.center);
		\draw (74.center) to (75.center) to (80.center) to (73.center);	\draw [dotted, black!45] (83.center) to (84.center);
	\end{pgfonlayer}
\end{tikzpicture}

%% file: biblio.bib
@inproceedings{morancey2015approximate,
  title={Approximate controllability for a 2D Grushin equation with potential having an internal singularity},
  author={Morancey, Morgan},
  booktitle={Annales de l'Institut Fourier},
  volume={65},
  number={4},
  pages={1525--1556},
  year={2015}
}

@article{allonsius2021analysis,
  title={Analysis of the null controllability of degenerate parabolic systems of Grushin type via the moments method},
  author={Allonsius, Damien and Boyer, Franck and Morancey, Morgan},
  journal={Journal of Evolution Equations},
  volume={21},
  number={4},
  pages={4799--4843},
  year={2021},
  publisher={Springer}
}

@book{watson1922treatise,
  title={A treatise on the theory of Bessel functions},
  author={Watson, George Neville},
  volume={3},
  year={1922},
  publisher={The University Press}
}

@article{LaurentLeauthaud23,
     author = {Camille Laurent and Matthieu L\'eautaud},
     title = {On uniform controllability of {1D} transport equations in the vanishing viscosity limit},
     journal = {Comptes Rendus. Math\'ematique},
     pages = {265--312},
     publisher = {Acad\'emie des sciences, Paris},
     volume = {361},
     year = {2023},
     doi = {10.5802/crmath.405},
     language = {en},
}

@book{coron2007control,
  title={Control and nonlinearity},
  author={Coron, Jean-Michel},
  number={136},
  year={2007},
  publisher={American Mathematical Soc.}
}

@inproceedings{beauchard2020minimal,
  title={Minimal time issues for the observability of Grushin-type equations},
  author={Beauchard, Karine and Dard{\'e}, J{\'e}r{\'e}mi and Ervedoza, Sylvain},
  booktitle={Annales de l'Institut Fourier},
  volume={70},
  number={1},
  pages={247--312},
  year={2020}
}

@article{romankummer2025,
  title={Spectral asymptotics for a class of singular Sturm--Liouville operators with applications to magnetic Laplacian and a-zeros of Kummer functions},
  author={Vanlaere, Roman},
  journal={Journal of Differential Equations},
  volume={475},
  pages={114450},
  year={2026},
  publisher={Elsevier}
}

@article{martinez2019cost,
  title={The cost of boundary controllability for a parabolic equation with inverse square potential},
  author={Martinez, Patrick and Vancostenoble, Judith},
  journal={Evolution Equations and Control Theory (EECT)},
  year={2019}
}

@article{cannarsa2014null,
  title={Null controllability in large time for the parabolic Grushin operator with singular potential},
  author={Cannarsa, Piermarco and Guglielmi, Roberto},
  journal={Geometric control theory and sub-Riemannian geometry},
  pages={87--102},
  year={2014},
  publisher={Springer}
}

@article{beauchard20152d,
  title={2D Grushin-type equations: minimal time and null controllable data},
  author={Beauchard, Karine and Miller, Luc and Morancey, Morgan},
  journal={Journal of Differential Equations},
  volume={259},
  number={11},
  pages={5813--5845},
  year={2015},
  publisher={Elsevier}
}

@article{vanlaere2025grushinlike,
      title={Non Null-Controllability Properties of the Grushin-Like Heat Equation on 2D-Manifolds},
        author={Vanlaere, Roman},
        journal={arXiv preprint arXiv:2503.00997},
        year={2025},
        note={Accepted for publication in ESAIM: Control,
        Optimisation and Calculus of Variations}
}

@article{beauchard2014null,
  title={Null controllability of Grushin-type operators in dimension two.},
  author={Beauchard, Karine and Cannarsa, Piermarco and Guglielmi, Roberto},
  journal={Journal of the European Mathematical Society (EMS Publishing)},
  volume={16},
  number={1},
  year={2014}
}

@article{duprez2020control,
  title={Control of the Grushin equation: non-rectangular control region and minimal time},
  author={Duprez, Michel and Koenig, Armand},
  journal={ESAIM: Control, Optimisation and Calculus of Variations},
  volume={26},
  pages={3},
  year={2020},
  publisher={EDP Sciences}
}

@article{darde2023null,
  title={Null-controllability properties of the generalized two-dimensional Baouendi-Grushin equation with non-rectangular control sets},
  author={Dard{\'e}, J{\'e}r{\'e}mi and Koenig, Armand and Royer, Julien},
  journal={Annales Henri Lebesgue},
  volume={6},
  pages={1479--1522},
  year={2023}
}

@article{CRMATH_2017__355_12_1215_0,
     author = {Koenig, Armand},
     title = {Non-null-controllability of the {Grushin} operator in {2D}},
     journal = {Comptes Rendus. Math\'ematique},
     pages = {1215--1235},
     publisher = {Elsevier},
     volume = {355},
     number = {12},
     year = {2017},
     doi = {10.1016/j.crma.2017.10.021},
     language = {en},
     url = {http://www.numdam.org/articles/10.1016/j.crma.2017.10.021/}
}

@article{LISSYprolate2025,
title = {A non-controllability result for the half-heat equation on the whole line based on the prolate spheroidal wave functions and its application to the Grushin equation},
journal = {Journal of Differential Equations},
volume = {433},
pages = {113306},
year = {2025},
issn = {0022-0396},
doi = {https://doi.org/10.1016/j.jde.2025.113306},
url = {https://www.sciencedirect.com/science/article/pii/S002203962500333X},
author = {Pierre Lissy}
}

@article{anh2016null,
  title={Null controllability in large time of a parabolic equation involving the Grushin operator with an inverse-square potential},
  author={Anh, Cung The and Toi, Vu Manh},
  journal={Nonlinear Differential Equations and Applications NoDEA},
  volume={23},
  number={2},
  pages={20},
  year={2016},
  publisher={Springer}
}

@article{beauchard2014inverse,
  title={Inverse source problem and null controllability for multidimensional parabolic operators of Grushin type},
  author={Beauchard, Karine and Cannarsa, Piermarco and Yamamoto, Masahiro},
  journal={Inverse Problems},
  volume={30},
  number={2},
  pages={025006},
  year={2014},
  publisher={IOP Publishing}
}

@article{tamekue2022null,
  title={Null controllability of the parabolic spherical Grushin equation},
  author={Tamekue, Cyprien},
  journal={ESAIM: Control, Optimisation and Calculus of Variations},
  volume={28},
  pages={70},
  year={2022},
  publisher={EDP Sciences}
}

@article{vancostenoble2008null,
  title={Null controllability for the heat equation with singular inverse-square potentials},
  author={Vancostenoble, Judith and Zuazua, Enrique},
  journal={Journal of Functional Analysis},
  volume={254},
  number={7},
  pages={1864--1902},
  year={2008},
  publisher={Elsevier}
}

@article{ervedoza2008control,
  title={Control and stabilization properties for a singular heat equation with an inverse-square potential},
  author={Ervedoza, Sylvain},
  journal={Communications in Partial Differential Equations},
  volume={33},
  number={11},
  pages={1996--2019},
  year={2008},
  publisher={Taylor \& Francis}
}

@article{beauchard2017heat,
  title={Heat equation on the Heisenberg group: observability and applications},
  author={Beauchard, Karine and Cannarsa, Piermarco},
  journal={Journal of Differential Equations},
  volume={262},
  number={8},
  pages={4475--4521},
  year={2017},
  publisher={Elsevier}
}

@article{LEROUSSEAU20163193,
title = {Null-controllability of the Kolmogorov equation in the whole phase space},
journal = {Journal of Differential Equations},
volume = {260},
number = {4},
pages = {3193-3233},
year = {2016},
issn = {0022-0396},
doi = {https://doi.org/10.1016/j.jde.2015.09.062},
url = {https://www.sciencedirect.com/science/article/pii/S0022039615005367},
author = {Jérôme {Le Rousseau} and Iván Moyano},
}

@article{beauchardKolmog2015,
	author = {{Beauchard, Karine} and {Helffer, Bernard} and {Henry, Raphael} and {Robbiano, Luc}},
	title = {Degenerate parabolic operators of Kolmogorov type with a
          geometric control condition},
	DOI= "10.1051/cocv/2014035",
	url= "https://doi.org/10.1051/cocv/2014035",
	journal = {ESAIM: COCV},
	year = 2015,
	volume = 21,
	number = 2,
	pages = "487-512",
	month = "",
}

@article{koenigfractionalheat,
author = {Koenig, Armand},
title = {Lack of Null-Controllability for the Fractional Heat Equation and Related Equations},
journal = {SIAM Journal on Control and Optimization},
volume = {58},
number = {6},
pages = {3130-3160},
year = {2020},
doi = {10.1137/19M1256610},

URL = { https://doi.org/10.1137/19M1256610},
eprint = {},
}

@article{beauchardkolmogorov,
  title={Null controllability of Kolmogorov-type equations},
  author={Beauchard, Karine},
  journal={Mathematics of Control, Signals, and Systems},
  volume={26},
  number={1},
  pages={145--176},
  year={2014},
  publisher={Springer}
}

@article{DardeRoyerKolmogorov,
     author = {Jeremi Darde and Julien Royer},
     title = {Critical time for the observability of {Kolmogorov-type} equations},
     journal = {Journal de l'Ecole polytechnique Mathematiques},
     pages = {859--894},
     publisher = {Ecole polytechnique},
     volume = {8},
     year = {2021},
     doi = {10.5802/jep.160},
     language = {en},
     url = {https://jep.centre-mersenne.org/articles/10.5802/jep.160/}
}

@incollection{berger2006spectre,
  title={Le spectre d'une vari{\'e}t{\'e} riemannienne},
  author={Berger, Marcel and Gauduchon, Paul and Mazet, Edmond},
  booktitle={Le spectre d’une vari{\'e}t{\'e} Riemannienne},
  pages={141--241},
  year={2006},
  publisher={Springer}
}

@article{cazacu2014controllability,
  title={Controllability of the heat equation with an inverse-square potential localized on the boundary},
  author={Cazacu, Cristian},
  journal={SIAM Journal on Control and Optimization},
  volume={52},
  number={4},
  pages={2055--2089},
  year={2014},
  publisher={SIAM}
}

@article{DardeTrabut2026,
  title={On the critical time of observability of the multi-dimensional Baouendi-Grushin equation},
  author={Dard{\'e}, J{\'e}r{\'e}mi and Trabut, Mathilda},
  journal={arXiv preprint arXiv:2604.01419},
  year={2026}
}

@article{miller2010direct,
  title={A direct Lebeau-Robbiano strategy for the observability of heat-like semigroups},
  author={Miller, Luc},
  journal={Discrete and Continuous Dynamical Systems-Series B},
  volume={14},
  number={4},
  pages={1465--1485},
  year={2010}
}

@article{jerison1999nodal,
  title={Nodal Sets of Sums of Eigenfunctions},
  author={Jerison, David and Lebeau, Gilles},
  journal={Harmonic Analysis and Partial Differential Equations: Essays in Honor of Alberto P. Calderon},
  pages={223},
  year={1999},
  publisher={University of Chicago Press}
}

@book{hardy1952inequalities,
  title={Inequalities},
  author={Hardy, Godfrey Harold and Littlewood, John Edensor and P{\'o}lya, George},
  year={1952},
  publisher={Cambridge university press}
}

@article{dietze2024concentration,
  title={Concentration of eigenfunctions on singular Riemannian manifolds},
  author={Dietze, Charlotte and Read, Larry},
  journal={arXiv preprint arXiv:2410.20563},
  year={2024}
}

@article{dietze2025critical,
  title={The critical case for the concentration of eigenfunctions on singular Riemannian manifolds},
  author={Dietze, Charlotte},
  journal={arXiv preprint arXiv:2510.23520},
  year={2025}
}

@article{helffer2026spectralasymptoticssingularharmonic,
      title={Spectral Asymptotics of a Singular Harmonic Oscillator with Aharonov--Bohm Flux}, 
      author={Bernard Helffer and Ayman Kachmar and François Nicoleau},
      journal={arXiv preprint arXiv:2609.28076v1},
      year={2026},
      primaryClass={math.SP}, 
}

@article{schrodinger1940method,
  title={A method of determining quantum-mechanical eigenvalues and eigenfunctions},
  author={Schr{\"o}dinger, Erwin},
  journal={Proceedings of the Royal Irish Academy. Section A: Mathematical and Physical Sciences},
  volume={46},
  pages={9--16},
  year={1940},
  organization={JSTOR}
}

@inproceedings{cabrera2008beyond,
  title={Beyond conventional factorization: Non-Hermitian Hamiltonians with radial oscillator spectrum},
  author={Cabrera-Munguia, Ivan and Rosas-Ortiz, Oscar},
  booktitle={Journal of Physics: Conference Series},
  volume={128},
  number={1},
  pages={012042},
  year={2008}
}

@article{quesne2009solvable,
  title={Solvable rational potentials and exceptional orthogonal polynomials in supersymmetric quantum mechanics},
  author={Quesne, Christiane and others},
  journal={SIGMA. Symmetry, Integrability and Geometry: Methods and Applications},
  volume={5},
  pages={084},
  year={2009},
  publisher={SIGMA. Symmetry, Integrability and Geometry: Methods and Applications}
}

@article{tamekue2026minimal,
  title={Minimal time for null controllability of the parabolic spherical Baouendi-Grushin equation},
  author={Tamekue, Cyprien},
  journal={arXiv preprint arXiv:2604.19700},
  year={2026}
}

@unpublished{boyer:hal-02470625,
  TITLE = {{Controllability of linear parabolic equations and systems}},
  AUTHOR = {Boyer, Franck},
  URL = {https://hal.science/hal-02470625},
  NOTE = {Lecture},
  TYPE = {Master},
  ADDRESS = {France},
  YEAR = {2022},
  MONTH = Feb,
  HAL_ID = {hal-02470625},
  HAL_VERSION = {v4},
}

@article{le2012carleman,
  title={On Carleman estimates for elliptic and parabolic operators. Applications to unique continuation and control of parabolic equations},
  author={Le Rousseau, J{\'e}r{\^o}me and Lebeau, Gilles},
  journal={ESAIM: Control, Optimisation and Calculus of Variations},
  volume={18},
  number={3},
  pages={712--747},
  year={2012}
}

@article{vanlaere2026observability,
  title={Observability properties of the singular Grushin equation},
  author={Vanlaere, Roman},
  journal={arXiv preprint arXiv:2602.08044},
  Note={To appear in Advances in Differential Equations, 32(3--4), 2027.},
  year={2026}
}

@article{boulard2026minimal,
  title={Minimal time for heat observability on almost-Riemannian surfaces},
  author={Boulard, Vincent and Hayat, Amaury and Tr{\'e}lat, Emmanuel},
  journal={arXiv preprint arXiv:2610.16568},
  year={2026}
}

@article{letrouit2023subelliptic,
  title={Subelliptic wave equations are never observable},
  author={Letrouit, Cyril},
  journal={Analysis \& PDE},
  volume={16},
  number={3},
  pages={643--678},
  year={2023},
  publisher={Mathematical Sciences Publishers}
}

@article{burq2022time,
  title={Time optimal observability for Grushin Schr{\"o}dinger equation},
  author={Burq, Nicolas and Sun, Chenmin},
  journal={Analysis \& PDE},
  volume={15},
  number={6},
  pages={1487--1530},
  year={2022},
  publisher={Mathematical Sciences Publishers}
}

@article{allibert1998controle,
  title={Controle analytique de l'equation des ondes et de l'equation de schrodinger sur des surfaces de revolution},
  author={Allibert, Brice},
  journal={Communications in partial differential equations},
  volume={23},
  number={9-10},
  pages={1493--1556},
  year={1998},
  publisher={Taylor \& Francis}
}

@article{laurent2023uniform,
  title={Uniform observation of semiclassical Schr{\"o}dinger eigenfunctions on an interval},
  author={Laurent, Camille and L{\'e}autaud, Matthieu},
  journal={Tunisian Journal of Mathematics},
  volume={5},
  number={1},
  pages={125--170},
  year={2023},
  publisher={Mathematical Sciences Publishers}
}

@book{laurent2022tunneling,
  title={Tunneling estimates and approximate controllability for hypoelliptic equations},
  author={Laurent, Camille and L{\'e}autaud, Matthieu},
  volume={276},
  number={1357},
  year={2022},
  publisher={American Mathematical Society}
}

@book{abramowitz1948handbook,
  title={Handbook of mathematical functions with formulas, graphs, and mathematical tables},
  author={Abramowitz, Milton and Stegun, Irene A},
  volume={55},
  year={1948},
  publisher={US Government printing office}
}

@article{lebeau1995controle,
  title={Contr{\^o}le exact de l{\'e}quation de la chaleur},
  author={Lebeau, Gilles and Robbiano, Luc},
  journal={Communications in Partial Differential Equations},
  volume={20},
  number={1-2},
  pages={335--356},
  year={1995},
  publisher={Taylor \& Francis}
}

@article{ervedoza2011observability,
    author  = {Ervedoza, Sylvain and Zuazua, Enrique},
    title   = {Observability of heat processes by transmutation without geometric restrictions},
    journal = {Mathematical Control and Related Fields},
    year    = {2011},
    volume  = {1},
    number  = {2},
    pages   = {177--187},
    doi     = {10.3934/mcrf.2011.1.177}
}

@article{ervedoza2011sharp,
  title={Sharp observability estimates for heat equations},
  author={Ervedoza, Sylvain and Zuazua, Enrique},
  journal={Archive for rational mechanics and analysis},
  volume={202},
  number={3},
  pages={975--1017},
  year={2011},
  publisher={Springer}
}

@article{vancostenoble2009hardy,
  title={Hardy inequalities, observability, and control for the wave and Schr{\"o}dinger equations with singular potentials},
  author={Vancostenoble, Judith and Zuazua, Enrique},
  journal={SIAM journal on mathematical analysis},
  volume={41},
  number={4},
  pages={1508--1532},
  year={2009},
  publisher={SIAM}
}

@article{huang2026observability,
  title={Observability and stabilizability for Grushin-Type parabolic equations in the whole space},
  author={Huang, Shanlin and Wang, Gengsheng and Yu, Jiaqi},
  journal={arXiv preprint arXiv:2611.05821},
  year={2026}
}

@article{alphonse2024quantitative,
  title={Quantitative spectral inequalities for the anisotropic Shubin operators and applications to null-controllability},
  author={Alphonse, Paul and Seelmann, Albrecht},
  journal={Comptes Rendus. Math{\'e}matique},
  volume={362},
  number={G12},
  pages={1635--1659},
  year={2024}
}

@article{alphonse2025unique,
  title={Unique continuation estimates for Baouendi--Grushin equations on cylinders},
  author={Alphonse, Paul and Seelmann, Albrecht},
  journal={Pure and Applied Analysis},
  volume={7},
  number={3},
  pages={733--770},
  year={2025},
  publisher={Mathematical Sciences Publishers}
}

@article{jaming2025null,
  title={Null-controllability of the generalized Baouendi--Grushin heat like equations},
  author={Jaming, Philippe and Wang, Yunlei},
  journal={Journal of Evolution Equations},
  volume={25},
  number={2},
  pages={43},
  year={2025},
  publisher={Springer}
}

@article{infeld1951factorization,
  title={The factorization method},
  author={Infeld, Leopold and Hull, TE},
  journal={Reviews of modern Physics},
  volume={23},
  number={1},
  pages={21},
  year={1951},
  publisher={APS}
}
